\documentclass[12pt]{article}
\usepackage{latexsym,amssymb,amsmath}
\begin{document}
\baselineskip=20pt

\newcommand{\la}{\langle}
\newcommand{\ra}{\rangle}
\newcommand{\psp}{\vspace{0.4cm}}
\newcommand{\pse}{\vspace{0.2cm}}
\newcommand{\ptl}{\partial}
\newcommand{\dlt}{\delta}
\newcommand{\sgm}{\sigma}
\newcommand{\al}{\alpha}
\newcommand{\be}{\beta}
\newcommand{\G}{\Gamma}
\newcommand{\gm}{\gamma}
\newcommand{\vs}{\varsigma}
\newcommand{\Lmd}{\Lambda}
\newcommand{\lmd}{\lambda}
\newcommand{\td}{\tilde}
\newcommand{\vf}{\varphi}
\newcommand{\yt}{Y^{\nu}}
\newcommand{\wt}{\mbox{wt}\:}
\newcommand{\rd}{\mbox{Res}}
\newcommand{\ad}{\mbox{ad}}
\newcommand{\stl}{\stackrel}
\newcommand{\ol}{\overline}
\newcommand{\ul}{\underline}
\newcommand{\es}{\epsilon}
\newcommand{\dmd}{\diamond}
\newcommand{\clt}{\clubsuit}
\newcommand{\vt}{\vartheta}
\newcommand{\ves}{\varepsilon}
\newcommand{\dg}{\dagger}
\newcommand{\tr}{\mbox{Tr}}
\newcommand{\ga}{{\cal G}({\cal A})}
\newcommand{\hga}{\hat{\cal G}({\cal A})}
\newcommand{\Edo}{\mbox{End}\:}
\newcommand{\for}{\mbox{for}}
\newcommand{\kn}{\mbox{ker}}
\newcommand{\Dlt}{\Delta}
\newcommand{\rad}{\mbox{Rad}}
\newcommand{\rta}{\rightarrow}
\newcommand{\mbb}{\mathbb}
\newcommand{\lra}{\Longrightarrow}

\begin{center}{\Large \bf Etingof Trace,  Path Hypergeometric}\end{center}
\begin{center}{\Large \bf Functions and Integrable Systems}\footnote{Research
supported by CNSF Grant 10371121}\end{center}

\vspace{0.2cm}

\begin{center}{\large Xiaoping Xu}\end{center}
\begin{center}{Institute of Mathematics, Academy of Mathematics \& System Sciences}\end{center}
\begin{center}{Chinese Academy of Sciences, Beijing 10080, P.R.China}\end{center}
\vspace{1cm}

\begin{center}{\Large \bf Abstract}\end{center}

{\small Under a certain condition, we find the explicit formulas
for the trace functions of certain intertwining operators among
$gl(n)$-modules, introduced by Etingof in connection with the
solutions of the Calogero-Sutherland model. If $n=2$, the master
function of the trace function is exactly the classical Gauss
hypergeometric function. When $n>2$, the master functions of the
trace functions are a new family of multiple hypergeometric
functions, whose differential property and integral representation
are  dominated by certain polynomials of integral paths connecting
pairs of positive integers. Moreover, we define and explicitly
find similar trace functions for $sp(2n)$, which give rise to
solutions of the Olshanesky-Perelomov model of type C. The master
functions of the trace functions for $sp(2n)$ are similar new
multiple path hypergeometric functions. Analogous multiple path
hypergeometric functions for orthogonal Lie algebras are defined
and studied}.

\section{Introduction}

The Calogero-Sutherland model is an exactly solvable quantum many-body system in one-dimension
 (cf. [C], [S]) , whose Hamiltonian is given by
$$H_{CS}=\sum_{i=1}^n\ptl_{x_i}^2+K\sum_{1\leq i<j\leq n}\frac{1}{\sinh^2(x_i-x_j)},\eqno(1.1)$$
where $K$ is a constant. The model was used to study long-range
interactions of $n$ particles. It has become an important
mathematical object partly because its potential is related to the
root
  system of the special linear algebra $sl(n)$. Olshanesky and Perelomov [OP] generalized the
  Calogero-Sutherland model to the system whose Hamiltonian is given by
\begin{eqnarray*}\qquad H_{OP}&=&\sum_{i=1}^n\ptl_{x_i}^2+\sum_{1\leq i<j\leq n}\left(\frac{K_1}
{\sinh^2(x_i-x_j)}+\frac{K_2}{\sinh^2(x_i+x_j)}\right)\\ &
&+\sum_{i=1}^n\left(\frac{K_3}
{\sinh^2x_i}+\frac{K_4}{\sinh^22x_i}\right),\hspace{7cm}(1.2)\end{eqnarray*}
where $K_1,K_2,K_3,K_4$ are constants. The potential of the
Olshanesky- Perelomov is related to the root systems of all four
families of finite-dimensional classical simple Lie algebras.
Solving a system is to find eigenfunctions of its Hamiltonians.
Etingof [Ep] found that the trace functions of certain
intertwining operators among $gl(n)$-modules giving rise to
eigenfunctions of $H_{CS}$.

    The aim of this paper is to find the functional implication of
   the Etingof trace.  In the case of $gl(2)$, the master function of the Etingof trace
function is exactly the classical Gauss hypergeometric function.
When $n>2$,  under a certain condition on weights, we find that
the master functions of the Etingof trace functions are a new
family of multiple hypergeometric functions, whose differential
property and integral representation are dominated by certain
polynomials of integral paths connecting pairs of positive
integers. Moreover, we define and explicitly find similar trace
functions for $sp(2n)$, which give rise to solutions of the
Olshanesky-Perelomov model (1.2) of type C ($K_3=0$). The master
functions of the trace functions for $sp(2n)$ are similar new
multiple path hypergeometric functions. Based on the combinatorial
feature of the path hypergeometric functions for $gl(n)$ and
$sp(2n)$, we define and study analogous hypergeometric functions
for orthogonal Lie algebras. Furthermore, we prove  that certain
variations of the Weyl functions of classical simple Lie algebras
naturally give solutions of  the Olshanesky-Perelomov model with
different sets of constants.

Below we give more detailed outline of the background and the results of this paper.

    The term of ``hypergeometric'' was first used by Wallis in Oxford as early as 1655 in his
    work {\it Arithmetrica Infinitorm} when referring  to any series which could be regarded as
     a generalization of the ordinary geometric series
$\sum_{n=0}^{\infty}z^n$. The well-known hypergeometric series of one-variable
$_2F_1(a,b;c;z)$ was introduced and studied  by Gauss in his thesis
presented at G\"{o}ttingen in 1812. Apell (1880) first systematically studied double
hypergeometric functions and defined his famous four functions of types $F_1$-$F_4$. Horn more
 completely  investigated  all the double hypergeometric functions of second order, in 1889,
 1931 and 1938, by studying corresponding second
order hypergeometric partial differential systems with two independent variables. He found
additional ten double hypergeometric functions of types $G_1$-$G_3$ and $H_1$-$H_7$, known as
 the Horn functions. Lauricella (1893) generalized Apell's four double hypergeometric functions
 to four families of  multiple  hypergeometric functions of types $F_A^{(n)},\;F_B^{(n)},
 \;F_C^{(n)}$ and $F_D^{(n)}$. Moreover, he conjectured the existence of fourteen complete
  triple hypergeometric functions of second order. Saran (1954) found the other ten triple
   hypergeometric functions of types $F_E$-$F_T$. However, Srivastava (1964) discovered three
more complete triple hypergeometric functions of second order of
types $H_A,\;H_B$ and $H_C$. Saran (1955) gave integrals of
Pochhammer type of the functions of types  $F_E$-$F_T$. Pandey
(1965) transformed the integral representations of $F_F$ and $F_G$
such that two interesting new triple hypergeometric functions of
the generalized Horn type were obtained. Exton (1972, 1973)
defined and studied twenty-six quadruple  hypergeometric functions
of types $K_1$-$K_{21}$ and $D_1$-$D_5$, which are generalizations
of the Horn functions. Furthermore, he introduced two families of
multiple hypergeometric functions in the line of his quadruple
functions. The multiple hypergeometric functions we find this this
paper are in the stream of the above hypergeometric functions of
classical type. We refer [Eh] for the details of the
hypergeometric functions of classical type and the references
therein.

 Throughout this paper, all the vector spaces (algebras) are assumed over $\mbb{C}$, the field
 of complex numbers. Denote by $\mbb{Z}$ the ring of integers and by $\mbb{N}$ the additive
 semigroup of nonnegative integers. For two integers $m_1,m_2$, we use the notation:
$$\ol{m_1,m_2}=\left\{\begin{array}{ll}\{m_1,m_1+1,...,m_2\}&\mbox{if}\;\;m_1\leq m_2\\
\emptyset & \mbox{if}\;\;m_1> m_2.\end{array}\right.\eqno(1.3)$$
The general Lie algebra
$$gl(n)=\mbox{Span}\{E_{i,j}\mid i,j\in\ol{1,n}\}\eqno(1.4)$$
with the Lie bracket
$$[A,B]=AB-BA\qquad\for\;\;A,B\in gl(n),\eqno(1.5)$$
where $E_{i,j}$ is the $n\times n$ matrix whose $(i,j)$-entry is 1 and the others are 0. For
$\mu\in\mbb{C}$, we denote
$$V_\mu=\mbox{Span}\{(x_1x_2\cdots x_n)^\mu x_1^{i_1}x_2^{i_2}\cdots x_n^{i_n}
\mid i_1,i_2,...,i_n\in\mbb{Z}\},\eqno(1.6)$$
where $x_1,x_2,..., x_n$ are $n$ indeterminates. Then $V_\mu$
becomes a $gl(n)$-module with the action
$$E_{i,j}=x_i\ptl_{x_j}-\mu\dlt_{i,j}\qquad \for\;\;i,j\in\ol{1,n}.\eqno(1.7)$$
The representation of $gl(n)$ in (1.7) is called a {\it point representation}.

The subspace
$$T=\sum_{i=1}^n\mbb{C}E_{i,i}\eqno(1.8)$$
is a toral Cartan subalgebra of the Lie algebra $gl(n)$. A {\it weight} $\lmd$ of $gl(n)$ is a
linear function on $T$ determined by
$$\lmd(E_{i,i})=\lmd_i\qquad \for\;\;i\in\ol{1,n}.\eqno(1.9)$$
Let $M_{\lmd}$ be an irreducible Verma module of $gl(n)$ with the highest weight $\lmd$ (e.g.
 cf. [Hj]) and let $M_\lmd \hat{\otimes}V_\mu$ be the completion of the tensor $M_\lmd \otimes
  V_\mu$ with respect to the topology given by weight subspaces. Suppose that $\Phi: M_\lmd\rta
   M_\lmd \hat{\otimes}V_\mu $ is a $gl(n)$-module homomorphism. View $\Phi$ as a function in
$\{x_1,x_2,...,x_n\}$ taking values in $\mbox{End}\:M_\lmd$, the
space of linear transformations on $M_\lmd$. The {\it Etingof
trace function}
$$E(z_1,z_2,...,z_n)=(x_1x_2\cdots x_n)^{-\mu}\mbox{tr}_{_{M_\lmd}}\:\Phi z_1^{E_{1,1}}
z_2^{E_{2,2}}\cdots z_n^{E_{n,n}}\eqno(1.10)$$ (cf. [Ep]).
Changing variables
$$z_i=e^{2x_i}\qquad\for\;\;i\in\ol{1,n},\eqno(1.11)$$
we have the equation for the Calogero-Sutherland model (1.1):
$$\sum_{i=1}^n(z_i\ptl_{z_i})^2(\Psi)+K\left(\sum_{1\leq i<j\leq n}\frac{z_iz_j}{(z_i-z_j)^2}
\right)\Psi=\nu \Psi\eqno(1.12)$$
and  the equation for the Olshanesky-Perelomov model (1.2):
\begin{eqnarray*}\qquad& & \sum_{i=1}^n\ptl_{z_i}^2(\Psi)+[\sum_{1\leq i<j\leq n}
\left(K_1\frac{z_iz_j}{(z_i-z_j)^2}+K_2\frac{z_iz_j}{(z_iz_j-1)^2}\right)\\ & &+\sum_{i=1}^n
\left(K_3\frac{z_i}{(z_i-1)^2}+K_4\frac{z_i^2}{(z_i^2-1)^2}\right)]\Psi=\nu\Psi,\hspace{5.2cm}
(1.13)\end{eqnarray*}
Etingof [Ep] proved that the function
$$\Psi(z_1,z_2,...,z_n)=\frac{E(z_1,z_2,...,z_n)}{\mbox{tr}_{_{M_{-\rho}}}z_1^{E_{1,1}}
z_2^{E_{2,2}}\cdots z_n^{E_{n,n}}}\eqno(1.14)$$
is a solution of (1.12) for suitable constants $K$ and $\nu$, where
$$\rho(E_{i,i})=\frac{n+1}{2}-i\qquad\for\;\;i\in\ol{1,n},\eqno(1.15)$$

Denote
$$(c)_i=c(c+1)(c+2)\cdots (c+i-1)\qquad \for\;\;c\in\mbb{C},\;i\in\mbb{N}.\eqno(1.16)$$
The  classical Gauss hypergeometric function
$$_2F_1(a,b;c; z)=\sum_{m=0}^{\infty}\frac{(a)_m(b)_m}{m!(c)_m}z^m,\qquad\mbox{where}\;\;a,b,c
\in\mbb{C},\;-c\not\in\mbb{N}.\eqno(1.17)$$
In the case $n=2$, we find that
$$E(z_1,z_2)=\frac{z_1^{\lmd_1+1}z_2^{\lmd_2}}{z_1-z_2}\:_2F_1\left(\mu+1,-\mu;\lmd_2-\lmd_1,
\frac{z_2}{z_2-z_1}\right)\eqno(1.18)$$
up to a scalar multiple. Thus we can view $_2F_1(\tau_1,\tau_2;\sgm; z)$ as the master function
 of $E(z_1,z_2)$.

Let
$$\G_A=\sum_{1\leq j<i\leq n}\mbb{N}\es_{i,j}\eqno(1.19)$$
be the torsion-free additive semigroup  of  rank $n(n-1)/2$ with $\es_{i,j}$ as base elements.
 For $\al=\sum_{1\leq j<i\leq n}\al_{i,j}\es_{i,j}\in \G_A$, we denote
$$\al_{_{\ul{1}}}=\al_{\ol{n}}=0,\;\;\al_{_{\ul{i}}}=\sum_{r=1}^{i-1}\al_{i,r},\;\;
\al_{_{\ol{j}}}=\sum_{s=j+1}^n\al_{s,j}\eqno(1.20)$$
 Moreover, we define our $(n(n-1)/2)$-variable hypergeometric function of type A  by
$${\cal X}_A(\tau_1,..,\tau_n;\vt)\{z_{j,k}\}=\sum_{\be\in\G_A}\frac{\left[\prod_{s=1}^{n-1}
(\tau_s-\be_{\ul{s}})_{_{\be_{\ol{s}}}}\right](\tau_n)_{_{\be_{\ul{n}}}}}
{\be!(\vt)_{_{\be_{\ul{n}}}}}z^\be,\eqno(1.21)$$
where
$$\be!=\prod_{1\leq k<j\leq n}\be_{j,k}!,\qquad z^\be=\prod_{1\leq k<j\leq n}
z_{j,k}^{\be_{j,k}}.\eqno(1.22)$$
 Set
$$\xi_{r_2,r_1}^A=\prod_{s=r_1}^{r_2-1}\frac{z_{s}}{z_{r_2}-z_s}\qquad
 \for\;\;1\leq r_1<r_2\leq n.\eqno(1.23)$$
When $n>2$, we assume
$$\mu=\lmd_i-\lmd_{i+1}\qquad\for\;\;i\in\ol{1,n-2}.\eqno(1.24)$$
 Under the condition (1.24), we find the
  Etingof trace function
$$E(z_1,z_2,...,z_n)=\frac{\prod_{r=1}^nz_r^{\lmd_r+n-r}}{\prod_{1\leq k<j\leq n}(z_k-z_j)}
{\cal X}_A(\mu+1,..,\mu+1,-\mu;\lmd_n-\mu)\{\xi_{r_2,r_1}^A\}.\eqno(1.25)$$
The  differential property and integral representation of ${\cal X}_A$ are dominated by certain
 polynomials of integral paths connecting pairs of positive integers (cf. Theorems 5.1 and 5.2).
  Moreover, the system of partial differential
equations for ${\cal X}_A$ is a natural multi-variable analogue of the classical hypergeometric
equation (e.g. cf. [AAR]).

Fernando [F] proved that point representations of finite-dimensional simple Lie algebras exist
only for $sl(n)$ and $sp(2n)$. Using the similar point representation of $sp(2n)$ as that of
 $gl(n)$ in (1.7),   we define and find analogues of the Etingof trace functions for $sp(2n)$,
 which give rise to solutions of the equation (1.13) with suitable constants $K_1,K_2,K_4,\nu$
 and $K_3=0$. The master functions of the trace functions for $sp(2n)$ are similar to the
 function ${\cal X}_A$ in (1.21).  Furthermore, we prove  that certain variations of the Weyl
 function of classical simple Lie algebras naturally give solutions of the equation (1.13) with
  different sets of constants. The results may yield more solutions of the equation (1.13)
  related to the solutions of the systems of partial differential equations for our path
  hypergeometric functions, as we show for the case $gl(2)$ (cf. Theorem 2.2).

Our path hypergeometric functions resemble the following
well-known multiple hypergeometric
 functions:
Apell's first double hypergeometric function (1880)
$$F_1(a_1,a_2,a_3:c;z_1,z_2)=\sum_{m_1,m_2=0}^{\infty}\frac{(a_1)_{m_1+m_2}(a_2)_{m_1}
(a_3)_{m_2}}{m_1!m_2!(c)_{m_1+m_2}}z_1^{m_1}z_2^{m_2},\eqno(1.26)$$
Lauricella's  second family of  hypergeometric functions (1893)
\begin{eqnarray*} \qquad &&F^{(k)}_B(a_1,...,a_k,b_1,...,b_k;c;z_1,...,z_k)
\\ &=&\sum_{m_1,...,m_k=0}^{\infty}\frac{(a_1)_{m_1}\cdots (a_k)_{m_k}(b_1)_{m_1}\cdots
(b_k)_{m_k}}{(c)_{m_1+\cdots m_k}}\frac{z_1^{m-1}\cdots z_k^{m_k}}{m_1!\cdots m_k!}
\hspace{3.4cm}(1.27)\end{eqnarray*}
and Horn's eighth double hypergeometric functions:
$$H_5(a_1,a_2;c;z_1,z_2)=\sum_{m_1,m_2=0}^{\infty}\frac{(a_1)_{2m_1+m_2}(a_2)_{m_2-m_1}}
{m_1!m_2!(c)_{m_2}}z_1^{m_1}z_2^{m_2}\eqno(1.28)$$
(cf. [Eh] for details). Heckman and Opdam [HO, Hg1,O1, O2, BO] introduced hypergeometric
equations related to root systems and analogous to (1.12). They proved the existence of
 solutions (hypergeometric functions) of their equations.
In our case, the interesting functions are the path hypergeometric functions like ${\cal X}_A$
 in (1.21) and the functions like $\Psi$ in (1.14) are not interesting from pure function
  point of view. Gel'fand and Graev studied analogues of classical hypergeometric functions
  (so called GG-functions) by generalizing the differential property of the classical
  hypergeometric functions (e.g. cf. [GG]).

The paper is organized as follows. In Section 2, we calculate the
Etingof trace function for $gl(2)$ and show a connection between
the equation (1.12) with $n=2$ and the classical hypergeometric
equation for the Gauss hypergeometric function. This section
serves as a contraction of the whole paper. We obtain (1.25) in
Section 3. Section 4 is devoted to defining and
 calculating the analogue of the Etingof trace function for
 $sp(2n)$. In Section 5, we find the differential property and
 systems of partial differential equations for the multiple
 hypergeometric functions of types A and C obtained in Section 3
 and Section 4, as the master functions of the Etingof trace
 functions. The integral representations of our multiple
 hypergeometric functions of type A are obtained. Analogous  multiple
 hypergeometric functions associated with the root systems of the
 orthogonal Lie algebras are introduced and studied.
In Section 6, we show a connection between the equation (1.13) and Weyl functions of classical
 Lie algebras. Section 7 is the proof of that the trace functions of certain intertwining
 operators among $sp(2n)$-modules give rise to solutions of the equation (1.13) with $K_3=0$.

\section{Preliminary}

In this section, we find Etingof trace functions for $gl(2)$ in
terms of classical Gauss hypergeometric functions. Moreover, we
show a connection between the equation (1.12) with
 $n=2$ and the classical hypergeometric equation.

Let $\lmd$ be a weight of $gl(2)$
such that
$$\sgm=\lmd_1-\lmd_2\not\in\mbb{N}\eqno(2.1)$$
(cf. (1.9)). Choose $E_{1,2}$ as  the positive root vector.  The Verma
$gl(2)$-module with the highest-weight vector $v_\lmd$ of weight $\lmd$ is given by
$$M_\lmd =\mbox{Span}\{E_{2,1}^iv_\lmd\mid i\in\mbb{N}\}\eqno(2.2)$$
with the action
$$E_{1,,2}(E_{2,1}^iv_\lmd)=i(\sgm+1-i)E_{2,1}^{i-1}v_\lmd,\qquad E_{2,1}(E_{2,1}^iv_\lmd)
=E_{2,1}^{i+1}v_\lmd,\eqno(2.3)$$
$$E_{1,1}(E_{2,1}^iv_\lmd)=(\lmd_1-i)E_{2,1}^iv_\lmd,\qquad E_{2,2}(E_{2,1}^iv_\lmd)
=(\lmd_2+i)E_{2,1}^iv_\lmd.\eqno(2.4)$$
Under the assumption (2.1), $M_\lmd$ is an irreducible
$gl(2)$-module.

 Let $\mu\in\mbb{C}$ be a fixed constant and let
$$\td{M}=\{\sum_{i_1,i_2\in\mbb{Z}}v_{i_1,i_2}x_1^{i_1+\mu}x_2^{i_2+\mu}\mid v_{i_1,i_2}\in
M_\lmd\},\eqno(2.5)$$
a space of Laurent series with coefficients in $M_\lmd$. Define the action of $gl(2)$ on
$\td{M}$ by
\begin{eqnarray*}E_{j_1,j_2}(\sum_{i_1,i_2\in\mbb{Z}}v_{i_1,i_2}x_1^{i_1+\mu}x_2^{i_2+\mu})
&=&\sum_{i_1,i_2\in\mbb{Z}}[E_{j_1,j_2}(v_{i_1,i_2})x_1^{i_1+\mu}x_2^{i_2+\mu}
\\ & &+v_{i_1,i_2}(x_{j_1}\ptl_{x_{j_2}}-\mu\dlt_{j_1,j_2})(x_1^{i_1+\mu}x_2^{i_2+\mu})]
\hspace{2.7cm}(2.6)\end{eqnarray*}
for $j_1,j_2\in\{1,2\}$. Then $\td{M}$ becomes a $gl(2)$-module that is isomorphic
to $M_\lmd \hat{\otimes}V_\mu$ (cf. (1.6)). Recall that a singular vector in a Verma module is
 a nonzero weight vector annihilated by positive root vectors.
 Any singular vector of weight $\lmd$ in $\td{M}$ is of the form:
$$u_\lmd=(x_1x_2)^\mu\sum_{i=0}^{\infty}a_iE_{2,1}^iv_\lmd x_1^ix_2^{-i}\eqno(2.7)$$
Observe that
\begin{eqnarray*}0&=&E_{1,2}(u_\lmd)\\ &=&(x_1x_2)^\mu\sum_{i=0}^{\infty}a_i[i(\sgm+1-i)
E_{2,1}^{i-1}v_\lmd x_1^ix_2^{-i}+(\mu-i)E_{2,1}^iv_\lmd x_1^{i+1}x_2^{-i-1}]
\\ &=&(x_1x_2)^\mu\sum_{i=0}^{\infty}((\mu-i)a_i+(i+1)(\sgm-i)a_{i+1})v_\lmd
x_1^{i+1}x_2^{-i-1}.\hspace{4.1cm}(2.8)\end{eqnarray*}
Thus
$$(\mu-i)a_i+(i+1)(\sgm-i)a_{i+1}=0\qquad\for\;\;i\in\mbb{N},\eqno(2.9)$$
equivalently,
$$a_{i+1}=-\frac{(\mu-i)a_i}{(i+1)(\sgm-i)}\qquad\for\;\;i\in\mbb{N},\eqno(2.10)$$
For any constant $c$ and nonnegative integer $i$, we denote
$$\la c\ra_i=c(c-1)\cdots (c-(i-1)).\eqno(2.11)$$
We normalize $u_\lmd$ by taking
$$a_0=1.\eqno(2.12)$$
By induction,
$$a_i=\frac{(-1)^i\la \mu\ra_i}{i!\la \sgm\ra_i}.\eqno(2.13)$$
Thus
$$u_\lmd =\sum_{i=0}^{\infty}\frac{(-1)^i\la \mu\ra_i}{i!\la \sgm\ra_i}E_{2,1}^iv_\lmd
x_1^{i+\mu}x_2^{-i+\mu}.\eqno(2.14)$$

Let $\Phi: M_\lmd\rta \td{M}$ be the Lie algebra module homomorphism
such that $\Phi(v_\lmd)=u_\lmd$.
Note that
$$E_{2,1}^j(E_{2,1}^iv_\lmd  x_1^{i+\mu}x_2^{-i+\mu})=\sum_{k=0}^j\left(\!\!\begin{array}{c}j\\
 k\end{array}\!\!\right)\la\mu+i\ra_kE_{2,1}^{i+j-k}v_\lmd x_1^{i-k+\mu}x_2^{-i+k+\mu}.
 \eqno(2.15)$$
Moreover, for any $\iota\in\mbb{N}$, we have
$$\sum_{k=\iota}^{\infty}\left(\!\!\begin{array}{c}k\\ \iota\end{array}\!\!\right)z^k=
z^{\iota}\sum_{k=\iota}^{\infty}\left(\!\!\begin{array}{c}k\\ \iota\end{array}\!\!\right)
z^{k-\iota}=\frac{z^{\iota}}{\iota !}\frac{d^\iota}{dz^\iota}(\sum_{k=0}^{\infty}z^k)
=\frac{z^{\iota}}{\iota !}\frac{d^\iota}{dz^\iota}\left(\frac{1}{1-z}\right)=\frac{z^{\iota}}
{(1-z)^{\iota+1}}.\eqno(2.16)$$
View $\Phi$ as a ``function taking values in $\mbox{End}\:M_\lmd$.''
By (2.15) and (2.16), the Etingof trace function
\begin{eqnarray*}& &E(z_1,z_2)=(x_1x_2)^{-\mu}\mbox{tr}_{_{M_\lmd}}\:\Phi
z_1^{E_{1,1}}z_2^{E_{2,2}}\\
&=&z_1^{\lmd_1}z_2^{\lmd_2}\sum_{j=0}^{\infty}\sum_{i=0}^\infty
\left(\!\!\begin{array}{c}j\\
i\end{array}\!\!\right)\la\mu+i\ra_i\frac{(-1)^i\la
\mu\ra_i}{i!\la
\sgm\ra_i}\left(\frac{z_2}{z_1}\right)^j\hspace{8cm}
\end{eqnarray*}
\begin{eqnarray*} &=&
z_1^{\lmd_1}z_2^{\lmd_2}\sum_{i=0}^\infty\la\mu+i\ra_i\frac{(-1)^i\la
\mu\ra_i}{i!\la \sgm\ra_i}
\sum_{j=0}^{\infty}\left(\!\!\begin{array}{c}j\\
i\end{array}\!\!\right)\left(\frac{z_2}{z_1}
\right)^j\\&=&z_1^{\lmd_1}z_2^{\lmd_2}\sum_{i=0}^\infty\la\mu+i\ra_i\frac{(-1)^i\la
\mu\ra_i}{i!\la
\sgm\ra_i}\left(\frac{z_2}{z_1}\right)^i\left(\frac{1}{1-z_2/z_1}\right)^{i+1}\\
&=&\frac{z_1^{\lmd_1+1}z_2^{\lmd_2}}{z_1-z_2}\sum_{i=0}^\infty\frac{\la\mu+i\ra_i\la
\mu\ra_i}{i!\la
\sgm\ra_i}\left(\frac{z_2}{z_2-z_1}\right)^i.\hspace{7.2cm}(2.17)\end{eqnarray*}

Recall the notion $(c)_n$ in (1.16), which is widely used in  the theory of  hypergeometric
functions. Note
$$ \la c+n\ra_n=(c+1)_n,\;\;\la c\ra_n=(-1)^n(-c)_n\qquad\for\;\;c\in\mbb{C},\;n\in\mbb{N}.
\eqno(2.18)$$
By (1.17), (2.1), (2.17) and (2.18), we have:
\psp

{\bf Theorem 2.1}. {\it The Etingof trace function for $gl(2)$ is}
$$E(z_1,z_2)=\frac{z_1^{\lmd_1+1}z_2^{\lmd_2}}{z_1-z_2}\:_2F_1\left(\mu+1,-\mu;\lmd_2-\lmd_1:
\frac{z_2}{z_2-z_1}\right).\eqno(2.19)$$
\pse

The function in (1.14)
$$\Psi(z_1,z_2)=z_1^{\lmd_1+1/2}z_2^{\lmd_2-1/2}\:_2F_1\left(\mu+1,-\mu;\lmd_2-\lmd_1:
\frac{z_2}{z_2-z_1}\right)\eqno(2.20)$$

Next, we will give a connection between the equation (1.12) with $n=2$ and the classical
hypergeometric equation. For $\mu_1,\mu_2\in\mbb{C}$, we set
$$\phi_{\mu_1,\mu_2}(z_1,z_2)=(z_1z_2)^{\mu_1}(z_1-z_2)^{\mu_2}.\eqno(2.21)$$
Then
$$z_1\ptl_{z_1}(\phi_{\mu_1,\mu_2})=\left(\mu_1+\mu_2\frac{z_1}{z_1-z_2}\right)
\phi_{\mu_1,\mu_2},\eqno(2.22)$$
$$z_2\ptl_{z_2}(\phi_{\mu_1,\mu_2})=\left(\mu_1-\mu_2\frac{z_2}{z_1-z_2}\right)
\phi_{\mu_1,\mu_2},\eqno(2.23)$$
$$(z_1\ptl_{z_1})^2(\phi_{\mu_1,\mu_2})=\left(\mu_1^2+(2\mu_1+1)\mu_2\frac{z_1}
{z_1-z_2}+\mu_2(\mu_2-1)\frac{z_1^2}{(z_1-z_2)^2}\right)\phi_{\mu_1,\mu_2},\eqno(2.24)$$
$$(z_2\ptl_{z_2})^2(\phi_{\mu_1,\mu_2})=\left(\mu_1^2-(2\mu_1+1)\mu_2\frac{z_2}{z_1-z_2}
+\mu_2(\mu_2-1)\frac{z_2^2}{(z_1-z_2)^2}\right)\phi_{\mu_1,\mu_2}.\eqno(2.25)$$
Thus
$$[(z_1\ptl_{z_1})^2+(z_2\ptl_{z_2})^2](\phi_{\mu_1,\mu_2})-2\mu_2(\mu_2-1)\frac{z_1z_2}
{(z_1-z_2)^2}\phi_{\mu_1,\mu_2}=(2\mu_1^2+2\mu_1\mu_2+\mu_2^2)
\phi_{\mu_1,\mu_2},\eqno(2.26)$$
which is of the form (1.12) with $n=2$ and
$$K=-2\mu_2(\mu_2-1),\qquad \nu=2\mu_1^2+2\mu_1\mu_2+\mu_2^2.\eqno(2.27)$$

Take
$$g(z)=z^r\sum_{m=0}^{\infty}a_mz^m\qquad\mbox{with}\;\;r,a_m\in\mbb{C}.\eqno(2.28)$$
Then
$$z_1\ptl_{z_1}\left(g\left(\frac{z_2}{z_2-z_1}\right)\right)=\frac{z_1z_2}{(z_2-z_1)^2}
g'\left(\frac{z_2}{z_2-z_1}\right),\eqno(2.29)$$
$$z_2\ptl_{z_2}\left(g\left(\frac{z_2}{z_2-z_1}\right)\right)=-\frac{z_1z_2}{(z_2-z_1)^2}
g'\left(\frac{z_2}{z_2-z_1}\right),\eqno(2.30)$$
$$(z_1\ptl_{z_1})^2\left(g\left(\frac{z_2}{z_2-z_1}\right)\right)=\frac{z_1^2z_2^2}
{(z_1-z_2)^4}g''\left(\frac{z_2}{z_2-z_1}\right)+\frac{z_1z_2(z_1+z_2)}{(z_2-z_1)^3}
g'\left(\frac{z_2}{z_2-z_1}\right),\eqno(2.31)$$
$$(z_2\ptl_{z_2})^2\left(g\left(\frac{z_2}{z_1-z_2}\right)\right)=\frac{z_1^2z_2^2}
{(z_1-z_2)^4}g''\left(\frac{z_2}{z_2-z_1}\right)-\frac{z_1z_2(z_1+z_2)}{(z_1-z_2)^3}
g'\left(\frac{z_2}{z_2-z_1}\right).\eqno(2.32)$$ Hence by (2.22),
(2.23), (2.26) and (2.29)-(2.32), we have
\begin{eqnarray*}& &[(z_1\ptl_{z_1})^2+(z_2\ptl_{z_2})^2]\left (\phi_{\mu_1,\mu_2}
g\left(\frac{z_2}{z_2-z_1}\right)\right)\\ &=&\phi_{\mu_1,\mu_2}[\left(
2\mu_2(\mu_2-1)\frac{z_1z_2}{(z_1-z_2)^2}+(2\mu_1^2+2\mu_1\mu_2+\mu_2^2)\right)
g\left(\frac{z_2}{z_2-z_1}\right)\\
& &+
2\left(\mu_1+\mu_2\frac{z_1}{z_1-z_2}\right)\frac{z_1z_2}{(z_1-z_2)^2}g'\left(\frac{z_2}
{z_2-z_1}\right)
\\ &
&-2\left(\mu_1-\mu_2\frac{z_2}{z_1-z_2}\right)\frac{z_1z_2}{(z_1-z_2)^2}g'\left(\frac{z_2}
{z_2-z_1}\right)\\
&&+2\frac{z_1^2z_2^2}{(z_1-z_2)^4}g''\left(\frac{z_2}{z_2-z_1}\right)+2\frac{z_1z_2(z_1+z_2)}
{(z_2-z_1)^3}g'\left(\frac{z_2}{z_1-z_2}\right)]\\
&=& \phi_{\mu_1,\mu_2}[\left(
2\mu_2(\mu_2-1)\frac{z_1z_2}{(z_1-z_2)^2}+(2\mu_1^2+2\mu_1\mu_2+\mu_2^2)\right)g
\left(\frac{z_2}{z_2-z_1}\right)\\ &
&+2(1-\mu_2)\frac{z_1z_2(z_1+z_2)}{(z_2-z_1)^3}
g'\left(\frac{z_2}{z_1-z_2}\right)+2\frac{z_1^2z_2^2}{(z_1-z_2)^4}g''\left(\frac{z_2}{z_2-z_1}
\right)].\hspace{2.1cm}(2.33)\end{eqnarray*}

Recall the classical Gauss hypergeometric equation
$$x(1-x)\frac{d^2y}{dx^2}+[c-(a+b+1)x]\frac{dy}{dx}-aby=0,\eqno(2.34)$$
whose fundamental solutions are
$$ _2F_1(a,b;c; x)\qquad\mbox{and}\qquad x^{1-c}\: _2F_1(a+1-c,b+1-c;2-c; x)\eqno(2.35)$$
(see (1.17) and e.g. cf. [AAR]).
Observe that
$$\frac{z_1+z_2}{z_2-z_1}=2\frac{z_2}{z_2-z_1}-1,\;\;\frac{z_2}{z_2-z_1}\left(1-\frac{z_2}
{z_2-z_1}\right)=-\frac{z_1z_2}{(z_1-z_2)^2},\eqno(2.36)$$
Therefore, a sufficient  condition for $\phi_{\mu_1,\mu_2}g(z_2/(z_2-z_1))$ to be a solution
of  the equation (1.12)  with $n=2$ for suitable $K$ and $\nu$  is
$$(1-\mu_2)\frac{z_1+z_2}{z_2-z_1}g'\left(\frac{z_2}{z_1-z_2}\right)+\frac{z_1z_2}{(z_1-z_2)^2}
g''\left(\frac{z_2}{z_2-z_1}\right)=\ell g\left(\frac{z_2}{z_2-z_1}\right)\eqno(2.37)$$
for some $\ell\in\mbb{C}$, equivalently,
$$(1-\mu_2)(2z-1)g'(z)-z(1-z)g''(z)=\ell g(z),\eqno(2.38)$$
which is of the form (2.34) with
$$a+b+1=2(1-\mu_2),\;\;c=1-\mu_2,\;\;ab=-\ell.\eqno(2.39)$$
Take $a\in\mbb{C}$ and
$$b=1-2\mu_2-a,\qquad \ell=a(a+2\mu_2-1).\eqno(2.40)$$

The second fundamental solution of (2.38) is
$$z^{\mu_2}\: _2F_1(a+\mu_2,1-\mu_2-a;\mu_2+1; z)\eqno(2.41)$$
(cf. (2.35)). In fact, by (2.20),
$$\Psi(z_1,z_2)=(-1)^{\mu_2}\phi_{\mu_1,\mu_2}\left(\frac{z_2}{z_2-1}\right)^{\mu_2}\:
_2F_1\left(a+\mu_2,1-\mu_2-a;\mu_2+1;\frac{z_2}{z_2-1}\right)\eqno(2.42)$$
(cf. (1.14)) with
$$\mu=\mu_2+a-1,\;\;\lmd_1=\mu_1-\frac{1}{2},\;\;\lmd_2=\mu_1+\mu_2+\frac{1}{2}.\eqno(2.43)$$
So the Etingof trace gives only one family of solutions of the form
$\phi_{\mu_1,\mu_2}g(z_2/(z_2-z_1))$ of the equation (1.12) with $n=2$ and
  $$K=2\mu_2(1-\mu_2)+2a(a+2\mu_2-1),\qquad \nu=2\mu_1^2+2\mu_1\mu_2+\mu_2^2.\eqno(2.44)$$
 \pse

{\bf Theorem 2.2}. {\it For  $g(z)$ in (2.28), the product  function
$\phi_{\mu_1,\mu_2}(z_1,z_2)g(z_2/(z_2-z_1))$ is a solution of  the equation
(1.12) with $n=2$ and suitable constants $K,\mu$ if $g(z)$ is a solution of the
 classical hypergeometric equation (2.38). In particular, the functions
$$ (z_1z_2)^{\mu_1}(z_1-z_2)^{\mu_2}\;_2F_1(a,1-2\mu_2-a;1-\mu_2;z)\eqno(2.45)$$
for $a,\mu_1,\mu_2\in\mbb{C}$ are a new  family of solutions of  the equation (1.12)
 with $n=2$ and $K,\nu$ in (2.44).}

\section{Etingof Trace for $gl(n)$ with $n>2$}

In this section, we find the Etingof trace functions for the
general case of $gl(n)$ with $n>2$ under a certain condition and
their master hypergeometric functions.

 In the Lie algebra $gl(n)$ given in (1.4) and (1.5), we choose
$$\{E_{i,j}\mid 1\leq i<j\leq n\}\;\;\mbox{as positive root vectors}.\eqno(3.1)$$
In particular, we have
$$\{E_{i,i+1}\mid i=1,2,...,n-1\}\;\;\mbox{as  positive simple root vectors}.\eqno(3.2)$$
Accordingly,
$$\{E_{i,j}\mid 1\leq j<i\leq n\}\;\;\mbox{are negative root vectors}\eqno(3.3)$$
and
$$\{E_{i+1,i}\mid i=1,2,...,n-1\}\;\;\mbox{are negative simple root vectors}.\eqno(3.4)$$

Recall  $\G_A$ in (1.19). Let $U(gl(n))$ be the universal enveloping algebra of $gl(n)$.
 For
$$\al =\sum_{1\leq j<i\leq n}\al_{ij}\es_{i,j}\in \G_A,\eqno(3.5)$$
we denote
$$E^\al=E_{2,1}^{\al_{2,1}}E_{3,1}^{\al_{3,1}}E_{3,2}^{\al_{3,2}}E_{4,1}^{\al_{4,1}}\cdots
E_{n,1}^{\al_{n,1}}\cdots E_{n,n-1}^{\al_{n,n-1}}\in U(gl(n)).\eqno(3.6)$$
Denote by  ${\cal G}_-$ the Lie subalgebra spanned by (3.3). Then
$$\{E^\al\mid \al\in\G_A\}\;\;\mbox{forms a PBW basis of}\;\;U({\cal G}_-).\eqno(3.7)$$

Let $\lmd$ be a weight of $gl(n)$
such that
$$\lmd_1-\lmd_2=\cdots =\lmd_{n-2}-\lmd_{n-1}=\mu\;\;\mbox{and}\;\;\lmd_{n-1}-\lmd_n
\not\in\mbb{N},
\eqno(3.8)$$ for some constant $\mu$. Set
$$\sgm_i=\lmd_i-\lmd_{i+1}\qquad\for\;\;i\in\ol{1,n-1}.\eqno(3.9)$$
 The Verma $gl(n)$-module with the highest-weight
vector $v_\lmd$ of weight $\lmd$ is given by
$$M_\lmd =\mbox{Span}\{E^\al v_\lmd\mid\al\in \G_A\},\eqno(3.10)$$
with the action determined by
\begin{eqnarray*}E_{i,i+1}(E^\al v_\lmd)&=&(\sum_{p=1}^{i-1}\al_{i+1,p}E^{\al+\es_{i,p}
-\es_{i+1,p}}-\sum_{p=i+2}^n\al_{p,i}E^{\al+\es_{p,i+1}-\es_{p,i}}\\
& &+\al_{i+1,i}
(\sgm_i+1-\sum_{p=i+1}^n\al_{p,i}+\sum_{p=i+2}^n\al_{p,i+1})E^{\al-\es_{i+1,i}})v_\lmd,
\hspace{2.2cm}(3.11)\end{eqnarray*}
$$E_{j,i}(E^\al v_\lmd)=(E^{\al+\es_{j,i}}+\sum_{p=1}^{i-1}\al_{i,p}E^{\al+\es_{j,p}
-\es_{i,p}})
v_\lmd,\eqno(3.12)$$
$$E_{k,k}(E^\al v_\lmd)=(\lmd_k+\sum_{p=1}^{k-1}\al_{k,p}-\sum_{q=k+1}^n\al_{q,k})
E^\al v_\lmd
\eqno(3.13)$$ for $1\leq i<j\leq n$ and $k=1,2,...,n$.

Let
$$\td{M}=\{\sum_{i_1,\cdots,i_n\in\mbb{Z}}v_{i_1,\cdots,i_n}x_1^{i_1+\mu}x_2^{i_2+\mu}
\cdots
 x_n^{i_n+\mu}\mid v_{i_1,\cdots,i_n}\in M_\lmd\},\eqno(3.14)$$
a space of Laurent series with coefficients in $M_\lmd$. Define the action of $gl(n)$ on
$\td{M}$ by
\begin{eqnarray*}& &E_{j_1,j_2}(\sum_{i_1,\cdots,i_n\in\mbb{Z}}v_{i_1,\cdots,i_n}
x_1^{i_1+\mu}
x_2^{i_2+\mu}\cdots
x_n^{i_n+\mu})=\sum_{i_1,\cdots,i_n\in\mbb{Z}}[E_{j_1,j_2}
(v_{i_1,\cdots,i_n})x_1^{i_1+\mu}x_2^{i_2+\mu}\cdots
x_n^{i_n+\mu}\\ & &+v_{i_1,\cdots,i_n}
(x_{j_1}\ptl_{x_{j_2}}-\mu\dlt_{j_1,j_2})(x_1^{i_1+\mu}x_2^{i_2+\mu}\cdots
x_n^{i_n+\mu})] \hspace{5.9cm}(3.15)\end{eqnarray*} for
$j_1,j_2\in\ol{1,n}$. Then $\td{M}$ becomes a $gl(n)$-module that
is isomorphic to $M_\lmd \hat{\otimes}V_\mu$,  the completion of
the tensor $M_\lmd \otimes V_\mu$ with respect to the topology
given by weight subspaces (cf. (1.6) and  (1.7)).

For
$$\vec i=(i_1,i_2,...,i_{n-1})\in\mbb{N}^{n-1},\eqno(3.16)$$
we denote
$$ x^{(\vec i)}=x_1^{i_1}x_2^{i_2-i_1}\cdots x_{n-1}^{i_{n-1}-i_{n-2}}x_n^{-i_{n-1}},
\qquad \es(\vec i)=\sum_{p=1}^{n-1}i_p\es_{p+1,p}.\eqno(3.17)$$
Then we have: \psp

{\bf Lemma 3.1}. {\it The vector
$$u_\lmd=(x_1x_2\cdots x_n)^{\mu}\sum_{i_1,...,i_{n-1}=0}^{\infty}\frac{(-1)^{i_1+i_2+\cdots
+i_{n-1}}\la \mu\ra_{i_{n-1}}}{i_1!i_2!\cdots i_{n-1}!\la
\sgm\ra_{i_{n-1}}}E^{\es(\vec i)} v_\lmd x^{(\vec i)}\eqno(3.18)$$
is a singular vector of weight $\lmd$, where
$\sgm=\sgm_{n-1}=\lmd_{n-1}-\lmd_n$.}

{\it Proof}. For $p\in\{1,2,...,n-2\}$, we have
\begin{eqnarray*} & &E_{p,p+1}(u_\lmd)\\ &=&
(x_1x_2\cdots
x_n)^{\mu}\sum_{i_1,...,i_{n-1}=0}^{\infty}\frac{(-1)^{i_1+i_2+\cdots+i_{n-1}}
\la \mu\ra_{i_{n-1}}}{i_1!i_2!\cdots i_{n-1}!\la
\sgm\ra_{i_{n-1}}}[(\mu+i_{p+1}-i_p) E^{\es(\vec i)}v_\lmd
x_px_{p+1}^{-1}x^{(\vec i)}\\ & &+(\mu+1-i_p+i_{p+1})i_p
E^{\es(\vec i)-\es_p}v_\lmd x^{(\vec i)}]\\ &=& (x_1x_2\cdots
x_n)^{\mu}[\sum_{i_1,...,i_{n-1}=0}^{\infty}\frac{(-1)^{i_1+i_2+\cdots+i_{n-1}}
\la \mu\ra_{i_{n-1}}}{i_1!i_2!\cdots i_{n-1}!\la
\sgm\ra_{i_{n-1}}}(\mu+i_{p+1}-i_p) E^{\es(\vec i)}v_\lmd
x_px_{p+1}^{-1}x^{(\vec i)}\\ &&+
\sum_{i_1,...,i_{n-1}=0}^{\infty}\frac{(-1)^{i_1+i_2+\cdots+i_{n-1}}\la
\mu\ra_{i_{n-1}}} {i_1!i_2!\cdots i_{n-1}!\la
\sgm\ra_{i_{n-1}}}(\mu+1-i_p+i_{p+1})i_pE^{\es(\vec i)-\es_p}
v_\lmd x^{(\vec i)}]\hspace{3cm}\end{eqnarray*}
\begin{eqnarray*}
&=&(x_1x_2\cdots x_n)^{\mu}[\sum_{i_1,...,i_{n-1}=0}^{\infty}
\frac{(-1)^{i_1+i_2+\cdots+i_{n-1}}\la
\mu\ra_{i_{n-1}}}{i_1!i_2!\cdots i_{n-1}! \la
\sgm\ra_{i_{n-1}}}(\mu+i_{p+1}-i_p)E^{\es(\vec i)}v_\lmd
x_px_{p+1}^{-1}x^{(\vec i)}
\\ &&-\sum_{i_1,...,i_{n-1}=0}^{\infty}\frac{(-1)^{i_1+i_2+\cdots+i_{n-1}}
\la \mu\ra_{i_{n-1}}}
{i_1!i_2!\cdots i_{n-1}!\la
\sgm\ra_{i_{n-1}}}(\mu+i_{p+1}-i_p)E^{\es(\vec i)}v_\lmd x_p
x_{p+1}^{-1}x^{(\vec i)}]=0.\hspace{1.1cm}(3.19)\end{eqnarray*}
Moreover,
\begin{eqnarray*}& & E_{n-1,n}(u_\lmd)\\ &=&
(x_1x_2\cdots
x_n)^{\mu}\sum_{i_1,...,i_{n-1}=0}^{\infty}\frac{(-1)^{i_1+i_2+\cdots+i_{n-1}}\la
\mu\ra_{i_{n-1}}}{i_1!i_2!\cdots i_{n-1}!\la
\sgm\ra_{i_{n-1}}}[(\mu-i_{n-1})E^{\es(\vec i)}v_\lmd
x_{n-1}x_n^{-1}x^{(\vec i)}\\ &
&+(\sgm+1-i_{n-1})i_{n-1}E^{\es(\vec i)-\es_{n-1}}v_\lmd x^{(\vec
i)}]\\ &=& (x_1x_2\cdots
x_n)^{\mu}[\sum_{i_1,...,i_{n-1}=0}^{\infty}\frac{(-1)^{i_1+i_2+\cdots+i_{n-1}}\la
\mu\ra_{i_{n-1}}}{i_1!i_2!\cdots i_{n-1}!\la
\sgm\ra_{i_{n-1}}}(\mu-i_{n-1})E^{\es(\vec i)}v_\lmd
x_{n-1}x_n^{-1}x^{(\vec i)}\\ &&+
\sum_{i_1,...,i_{n-1}=0}^{\infty}\frac{(-1)^{i_1+i_2+\cdots+i_{n-1}}\la
\mu\ra_{i_{n-1}}}{i_1!i_2!\cdots i_{n-1}!\la
\sgm\ra_{i_{n-1}}}(\sgm+1-i_{n-1})i_{n-1}E^{\es(\vec
i)-\es_{n-1}}v_\lmd x^{(\vec i)}]\\ &=&(x_1x_2\cdots
x_n)^{\mu}[\sum_{i_1,...,i_{n-1}=0}^{\infty}\frac{(-1)^{i_1+i_2+\cdots+i_{n-1}}\la
\mu\ra_{i_{n-1}}}{i_1!i_2!\cdots i_{n-1}!\la
\sgm\ra_{i_{n-1}}}(\mu-i_{n-1})E^{\es(\vec i)}v_\lmd
x_{n-1}x_n^{-1}x^{(\vec i)}\\
&&-\sum_{i_1,...,i_{n-1}=0}^{\infty}\frac{(-1)^{i_1+i_2+\cdots+i_{n-1}}\la
\mu\ra_{i_{n-1}+1}}{i_1!i_2!\cdots i_{n-1}!\la
\sgm\ra_{i_{n-1}}}E^{\es(\vec i)}v_\lmd x_px_{p+1}^{-1}x^{(\vec
i)}]=0.\hspace{3.1cm}(3.20)\end{eqnarray*} Hence $u_\lmd$ is a
singular vector. $\qquad\Box$ \psp

 Let $\Phi: M_\lmd\rta\td{M}$ be the Lie algebra module homomorphism
such that $\Phi(v_\lmd)=u_\lmd$. Let $\vec i$ as in (3.16) and
$\be\in\G_A$. By (3.11),  we have
\begin{eqnarray*}& &E_{n,n-1}^{\be_{n,n-1}}(E^{\es(\vec i)}v_\lmd)=
\sum_{p_{n,n-1}=0}^{\be_{n,n-1}}
\left(\!\!\begin{array}{c}\be_{n,n-1}\\ p_{n,n-1}\end{array}\!\!\right)\la i_{n-2}\ra_{p_{n,n-1}}
\\ & &E^{\es(\vec i)-p_{n,n-1}\es_{n-1,n-2}+p_{n,n-1}\es_{n,n-2}+(\be_{n,n-1}-p_{n,n-1})
\es_{n,n-1}}v_\lmd,\hspace{5cm}(3.21)\end{eqnarray*}
\begin{eqnarray*}& &E_{n,n-2}^{\be_{n,n-2}}E_{n,n-1}^{\be_{n,n-1}}(E^{\es(\vec i)}v_\lmd)
=\sum_{p_{n,n-2},p_{n,n-1}=0}^{\infty}
\left(\!\!\begin{array}{c}\be_{n,n-2}\\
p_{n,n-2}\end{array}\!\!\right)\left(\!\!\begin{array}
{c}\be_{n,n-1}\\ p_{n,n-1}\end{array}\!\!\right)\la
i_{n-3}\ra_{p_{n,n-2}}\la i_{n-2} \ra_{p_{n,n-1}}\\ & &
E^{\es(\vec
i)+\sum_{r=n-2,n-1}[-p_{n,r}\es_{r,r-1}+p_{n,r}\es_{n,r-1}+(\be_{n,r}-p_{n,r})
\es_{n,r}]}v_\lmd,\hspace{5.2cm}(3.22)\end{eqnarray*}
\begin{eqnarray*}\qquad& &\left[\prod_{r=1}^{n-1}E_{n,r}^{\be_{n,r}}\right]
(E^{\es(\vec
i)}v_\lmd)=\sum_{p_{n,2},...,,p_{n,n-1}=0}^{\infty}\left[\prod_{r=2}^{n-1}
\left(\!\!\begin{array}{c}\be_{n,r}\\
p_{n,r}\end{array}\!\!\right)\la i_{r-1}\ra_{p_{n,r}}\right]\\ & &
E^{\es(\vec
i)+\be_{n,1}\es_{n,1}+\sum_{s=2}^{n-1}(-p_{n,s}\es_{s,s-1}+p_{n,s}\es_{n,s-1}
+(\be_{n,s}-p_{n,s})\es_{n,s})}v_\lmd,\hspace{3.8cm}(3.23)\end{eqnarray*}
\begin{eqnarray*}&&\left[\prod_{r=1}^{n-2}E_{n-1,r}^{\be_{n-1,r}}\right]\left[\prod_{r=1}^{n-1}
E_{n,r}^{\be_{n,r}}\right](E^{\es(\vec
i)}v_\lmd)=\sum_{i=n-1,n}\;\sum_{j=1}^{i-1}
\sum_{p_{i,j}=0}^{\infty}\hspace{7cm}\end{eqnarray*}
\begin{eqnarray*} &
&\left[\prod_{r=1}^{n-2}\left(\!\!\begin{array}{c}\be_{n-1,r}
\\ p_{n-1,r}\end{array}\!\!\right)\right]
\left[\prod_{r=1}^{n-1}\left(\!\!\begin{array}{c}\be_{n,r}\\ p_{n,r}\end{array}\!\!\right)
\right]
\left[\prod_{r=1}^{n-3}\la i_r\ra_{p_{n,r+1}+p_{n-1,r+1}}\right]\la i_{n-2}\ra_{p_{n,n-1}}
\\ & &E^{\es(\vec i)+\sum_{r=n-1,n}\be_{r,1}\es_{r,1}+\sum_{r=n-1,n}\sum_{s=2}^{r-1}(-p_{r,s}
\es_{s,s-1}+p_{r,s}\es_{r,s-1}+(\be_{r,s}-p_{r,s})\es_{r,s})}v_\lmd,\hspace{2.3cm}(3.24)
\end{eqnarray*}
\begin{eqnarray*}& &E^\be(E^{\es(\vec i)}v_\lmd)=\sum_{2\leq k<j\leq n}\;
\sum_{p_{j,k}=0}^{\infty}
\left[\prod_{2\leq m<l\leq n}\left(\!\!\begin{array}{c}\be_{l,m}\\
p_{l,m}\end{array}\!\!\right) \right]\left[\prod_{s=1}^{n-2}\la
i_s\ra_{_{\sum_{r=s+2}^np_{r,s+1}}}\right]\\ & & E^{\es(\vec
i)+\sum_{1\leq m<l\leq
n}(\be_{l,m}+(1-\dlt_{l,m+1})p_{l,m+1}-p_{l,m}-\dlt_{l,m+1}
\sum_{s=l+1}^np_{s,l})\es_{l,m}}v_\lmd.\hspace{3.8cm}(3.25)\end{eqnarray*}
where we treat $p_{2,1}=p_{3,1}=\cdots=p_{n,1}=0$. Hence by (3.14)
and (3.25), we obtain
\begin{eqnarray*}& &E^\al(E^{\es(\vec i)}v_\lmd (x_1x_2\cdots x_n)^{\mu}x^{(\vec i)})\\
&=&\sum_{p_{2,1}^0,...,p_{n,1}^0=0}^{\infty} \;\sum_{2\leq k<j\leq
n}\;\sum_{p_{j,k}^0,p_{j,k}=0}^{\infty}\left[\prod_{m=2}^n\left(\!\!\begin{array}{c}\al_{m,1}\\
p_{m,1}^0\end{array}\!\!\right)\right]\left[\prod_{2\leq m<l\leq
n}\left(\!\!\begin{array}{c}\al_{l,m}\\
p_{l,m}^0,p_{l,m}\end{array}\!\!\right)\right]\\
&&\left[\prod_{s=1}^{n-2}\la
i_s\ra_{_{\sum_{r=s+2}^np_{r,s+1}}}\right]\left[\prod_{s=1}^{n-1}\la
\mu+i_s-i_{s-1}\ra_{_{\sum_{r=s+1}^np_{r,s}^0}}\right]\\ & &
E^{\es(\vec i)+\sum_{1\leq m<l\leq
n}(\al_{l,m}+p_{l,m+1}-p_{l,m}^0-p_{l,m}
-\dlt_{l,m+1}\sum_{s=l+1}^np_{s,l})\es_{l,m}}v_\lmd\\ & &
(x_1x_2\cdots
x_n)^{\mu}\prod_{r=1}^nx_r^{i_r-i_{r-1}+\sum_{s=1}^{r-1}p^0_{r,s}
-\sum_{q=r+1}^np_{q,r}^0}\hspace{6cm}(3.26)\end{eqnarray*} for
$\al\in\G_A$,  where we treat
$$i_{-1}=i_n=p_{2,1}=p_{3,1}=\cdots =p_{n,1}=p_{2,2}=p_{3,3}=\cdots p_{n,n}=0.\eqno(3.27)$$

View $\Phi$  as a function in $\{x_1,x_2,...,x_n\}$ taking values in $\mbox{End}\:M_\lmd$, the
space of linear transformations on $M_\lmd$. The {\it Etingof trace function}
$$E(z_1,z_2,...,z_n)=(x_1x_2\cdots x_n)^{-\mu}\mbox{tr}_{_{M_\lmd}}\:\Phi z_1^{E_{1,1}}
z_2^{E_{2,2}}\cdots z_n^{E_{n,n}}.\eqno(3.28)$$ In order to
calculate $E(z_1,z_2,...,z_n)$, we need to take
$$p_{l,m+1}=p_{l,m}^0+p_{l,m},\;\; i_k=p_{k+1,k}^0+p_{k+1,k}+\sum_{s=k+2}^np_{s,k+1}\eqno(3.29)$$
in (3.26) for $1\leq m<l-1\leq n-1$ and $k\in\ol{1,n-1}$. By
(3.27), (3.29) and induction,
$$p_{l,m}=\sum_{r=1}^{m-1}p_{l,r}^0\qquad\for\;\;1\leq m<l\leq n.\eqno(3.30)$$
Moreover,
$$i_k=\sum_{r=k+1}^n(p_{r,k}^0+p_{r,k})=\sum_{r=k+1}^n\sum_{s=1}^kp^0_{r,s}\qquad\for\;\;
k\in\ol{1,n-1}.\eqno(3.31)$$
$$i_k-i_{k-1}=\sum_{r=k+1}^n\sum_{s=1}^kp^0_{r,s}-\sum_{r=k}^n\sum_{s=1}^{k-1}p^0_{r,s}=
\sum_{r=k+1}^np^0_{r,k}-\sum_{s=1}^{k-1}p^0_{k,s}\eqno(3.32)$$ for
$k\in\ol{1,k}$. Furthermore,
$$\left(\!\!\begin{array}{c}\al_{j,k}\\ p_{j,k}^0,p_{j,k}\end{array}\!\!\right)=\left(\!\!
\begin{array}{c}p^0_{j,k}+p_{j,k}\\ p_{j,k}\end{array}\!\!\right)\left(\!\!\begin{array}{c}
\al_{j,k}\\ p^0_{j,k}+p_{j,k}\end{array}\!\!\right) \eqno(3.33)$$
and
$$ \left(\!\!\begin{array}{c}p^0_{j,k}+p_{j,k}\\ p_{j,k}\end{array}\!\!\right)=
\frac{(p^0_{j,k}+p_{j,k})!}{p^0_{j,k}!p_{j,k}!}=\frac{(p^0_{j,k}+p_{j,k})!}{p^0_{j,k}!
(p^0_{j,k-1}+p_{j,k-1})!}\eqno(3.34)$$ by (3.29).

For convenience, we denote
$$\al_{_{\ul{1}}}=\al_{\ol{n}}=0,\;\;\al_{_{\ul{i}}}=\sum_{r=1}^{i-1}\al_{i,r},\;\;
\al_{_{\ol{j}}}=\sum_{s=j+1}^n\al_{s,j}\eqno(3.35)$$ for
$\al\in\G_A,\;i\in \ol{2,n}$ and $j\in\ol{1,n-1}$. By (2.16),
(3.18) and (3.26)-(3.35), the Etingof's trace function
\begin{eqnarray*}& &E(z_1,z_2,...,z_n)
\\ &=&\left(\prod_{r=1}^nz_r^{\lmd_r}\right)\sum_{\al\in\G_A}\;\sum_{p_{r_1,r_2}^0\in\mbb{N},
\;1\leq r_2<r_1\leq n}(-1)^{\sum_{s=1}^{n-1}i_s}\left[\prod_{2\leq
k<j\leq n}\left(\!\!
\begin{array}{c}p^0_{j,k}+p_{j,k}\\ p_{j,k}\end{array}\!\!\right)\right]\\ &&\times
\frac{\prod_{s=1}^{n-1}\la\mu+i_s-i_{s-1}\ra_{_{\sum_{r=s+1}^np^0_{r,s}}}}
{\la\sgm\ra_{_{i_{n-1}}}\prod_{s=1}^{n-1}(p_{s+1,s}^0+p_{s+1.s})!}
\left[\prod_{1\leq k<j\leq n}\left(\!\!\begin{array}{c}\al_{j,k}\\ p^0_{j,k}+p_{j,k}
\end{array}\!\!\right)\right] \la\mu\ra_{i_{n-1}}
\prod_{1\leq k<j\leq n}\left(\frac{z_j}{z_k}\right)^{\al_{j,k}}\\
&=&\left(\prod_{r=1}^nz_r^{\lmd_r}\right)\sum_{p_{r_1,r_2}^0\in\mbb{N},\;1\leq
r_2<r_1\leq
n}(-1)^{\sum_{s=1}^{n-1}i_s}\frac{\prod_{s=1}^{n-1}\la\mu+i_s-i_{s-1}\ra_{_{\sum_{r=s+1}^n
p^0_{r,s}}}\la\mu\ra_{i_{n-1}}}{\la\sgm\ra_{_{i_{n-1}}}\prod_{1\leq
s_1<s_2\leq n}p_{s_2,s_1}^0!}\\ & &\times \left[\prod_{1\leq
r_1<r_2\leq n}\left(\frac{z_{r_2}}{z_{r_1}}\right)^{p_{r_2,r_1}^0
+p_{r_2,r_1}}\right]\left[\prod_{1\leq s_1<s_2\leq
n}\left(\frac{1}{1-z_{s_2}/z_{s_1}}
\right)^{p_{s_2,s_1}^0+p_{s_2,s_1}+1}\right] \\
&=&\frac{\prod_{r=1}^nz_r^{\lmd_r+n-r}}{\prod_{1\leq k<j\leq
n}(z_k-z_j)}\sum_{\be\in\G_A}\frac{\left[\prod_{s=1}^{n-1}\la\mu+\be_{\ol{s}}-\be_{\ul{s}}
\ra_{_{\be_{\ol{s}}}}\right]\la\mu\ra_{_{\be_{\ul{n}}}}}{\la\sgm\ra_{_{\be_{\ul{n}}}}
\prod_{1\leq s_1<s_2\leq n}\be_{s_2,s_1}!}\\ & &\times
\left[\prod_{1\leq k<j\leq
n}\left(\frac{-z_k}{z_k-z_j}\right)^{\sum_{r=1}^k\be_{j,r}}\right].\hspace{8.2cm}(3.36)
\end{eqnarray*}
Set
$$\xi_{r_2,r_1}^A=\prod_{s=r_1}^{r_2-1}\frac{z_{r_2}}{z_{r_2}-z_s}
\qquad \for\;\;1\leq r_1<r_2\leq n.\eqno(3.37)$$

Motivated by (2.17)-(2.19) and (3.36), we define our
$(n(n-1)/2)$-variable hypergeometric function of type A  by
$${\cal X}_A(\tau_1,..,\tau_n;\vt)\{z_{j,k}\}=\sum_{\be\in\G_A}\frac{\left[\prod_{s=1}^{n-1}
(\tau_s-\be_{\ul{s}})_{_{\be_{\ol{s}}}}\right](\tau_n)_{_{\be_{\ul{n}}}}}
{\be!(\vt)_{_{\be_{\ul{n}}}}}z^\be\eqno(3.38)$$ (cf. (3.35)),
where
$$\be!=\prod_{1\leq k<j\leq n}\be_{j,k}!,\qquad z^\be=\prod_{1\leq k<j\leq n}z_{j,k}^{\be_{j,k}}.
\eqno(3.39)$$
 According to (2.18) and (3.36)-(3.38), we have:
\psp

{\bf Theorem 3.3}. {\it Under the condition (3.8), the Etingof
trace function   is}
$$E(z_1,z_2,...,z_n)=\left[\frac{\prod_{r=1}^nz_r^{\lmd_r+n-r}}{\prod_{1\leq k<j\leq n}
(z_k-z_j)}\right]{\cal
X}_A(\mu+1,..,\mu+1,-\mu;-\sgm)\{\xi_{r_2,r_1}^A\}.\eqno(3.40)$$
\pse

The function in (1.14)
$$\Psi(z_1,z_2,\cdots,z_n)=\prod_{r=1}^nz_r^{\lmd_r+(n+1)/2-r}
{\cal
X}_A(\mu+1,..,\mu+1,-\mu;-\sgm)\{\xi_{r_2,r_1}^A\}.\eqno(3.41)$$

\section{Etingof Trace for $sp(2n)$}

In this section, we want to introduce and calculate  the Etingof trace function of $sp(2n)$.

The symplectic Lie algebra
\begin{eqnarray*}\hspace{1cm}sp(2n)&=&\sum_{i,j=1}^n\mbb{C}(E_{i,j}-E_{n+j,n+i})
+\sum_{i=1}^n(\mbb{C}E_{i,n+i}+\mbb{C}E_{n+i,i})\\  & &+\sum_{1\leq i<j\leq n }[\mbb{C}
(E_{i,n+j}+E_{n+j,i})+\mbb{C}(E_{n+i,j}+E_{n+j,i})]\hspace{3cm}(4.1)\end{eqnarray*}
is a Lie subalgebra of $gl(2n)$. Set
$$h_n=E_{n,n}-E_{2n,2n},\;\;h_i=E_{i,i}-E_{i+1,i+1}-E_{n+i,n+i}+E_{n+i+1,n+i+1}\eqno(4.2)$$
for $i=1,2,...,n-1.$
The subspace
$$T=\sum_{i=1}^n\mbb{C}h_i\eqno(4.3)$$
forms a Cartan subalgebra of $sp(2n)$. We choose positive root vectors
$$\{E_{i,j}-E_{n+j,n+i},\;E_{i,n+j}+E_{j,n+i},\;E_{k,n+k} \mid 1\leq i<j\leq n,\; k\ol{1,n}\}.
\eqno(4.4)$$
In particular, we take
$$\{E_{i,i+1}-E_{n+i+1,n+i},\;E_{n,2n}\mid i=1,2,...,n-1\}\;\;\mbox{as positive simple root
vectors}.\eqno(4.5)$$
Accordingly,
$$\{E_{i,j}-E_{n+j,n+i},\;E_{i,n+j}+E_{j,n+i},\;E_{n+k,k} \mid 1\leq j<i\leq n,\;k\in\ol{1,n}\}
\eqno(4.6)$$
are negative root vectors and
$$\{E_{i+1,i}-E_{n+i,n+i+1},\;E_{2n,n}\mid i\in\ol{1,n-1}\}\;\;\mbox{are negative simple root
vectors}.\eqno(4.7)$$

For convenience, we denote
$$C_{i,j}=E_{i,j}-E_{n+j,n+i}\qquad\for\;\;i,j\in\ol{1.n},\eqno(4.8)$$
$$C_{j,n+k}=E_{j,n+k}+E_{k,n+j},\qquad C_{n+j,k}=E_{n+j,k}+E_{n+k,j}\eqno(4.9)$$
for $j,k\in\ol{1,n}$ such that $j\neq k$, and
$$C_{r,n+r}=E_{r,n+r},\qquad C_{n+r,r}=E_{n+r,r}\qquad\for\;\;r\in\ol{1,n}.\eqno(4.10)$$
Let
$$\G_C=\sum_{1\leq j<i\leq n}\mbb{N}\es_{i,j}+\sum_{1\leq j\leq i\leq n}\mbb{N}\es_{n+i,j}
\eqno(4.11)$$
be the torsion-free additive semigroup  of  rank $n^2$ with $\es_{p,q}$ as base elements, and
let $U(sp(2n))$ be the universal enveloping algebra of $sp(2n)$. For any
$$\al=\sum_{1\leq j<i\leq n}\al_{i,j}\es_{i,j}+\sum_{1\leq j\leq i\leq n}\al_{n+i,j}\es_{n+i,j}
\in \G_C,\eqno(4.12)$$
we denote
\begin{eqnarray*}\hspace{1cm}C^\al&=&C_{2,1}^{\al_{2,1}}C_{3,1}^{\al_{3,1}}C_{3,2}^{\al_{3,2}}
C_{4,1}^{\al_{4,1}}\cdots C_{n,1}^{\al_{n,1}}\cdots C_{n,n-1}^{\al_{n,n-1}}\\ & &\times
 C_{n+1,1}^{\al_{n+1,1}}C_{n+2,1}^{\al_{n+2,1}}C_{n+2,2}^{\al_{n+2,2}}C_{n+3,1}^{\al_{n+3,1}}
 \cdots C_{2n,1}^{\al_{2n,1}}\cdots C_{2n,n}^{\al_{2n,n}}\hspace{3.6cm}(4.13)\end{eqnarray*}
Denote by  ${\cal G}_-$  the Lie subalgebra spanned by (4.6). Then
$$\{C^\al\mid \al\in\G_C\}\;\;\mbox{forms a PBW basis of}\;\;U({\cal G}_-)).\eqno(4.14)$$

For convenience, we treat
$$\es_{n+i,j}=\es_{n+j,i}\qquad\for\;\;i,j=1,2,...,n.\eqno(4.15)$$
Note that
$$[C_{n,2n},C_{2n,n}]=h_n,\qquad [C_{i,i+1},C_{i+1,i}]=h_i,\eqno(4.16)$$
$$[h_i,C_{n+i+1,i}]=0,\qquad [h_i,C_{n+i,i}]=-2C_{n+i,i},\qquad [h_i,C_{n+i+1,i+1}]=2C_{n+i+1,i+1},\eqno(4.17)$$
$$[h_i,C_{n+k,i}]=-C_{n+k,i},\qquad[h_i,C_{n+k,i+1}]=C_{n+k,i+1},\eqno(4.18)$$
$$[C_{i,i+1},C_{n+p,i}]=-C_{n+p,i+1},\qquad [C_{i,i+1},C_{n+i+1,i}]=-2C_{n+i+1,i+1},\eqno(4.19)$$
$$[C_{n,2n},C_{2n,i}]=C_{n,i},\qquad [C_{n+j,r}, C_{r,q}]=(1+\dlt_{j,q})C_{n+j,q},
\qquad[C_{n+r,r},C_{r,q}]=C_{n+r,q}\eqno(4.20)$$
for $i\in\ol{1,n-1}$ and $ j,k,p,q,r\in\ol{1,n}$ such that  $k\neq i,i+1,\;p\neq i+1$ and $q< r$.

 Let $\lmd$ be a weight, which is a linear function on $T$, such that
$$\lmd(h_i)=-\frac{1}{2}\qquad\for\;\;i=1,2,...,n-1;\qquad \lmd(h_n)=\lmd_n\in\mbb{C}. \eqno(4.21)$$
Recall that $sp(2n)$ is generated by
$$\{C_{i,i+1},C_{i+1,i},C_{n,2n},C_{2n,n},\mid i=1,2,...,n-1\}\eqno(4.22)$$
as a Lie algebra. The Verma  $sp(2n)$-module with the highest-weight vector $v_\lmd$ of weight $\lmd$ is given by
$$M_\lmd =\mbox{Span}\{E^\al v_\lmd\mid\al\in \G_C\},\eqno(4.23)$$
with the action determined by
\begin{eqnarray*}& &C_{i,i+1}(C^\al v_\lmd)\\ &=&[\sum_{j=1}^{i-1}\al_{i+1,j}C^{\al+\es_{i,j}
-\es_{i+1,j}}-\sum_{j=i+2}^n\al_{j,i}C^{\al+\es_{j,i+1}-\es_{j,i}}-\sum_{k\neq i+1}\al_{n+k,i}
C^{\al-\es_{n+k,i}+\es_{n+k,i+1}}\\ & &-2\al_{n+i+1,i}C^{\al-\es_{n+i+1,i}+\es_{n+i+1,i+1}}+
\al_{i+1,i}(1/2-\sum_{j=i+1}^n\al_{j,i}+\sum_{j=i+2}^n\al_{j,i+1}\\ & &+\sum_{k\neq i,i+1}
(\al_{n+k,i+1}-\al_{n+k,i})-2\al_{n+i,i}+2\al_{n+i+1,i+1})C^{\al-\es_{i+1,i}}]v_\lmd
\hspace{2.6cm}(4.24)\end{eqnarray*}
for $i\in \ol{1,n-1}$,
\begin{eqnarray*}& &C_{n,2n}(C^\al v_\lmd)=[\sum_{i=1}^{n-1}\al_{2n,i}(C^{\al-\es_{2n,i}
+\es_{n,i}}+(\al_{2n,i}-1)C^{\al-2\es_{2n,i}+\es_{n+i,i}}\\ & &+\sum_{j=1}^{i-1}\al_{2n,j}
C^{\al-\es_{2n,i}-\es_{2n,j}+\es_{n+i,j}})+\al_{2n,n}(\lmd_n+1-\al_{2n,n})C^{\al-\es_{2n,n}}]
v_\lmd\hspace{2.5cm}(4.25)\end{eqnarray*}
and
$$C_{j,i}(E^\al v_\lmd)=(C^{\al+\es_{j,i}}+\sum_{p=1}^{i-1}\al_{i,p}C^{\al+\es_{j,p}
-\es_{i,p}})v_\lmd,\eqno(4.26)$$
\begin{eqnarray*}C_{n+j,i}(E^\al v_\lmd)&=&(C^{\al+\es_{n+j,i}}+\sum_{q=1}^{j-1}\al_{n+j,q}
C^{\al+\es_{n+q,i}-\es_{j,q}}\\ & &+
\sum_{p=1}^{i-1}\al_{i,p}(C^{\al+\es_{n+j,p}-\es_{i,p}}+\sum_{q=1}^{j-1}\al_{n+j,q}C^{\al+
\es_{n+q,p}-\es_{i,p}-\es_{j,q}})\hspace{2cm}(4.27)\end{eqnarray*}
for $1\leq i<j\leq n$,
\begin{eqnarray*}\hspace{1cm}C_{n+k,k}(C^\al v_\lmd)&=&[C^{\al+\es_{n+k,k}}+\sum_{r=1}^{k-1}
\al_{k,r}(C^{\al-\es_{k,r}+\es_{n+k,r}}+(\al_{k,r}-1)C^{\al-2\es_{k,r}+\es_{n+r,r}}\\ &&
+\sum_{s=r+1}^{k-1}\al_{k,s}C^{\al-\es_{k,r}-\es_{k,s}+\es_{n+s,r}})]v_\lmd\hspace{4.6cm}
(4.28)\end{eqnarray*}
 by (4.16)-(4.20).

Let $\{x_1,x_2,...,x_n\}$ be $n$ indeterminates. Set
$$x^{\vec i}=x_1^{i_1}x_2^{i_2}\cdots x_n^{i_n}\qquad\vec i=(i_1,i_2,...,i_n)\in\mbb{Z}^n
\eqno(4.29)$$
and
$$x^\ast=x_1^{-1/2}x_2^{-1/2}\cdots x_n^{-1/2}.\eqno(4.30)$$
Denote by
$$\td{M}=\{\sum_{\vec i\in\mbb{Z}^n}w_{\vec i}x^{\vec i}x^\ast\mid w_{\vec i}\in M_\lmd\},
\eqno(4.31)$$
the space of formal Laurent series in $\{x_1,x_2,...,x_n\}$ with the coefficients in $M_\lmd$.
 Moreover, we define the action of $sp(2n)$ on $\td{M}$ by
$$C_{p,q}(\sum_{\vec i\in\mbb{Z}^n}w_{\vec i}x^{\vec i}x^\ast)=\sum_{\vec i\in\mbb{Z}^n}
\left(C_{p,q}(w_{\vec i})x^{\vec i}x^\ast+w_{\vec i}x_p\ptl_{x_q}(x^{\vec i}x^\ast)
+\frac{\dlt_{p,q}}{2}w_{\vec i}x^{\vec i}x^\ast\right),\eqno(4.32)$$
$$C_{p,n+q}(\sum_{\vec i\in\mbb{Z}^n}w_{\vec i}x^{\vec i}x^\ast)=\sum_{\vec i\in\mbb{Z}^n}
\left(C_{p,n+q}(w_{\vec i})x^{\vec i}x^\ast-\frac{1}{1+\dlt_{p,q}}w_{\vec i}x_px_qx^{\vec i}
x^\ast\right),\eqno(4.33)$$
$$C_{n+p,q}(\sum_{\vec i\in\mbb{Z}^n}w_{\vec i}x^{\vec i}x^\ast)=\sum_{\vec i\in\mbb{Z}^n}
(C_{n+p,q}(w_{\vec i})x^{\vec i}x^\ast+\frac{1}{1+\dlt_{p,q}}w_{\vec i}\ptl_{x_p}\ptl_{x_q}
(x^{\vec i}x^\ast))\eqno(4.34)$$
for $1\leq p,q\leq n$. It can be verified that $\td{M}$ forms an $sp(2n)$-module (cf. [F]).

For $\vec i\in \mbb{N}\:^n$, we set
$$x^{(\vec i)}=x_1^{i_1}x_2^{i_2-i_1}x_3^{i_3-i_2}\cdots x_{n-1}^{i_{n-1}-i_{n-2}}x_n^{2i_n
-i_{n-1}},\eqno(4.35)$$
$$\es(\vec i)=\sum_{p=1}^{n-1}i_p\es_{p+1,p}+i_n\es_{2n,n}.\eqno(4.36)$$
Let
$$u_\lmd=\sum_{\vec i\in\mbb{N}\:^n}\frac{(-1)^{i_1+i_2+\cdots+i_{n-1}}}{i_1!i_2!\cdots i_n!
2^{i_n}\la \lmd_n\ra_{i_n}}C^{\es(\vec i)}v_\lmd x^{(\vec i)}x^\ast.\eqno(4.37)$$
It is verified that $u_\lmd$ is a singular vector of weight $\lmd$ in $\td{M}$, that is,
$$C_{i,j}(u_\lmd)=C_{p,n+q}(u_\lmd)=0\qquad\for\;\;1\leq i<j\leq n;\;1\leq p\leq q\leq n.
\eqno(4.38)$$
Hence there exists $sp(2n)$-module homomorphism $\Phi: M_\lmd\rta \td{M}$ such
that
$$\Phi(v_\lmd)=u_\lmd.\eqno(4.39)$$
View $\Phi$ as a function in $\{x_1,x_2,...,x_n\}$ taking value in
the spaces of  linear transformations on $M_\lmd$. The Etingof
trace function is defined as
$$E_C(z_1,z_2,...,z_n)=(x^\ast)^{-1}\mbox{tr}_{M_\lmd}\Phi z_1^{C_{1,1}} z_2^{C_{2,2}}\cdots
 z_n^{C_{n,n}}.\eqno(4.40)$$
In the rest of this section, we want to calculate $E_C(z_1,z_2,...,z_n)$.

For $\vec i\in\mbb{N}\:^n$ and $\be \in \G_C$ (cf. (4.11)), we
have
\begin{eqnarray*}& &C^\be(C^{\es(\vec i)}v_\lmd)=\sum_{2\leq j_1<j_2\leq n}\;\sum_{k_1=1}^n
\;\sum_{k_2=2}^n\;\sum_{k_3=k_2}^n(\sum_{p_{j_2,j_1},q_{k_1,k_2},q_{k_2,k_2}^1=0}^{\infty}\\
& &\left[\prod_{2\leq l_1<l_2\leq
n}\left(\!\!\begin{array}{c}\be_{l_2,l_1}\\
p_{l_2,l_1}\end{array}\!\!\right)\left(\!\!\begin{array}{c}\be_{n+l_2,l_1}\\
q_{l_1,l_2},q_{l_2,l_1},
q^1_{l_2,l_1}\end{array}\!\!\right)\right]
\left[\prod_{r=2}^n\left(\!\!\begin{array}{c}\be_{n+r,r}\\
q_{r,r},q_{r,r}^1\end{array}\!\!\right)\left(\!\!\begin{array}{c}\be_{m+r,1}\\
q_{1,r}\end{array}\!\!\right)\right]\\ & &\times
\left[\prod_{s=1}^{n-2}\la
i_s\ra_{_{\sum_{r=s+2}^np_{r,s+1}+\sum_{r=1}^nq_{r,s+1}+\sum_{r=2}^{s+1}q^1_{s+1,r}+
\sum_{r=s+1}^nq^1_{r,s+1}
}}\right]\la
i_{n-1}\ra_{_{\sum_{r=1}^nq_{r,n}+\sum_{r=2}^nq^1_{n,r}+q_{n,n}^1}}
\\ & &\times
E^{\sum_{1\leq l_1<l_2\leq
n}(p_{l_2,l_1+1}-p_{l_2,l_1}-\dlt_{l_2,l_1+1}(\sum_{r=l_2+1}^n
p_{r,l_2}+\sum_{s=1}^nq_{s,l_2}+\sum_{s=2}^{l_2}q^1_{l_2,s}+\sum_{s=l_2}^nq^1_{s,l_2}))
\es_{l_2,l_1}}
\\ & &\times E^{\sum_{1\leq s_1<s_2\leq n}(q_{s_2,s_1+1}+q_{s_1,s_2+1}+q_{s_2+1,s_1+1}^1
-q_{s_1.s_2}-q_{s_2,s_1}-q_{s_2,s_1}^1)\es_{n+s_2,s_1}}\\ &
&\times E^{\es(\vec i)+\be
+\sum_{r=1}^n(q_{r,r+1}+q_{r+1,r+1}^1-q_{r,r}-q_{r,r}^1)\es_{n+r,r}}v_\lmd)\hspace{6.4cm}
(4.41)\end{eqnarray*}
by the calculations in (3.21)-(3.26) and (4.26)-(4.28), where we
treat
$$p_{r,r}=q_{1,1}=q^1_{n+1,n+1}=q^1_{r,1}=q^1_{n+1,r}=q_{s,n+1}=q_{n+1,s}=0\eqno(4.42)$$
for $2\leq r\leq n$ and $1\leq s\leq n.$ Hence by (4.32), (4.34) and (4.41), we obtain
\begin{eqnarray*}& &C^\al(C^{\es(\vec i)}v_\lmd x^{(\vec i)}x^\ast)\\ &=&
\sum_{2\leq j_1<j_2\leq
n}\;\sum_{k_1=1}^n\;\sum_{k_2=2}^n\;\sum_{k_3=k_1}^n\;\sum_{k_4=k_2}^n\{\sum_{p^0_{k_2,1},
p_{j_2,j_1}^0,p_{j_2,j_1},q^0_{k_3,k_1},q^1_{k_3,k_1},q_{k_1,k_2}=0}^{\infty}\frac{1}
{2^{\sum_{r=1}^nq^0_{r,r}}}\\
& &\times \left[\prod_{2\leq l_1<l_2\leq
n}\left(\!\!\begin{array}{c}\al_{l_2,l_1}\\  p_{l_2,l_1}^0,
p_{l_2,l_1}\end{array}\!\!\right)\left(\!\!\begin{array}{c}\al_{n+l_2,l_1}\\
q_{l_2,l_1}^0,q_{l_2,l_1}^1,
q_{l_1,l_2},q_{l_2,l_1}\end{array}\!\!\right)\right]\\ &
&\times\left[\prod_{m=2}^n\left(\!\!\begin{array}{c}\al_{m,1}\\
p_{m,1}^0\end{array}\!\!\right)\right]\left(\!\!\begin{array}{c}\al_{n+1,1}\\
q_{1,1}^0\end{array}\!\!\right)\left[\prod_{r=2}^n\left(\!\!\begin{array}{c}\al_{n+r,r}\\
q_{r,r}^0,
q_{r,r}^1,q_{r,r}\end{array}\!\!\right)\left(\!\!\begin{array}{c}\al_{n+r,1}\\
q_{r,1}^0,q_{1,r}\end{array}\!\!\right)\right]\\ & &\times
\left[\prod_{s=1}^{n-2}\la
i_s\ra_{_{\sum_{r=s+2}^np_{r,s+1}+\sum_{r=1}^nq_{r,s+1}+\sum_{r=2}^{s+1}q^1_{s+1,r}+
\sum_{r=s+1}^nq^1_{r,s+1}
}}\right]\la
i_{n-1}\ra_{_{\sum_{r=1}^nq_{r,n}+\sum_{r=2}^nq^1_{n,r}+q_{n,n}^1}}\\
& &\times \left[\prod_{s=1}^{n-1}\la
i_s-i_{s-1}-1/2\ra_{_{\sum_{r=s+1}^np_{r,s}^0+\sum_{r=1}^sq_{s,r}^0+\sum_{r=s}^nq_{r,s}^0}}
\right]\la
2i_n-i_{n-1}-1/2\ra_{_{2q^0_{n,n}+\sum_{r=1}^{n-1}q^0_{n,r}}}
\\ & &\times
E^{\sum_{1\leq l_1<l_2\leq
n}(p_{l_2,l_1+1}-p_{l_2,l_1}-p_{l_2,l_1}^0-\dlt_{l_2,l_1+1}
(\sum_{r=l_2+1}^np_{r,l_2}+\sum_{s=1}^nq_{s,l_2}+\sum_{s=2}^{l_2}q^1_{l_2,s}+\sum_{s=l_2}^n
q^1_{s,l_2}))\es_{l_2,l_1}}
\\ & &\times E^{\sum_{1\leq s_1<s_2\leq n}(q_{s_2,s_1+1}+q_{s_1,s_2+1}+q_{s_2+1,s_1+1}^1
-q_{s_1.s_2}-q_{s_2,s_1}-q_{s_2,s_1}^0-q_{s_2,s_1}^1)\es_{n+s_2,s_1}}\\ & &\times
E^{\al+\es(\vec i)+\sum_{r=1}^n(q_{r,r+1}+q_{r+1,r+1}^1-q_{r,r}-q_{r,r}^0-q_{r,r}^1)
\es_{n+r,r}}v_\lmd\\ & &
\left[\prod_{s=1}^{n-1}x_s^{i_s-i_{s-1}-1/2+\sum_{r=1}^{s-1}p^0_{s,r}-\sum_{r=s+1}^np_{r,s}^0
-\sum_{r=1}^sq_{s,r}^0-\sum_{r=s}^nq_{r,s}^0}\right]\\ & &\times x_n^{2i_n-i_{n-1}-1/2+
\sum_{r=1}^{n-1}p^0_{n,r}-2q^0_{n,n}-\sum_{r=1}^{n-1}q^0_{n,r}}\},\hspace{6.3cm}(4.43)
\end{eqnarray*}
where we treat $i_{-1}=0$.

 In order to  calculate  the Etingof trace (4.40), we have to let in (4.43):
$$ p_{r_2,r_1}=p_{r_2,r_1-1}^0+p_{r_2,r_1-1}\qquad\for\;\;2\leq r_1<r_2\leq n,\eqno(4.44)$$
$$q_{r,r+1}+q_{r+1,r+1}^1=q^0_{r,r}+q_{r,r}++q_{r,r}^1\qquad\for\;\;r\in\ol{1,n-1},\eqno(4.45)$$
$$q_{s_2,s_1+1}+q_{s_1,s_2+1}+q_{s_2+1,s_1+1}^1=q_{s_1.s_2}+q_{s_2,s_1}+q_{s_2,s_1}^0
+q_{s_2,s_1}^1\eqno(4.46)$$
for $1\leq s_1<s_2\leq n$ and
$$i_r=p^0_{r+1,r}+p_{r+1,r}+\sum_{s=r+2}^np_{s,r+1}+\sum_{s=1}^nq_{s,r+1}
+\sum_{s=2}^{r+1}q_{r+1,s}^1+\sum_{s=r+1}^nq_{s,r+1}^1\eqno(4.47)$$
for $r\in\ol{1,n-1}$,
$$i_n=q_{n,n}^0+q_{n,n}+q^1_{n,n}.\eqno(4.48)$$
In particular, (4.44) implies
$$p_{r_2,r_1}=\sum_{s=1}^{r_1-1}p^0_{r_2,s}\qquad\for\;\;2\leq r_1<r_2\leq n.\eqno(4.49)$$

Set
$$\Im_r=\sum_{s=1}^nq_{s,r+1}+\sum_{s=2}^{r+1}q_{r+1,s}^1+\sum_{s=r+1}^nq_{s,r+1}^1\eqno(4.50)$$
for $r\in\ol{1,n-1}$.
By (4.42) and (4.45), we have
$$q_{1,2}+q^1_{2,2}=q_{1,1}^0.\eqno(4.51)$$
Moreover, (4.46) implies
$$q_{i,2}+q_{1,i+1}+q_{i+1,2}^1=q_{1.i}+q_{i,1}^0\qquad\for\;\;i\in\ol{2,n}.\eqno(4.52)$$
$$2(q_{1,2}+q^1_{2,2})+\sum_{i=2}^nq_{i,2}+\sum_{r=3}^n(q_{1,r}+q_{r,2}^1)=2q^0_{1,1}
+\sum_{i=2}^n(q_{1,i}+q^0_{i,1}),\eqno(4.53)$$
equivalently,
$$\sum_{i=1}^nq_{i,2}+q^1_{2,2}+\sum_{r=2}^nq_{r,2}^1=q^0_{1,1}+\sum_{i=1}^nq^0_{i,1}.\eqno(4.54)$$
Thus
$$\Im_1=q^0_{1,1}+\sum_{i=1}^nq^0_{i,1}.\eqno(4.55)$$

Let $r\in\ol{2,n-1}$. By (4.46), we have
$$q_{i,r+1}+q_{r,i+1}+q_{i+1,r+1}^1=q_{r.i}+q_{i,r}+q_{i,r}^0+q_{i,r}^1\eqno(4.56)$$
for $i\in\ol{r+1,n}$ and
$$q_{r,j+1}+q_{j,r+1}+q_{r+1,j+1}^1=q_{j.r}+q_{r,j}+q_{r,j}^0+q_{r,j}^1\eqno(4.57)$$
for $j\in\ol{1,r-1}$. Expressions (4.45), (4.56) and (4.57) imply
\begin{eqnarray*} & &2(q_{r,r+1}+q_{r+1,r+1}^1)+\sum_{i=r+1}^n(q_{i,r+1}+q_{r,i+1}
+q_{i+1,r+1}^1)+\sum_{j=1}^{r-1}(q_{r,j+1}+q_{j,r+1}+q_{r+1,j+1}^1)\\
&=& 2(q^0_{r,r}+q_{r,r}++q_{r,r}^1)
+q_{i,r}^1)+\sum_{j=1}^{r-1}(q_{j.r}+q_{r,j}+q_{r,j}^0+q_{r,j}^1)\\&
&
+\sum_{s=r+1}^n(q_{r,s}+q_{s,r}+q_{s,r}^0+q_{s,r}^1),\hspace{8.4cm}(4.58)\end{eqnarray*}
equivalently,
\begin{eqnarray*}\qquad & &\sum_{i=1}^nq_{i,r+1}+\sum_{j=2}^{r+1}q_{r+1,j}^1+\sum_{i=r+1}^n q_{i,r+1}
\\ &=&\sum_{j=1}^rq_{r,j}^0+\sum_{s=r}^nq_{s,r}^0+\sum_{s=1}^nq_{s,r}+\sum_{j=2}^rq_{r,j}^1
+\sum_{s=r}^nq_{s,r}^1\hspace{4.9cm}(4.59)\end{eqnarray*} by
(4.42). So
$$\Im_r-\Im_{r-1}=\sum_{j=1}^rq_{r,j}^0+\sum_{i=r}^nq_{i,r}^0.\eqno(4.60)$$
Furthermore, by (4.47), (4.49), (4.50), (4.60) and induction, we
have
$$
i_r=\sum_{s=r+1}^n\sum_{j=1}^rp_{s,j}^0+\sum_{s=1}^r(\sum_{j=1}^sq_{s,j}^0+\sum_{j=s}^nq_{j,s}^0)
\eqno(4.61)$$ for $r\in\ol{1,n-1}$. In particular,
$$
i_r-i_{r-1}=\sum_{s=r+1}^np^0_{s,r}-\sum_{s=1}^{r-1}p^0_{r,s}+\sum_{j=1}^rq_{r,j}^0
+\sum_{s=r}^nq_{s,r}^0\eqno(4.62)$$ for $r\in\ol{1,n-1}$.

Next, we want to calculate $i_n$. First, we have
$$q_{n.n}+q^1_{n,n}=q_{n,n-1}+q_{n-1,n}+q_{n,n}^1+q_{n,n-1}^1+q_{n,n-1}^0\eqno(4.63)$$
by (4.46). Moreover, for $r\in\ol{1,n-2}$, we have
\begin{eqnarray*}& &\sum_{i=0}^{2r}q_{n-i,n-2r+i}+\sum_{i=0}^{r-1}q_{n-i,n-2r+1+i}
\\ &=&q_{n,n-2r}+\sum_{i=1}^r(q_{n-i,n-2r+i}+q_{n-2r-1+i,n-i+1}+q_{n-i+1,n-2r+i}^1)\\ &=&
\sum_{i=1}^r(q_{n-i,n-2r-1+i}+q_{n-2r-1+i,n-i}+q_{n-i,n-2r-1+i}^0+q_{n-i,n-2r-1+i}^1)
\hspace{4cm}\end{eqnarray*}
\begin{eqnarray*}
 && +q_{n,n-2r-1}+q_{n-2r-1,n}+q_{n,n-2r-1}^0+q_{n,n-2r-1}^1
\\ &=& \sum_{i=0}^r(q_{n-i,n-2r-1+i}+q_{n-2r-1+i,n-i})+\sum_{i=0}^rq_{n-i,n-2r-1+i}^1+
 \sum_{i=0}^rq_{n-i,n-2r-1+i}^0\\ &=&\sum_{i=0}^{2r+1}q_{n-i,n-2r-1+i}+ \sum_{i=0}^r
 q_{n-i,n-2r-1+i}^0+ \sum_{i=0}^rq_{n-i,n-2r-1+i}^1.\hspace{3.7cm}(4.64)\end{eqnarray*}
by (4.46). Furthermore, for $r\in\ol{0,n-2}$, we get
\begin{eqnarray*}& &\sum_{i=0}^{2r+1}q_{n-i,n-2r-1+i}+\sum_{i=0}^rq_{n-i,n-2r+i}\\
&=&q_{n,n-2r-1}+q_{n-r-1,n-r}+q_{n-r,n-r}^1+\sum_{i=1}^r(q_{n-i,n-2r-1+i}\\
& & +q_{n-2r-2+i,n+1-i} +q_{n-i+1,n-2r-1+i}^1)
\\ &=& q_{n,n-2r-2}+q_{n-2r-2,n}+q^0_{n,n-2r-2}+q^1_{n,n-2r-2}+q_{n-r-1,n-r-1}
+q_{n-r-1,n-r-1}^0\\ & &+q_{n-r-1,n-r-1}^1
+\sum_{i=1}^r(q_{n-i,n-2r-2+i}+q_{n-2r-2+i,n-i}+q_{n-i,n-2r-2+i}^0+q_{n-i,n-2r-2+i}^1)
\\ &=&
q_{n-r-1,n-r-1}+\sum_{i=0}^r(q_{n-i,n-2r-2+i}+q_{n-2r-2+i,n-i})
\\ & &+\sum_{i=0}^{r+1}q_{n-i,n-2r-2+i}^0+\sum_{i=0}^{r+1}q_{n-i,n-2r-2+i}^1\\
&=& \sum_{i=0}^{2r+2}
q_{n-i,n-2r-2+i}+\sum_{i=0}^{r+1}q_{n-i,n-2r-2+i}^0+\sum_{i=0}^{r+1}q_{n-i,n-2r-2+i}^1
\hspace{3.9cm}(4.65)\end{eqnarray*} based on (4.45) and (4.46).

By (4.63)-(4.65) and induction, we obtain
$$q_{n,n}+q_{n,n}^1=\sum_{r=0}^{n-1}\sum_{s=0}^rq_{n-s,n-2r-1+s}^0+
\sum_{r=0}^{n-2}\sum_{s=0}^rq_{n-s,n-2r-2+s}^0=\sum_{1\leq r_1\leq
r_2\leq n}q^0_{r_2,r_1} -q_{n,n}^0.\eqno(4.66)$$ Therefore,
$$i_n=\sum_{1\leq r_1\leq r_2\leq n}q_{s_2,s_1}^0.\eqno(4.67)$$
Now
$$2i_n-i_{n-1}=\sum_{1\leq r_1\leq r_2\leq n}2q_{s_2,s_1}^0-\sum_{r=1}^{n-1}(p^0_{n,r}
+\sum_{p=1}^rq_{r,i}^0+\sum_{s=r}^nq_{s,r}^0)=q_{n,n}^0+\sum_{r=1}^nq_{n,r}^0
-\sum_{s=1}^{n-1}p^0_{n,s}.\eqno(4.68)$$

By (4.45) and induction on $r$, we get
$$q_{r,r}^1=\sum_{s=1}^{r-1}(q_{s,s}^0+q_{s,s}-q_{s,s+1})\qquad\for\;\;r\in\ol{2,n}.\eqno(4.69)$$
Moreover, (4.46) and induction imply
$$q_{n,r}=\sum_{s=1}^{r-1}(q^0_{n,s}+q_{s,n}+q_{n,s}^1)\qquad\for\;\;r\in\ol{2,n}\eqno(4.70)$$
and
$$q_{r_2,r_1}^1=\sum_{s=1}^{r_2-1}(q_{r_2-s,r_1-s}^0+q_{r_2-s,r_1-s}+q_{r_1-s,r_2-s})
-\sum_{s=1}^{r_1-1}(q_{r_2-s,r_1-s+1}+q_{r_1-s,r_2-s+1})\eqno(4.71)$$
for $2\leq r_1<r_2\leq n$. Furthermore, (4.70) and (4.71) show
\begin{eqnarray*}\qquad q_{n,r}&=&\sum_{s=1}^{r-1}(q_{n,i}^0+q_{i,n})+\sum_{s=2}^{r-1}\sum_{j=1}^{s-1}(q_{n-j,s-j}^0+
q_{n-j,s-j}+q_{s-j,n-j})\\ &
&-\sum_{s=2}^{r-1}\sum_{j=1}^{s-1}(q_{n-j,s-j+1}+q_{s-j,n-j+1})\hspace{6.4cm}
(4.72)\end{eqnarray*} for $r\in\ol{2,n}$. Set
$$\vec{p^0}!=\prod_{1\leq r_1<r_2\leq n}p_{r_2,r_1}^0!,\qquad \vec{q^0}!=
\prod_{1\leq r_1\leq r_2\leq n}q_{r_2,r_1}^0!.\eqno(4.73)$$

By (4.40), (4.43), (4.61), (4.67), (4.69), (4.71) and (4.72), we
have
\begin{eqnarray*}& &E_C(z_1,z_2,...,z_n)\\ &=& \sum_{\al\in
\G_C}\; \sum_{1\leq j_1<j_2\leq
n}\;\sum_{k_1=1}^n\;\sum_{k_2=2}^n\;\sum_{k_3=k_1}^n\;\sum_{k_4=1}^{n-1}\;
\sum_{p_{j_2,j_1}^0,
q_{k_3,k_1}^0,q_{k_4,k_2}=0}^{\infty}\frac{(-1)^{i_1+i_2+\cdots+i_{n-1}}}{i_1!i_2!\cdots
i_n! 2^{i_n+\sum_{r=1}^nq^0_{r,r}}\la \lmd_n\ra_{i_n}}
\\ & &\times \left[\prod_{2\leq l_1\leq l_2\leq n}\left(\!\!\begin{array}{c}\al_{l_2,l_1}\\
p^0_{l_2,l_1},
p_{l_2,l_1}\end{array}\!\!\right)\left(\!\!\begin{array}{c}\al_{n+l_2,l_1}\\
q^0_{l_2,l_1},q_{l_2,l_1}^1,q_{l_1,l_2},q_{l_2,l_1}\end{array}\!\!\right)\right]\\
& &\times\left[\prod_{m=2}^n\left(\!\!\begin{array}{c}\al_{m,1}\\
p^0_{m,1}\end{array}\!\!\right)\right]\left(\!\!\begin{array}{c}\al_{n+1,1}\\
q^0_{1,1}\end{array}\!\!\right)\left[\prod_{r=2}^n\left(\!\!\begin{array}{c}\al_{n+r,r}\\
q^0_{r,r},q_{r,r}^1,q_{r,r}\end{array}\!\!\right)\left(\!\!\begin{array}{c}\al_{n+r,1}\\
q^0_{r,1},q_{1,r}\end{array}\!\!\right)\right]\\ & &\times
\left[\prod_{s=1}^{n-2}\la
i_s\ra_{_{\sum_{r_1=s+2}^np_{r_1,s+1}+\sum_{r=1}^nq_{r,s+1}+\sum_{r=2}^{s+1}q_{s+1,r}^1
+\sum_{r=s+1}^nq_{r,s+1}^1}}\right] \la
i_{n-1}\ra_{_{\sum_{r=1}^nq_{r,n}+\sum_{r=2}^nq_{n,r}^1+q^1_{n,n}}}\\
& &\times \left[\prod_{s=1}^{n-1}\la i_s
-i_{s-1}-1/2\ra_{_{\sum_{r=s+1}^np_{r,s}^0+\sum_{r=1}^sq_{s,r}^0+\sum_{r=s}^nq_{r,s}^0}}\right]\la
2i_n-i_{n-1}-1/2\ra_{_{2q^0_{n,n}+\sum_{r=1}^{n-1}q^0_{n,r}}}\\
& &\times\left[\prod_{1\leq r_1<r_2\leq
n}\left(\frac{z_{r_2}}{z_{r_1}}\right)^{\al_{r_2,r_1}}\frac{1}{(z_{r_1}z_{r_2})
^{\al_{n+r_2,r_1}}}\right]\left[\prod_{r=1}^n\frac{1}{z_r^{\al_{n+r,r}}}\right](z_1z_2\cdots
z_{n-1})^{-1/2}z_n^{\lmd_n}\\ &=&\sum_{1\leq j_1<j_2\leq
n}\;\sum_{k_1=1}^n\;\sum_{k_2=2}^n\;\sum_{k_3=k_1}^n\;\sum_{k_4=1}^{n-1}
\;\sum_{p_{j_2,j_1}^0,
q_{k_3,k_1}^0,q_{k_4,k_2}=0}^{\infty}\frac{(-1)^{i_1+i_2+\cdots+i_{n-1}}}{\vec
p^0! 2^{i_n+\sum_{r=1}^nq^0_{r,r}}\la \lmd_n\ra_{i_n}}\\ &
&\times\frac{1}{q_{n,n}^0!q_{n,n}!q_{n,n}^1!} \left[\prod_{1\leq
l_1\leq l_2\leq n}
\left(\!\!\begin{array}{c}q^0_{l_2,l_1}+q_{l_2,l_1}^1+q_{l_1,l_2}+q_{l_2,l_1}
\\ q_{l_2,l_1}^1,q_{l_1,l_2},q_{l_2,l_1}\end{array}\!\!\right)\right]\left[\prod_{r=2}^{n-1}
\left(\!\!\begin{array}{c}q^0_{r,r}+q_{r,r}^1+q_{r,r}\\
q_{r,r}^1,q_{r,r}\end{array}\!\!\right)\right]\end{eqnarray*}
\begin{eqnarray*} & &\times
\left[\prod_{s=1}^{n-1}\la
i_s-i_{s-1}-1/2\ra_{_{\sum_{r=s+1}^np_{r,s}^0+\sum_{r=1}^sq_{s,r}^0+\sum_{r=s}^nq_{r,s}^0}}
\right]\la 2i_n-i_{n-1}-1/2\ra_{_{2q^0_{n,n}+\sum_{r=1}^{n-1}q^0_{n,r}}}\\
& &\times \left[\prod_{1\leq r_1<r_2\leq
n}\frac{z_{r_1}^2z_{r_2}^{1+\sum_{s=1}^{r_1}p_{r_2,s}^0}}{(z_{r_1}-z_{r_2})
^{1+\sum_{s=1}^{r_1}p_{r_2,s}^0}(z_{r_1}z_{r_2}-1)^{1+q^0_{r_2,r_1}+q_{r_2,r_1}^1
+q_{r_1,r_2}+q_{r_2,r_1}}}\right]\\ & &\times
\left[\prod_{r=1}^n\frac{z_r^2}{(z_r^2-1)^{1+q^0_{r,r}+q^1_{r,r}+q_{r,r}}}\right]
(z_1z_2\cdots
z_{n-1})^{-1/2}z_n^{\lmd_n}.\hspace{5cm}(4.74)\end{eqnarray*} Note
that
\begin{eqnarray*}& &\frac{1}{q_{n,n}^0!q_{n,n}!q_{n,n}^1!} \left[\prod_{1\leq
l_1\leq l_2\leq n}
\left(\!\!\begin{array}{c}q^0_{l_2,l_1}+q_{l_2,l_1}^1+q_{l_1,l_2}+q_{l_2,l_1}
\\ q_{l_2,l_1}^1,q_{l_1,l_2},q_{l_2,l_1}\end{array}\!\!\right)\right]\left[\prod_{r=2}^{n-1}
\left(\!\!\begin{array}{c}q^0_{r,r}+q_{r,r}^1+q_{r,r}\\
q_{r,r}^1,q_{r,r}\end{array}\!\!\right)\right]\\
&=&\frac{1}{q_{n,n}^0!q_{n,n}!q_{n,n}^1!} \left[\prod_{1\leq
l_1\leq l_2\leq n}
\left(\!\!\begin{array}{c}q_{l_2+1,l_1+1}^1+q_{l_1,l_2+1}+q_{l_2,l_1+1}
\\ q_{l_2,l_1}^1,q_{l_1,l_2},q_{l_2,l_1}\end{array}\!\!\right)\right]\left[\prod_{r=2}^{n-1}
\left(\!\!\begin{array}{c}q_{r+1,r+1}^1+q_{r,r+1}\\
q_{r,r}^1,q_{r,r}\end{array}\!\!\right)\right]\\ &=&\frac{1}{\vec
q^0!} \left[\prod_{1\leq l_1\leq l_2\leq n-1}
\left(\!\!\begin{array}{c}q_{l_2+1,l_1+1}^1+q_{l_1,l_2+1}+q_{l_2,l_1+1}
\\ q_{l_1,l_2+1},q_{l_2,l_1+1}\end{array}\!\!\right)\right]\left[\prod_{r=1}^{n-1}
\left(\!\!\begin{array}{c}q_{r+1,r+1}^1+q_{r,r+1}\\
q_{r,r+1}\end{array}\!\!\right)\right]\\ &=&\frac{1}{\vec q^0!}
\left[\prod_{1\leq l_1\leq l_2\leq n-1}
\left(\!\!\begin{array}{c}q^0_{l_2,l_1}+q_{l_2,l_1}^1+q_{l_1,l_2}+q_{l_2,l_1}
\\ q_{l_1,l_2+1},q_{l_2,l_1+1}\end{array}\!\!\right)\right]\left[\prod_{r=1}^{n-1}
\left(\!\!\begin{array}{c}q^0_{r,r}+q_{r,r}^1+q_{r,r}\\
q_{r,r+1}\end{array}\!\!\right)\right]\hspace{0.9cm}(4.75)\end{eqnarray*}
according to (4.45) and (4.46). Furthermore,
$$q_{r,r}^0+q^1_{r,r}+q_{r,r}=\sum_{s=1}^r(q_{s,s}^0+q_{s,s})-\sum_{s=1}^{r-1}q_{s,s+1}
\eqno(4.76)$$ by (4.69) and
\begin{eqnarray*}& & q_{n,r}^0+q_{n,r}^1+q_{r,n}+q_{n,r}\\
&=&\sum_{s=1}^r(q^0_{n,s}+q_{s,n})+\sum_{s=0}^{r-2}\sum_{j=1}^{r-s-1}(q^0_{n-j,r-s-j}
+q_{n-j,r-s-j}+q_{r-s-j,n-j})\\ &
&-\sum_{s=0}^{r-2}\sum_{j=1}^{r-s-1}(q_{n-j,r-s-j+1}+q_{r-s-j,n-j+1})\\
&=&\sum_{s=0}^{r-1}\sum_{j=0}^{r-s-1}q^0_{n-j,r-s-j}+\sum_{j=0}^{r-1}q_{r-j,n-j}-\sum_{j=1}^{r-1}
q_{r-j,n-j+1}\hspace{5cm}(4.77)\end{eqnarray*} by (4.70) and
(4.72) for $r\in\ol{1,n-1}$, and
\begin{eqnarray*}& &
q_{r_2,r_1}^0+q_{r_2,r_1}^1+q_{r_1,r_2}+q_{r_2,r_1}\\
&=&\sum_{s=0}^{r_1-1}(q_{r_2-s,r_1-s}^0+q_{r_2-s,r_1-s}+q_{r_1-s,r_2-s})
-\sum_{s=1}^{r_1-1}(q_{r_2-s,r_1-s+1}+q_{r_1-s,r_2-s+1})
\hspace{1.2cm}(4.78)\end{eqnarray*} by (4.71) for $1\leq r_1<
r_2\leq n-1.$

Let
$$ y_{r_2,r_1}=\frac{1}{z_{r_1}z_{r_2}-1}\qquad\for\;\;1\leq
r_1\leq r_2\leq n.\eqno(4.79)$$ Set
$$w_{r_2,r_1}=\frac{y_{r_2,r_1}y_{r_2+1,r_1+1}\cdots
y_{n-1,n+r_1-r_2-1}}{y_{r_2+1,r_1}y_{r_2+2,r_1+1}\cdots
y_{n,n+r_1-r_2-1}}\qquad\for\;\;2\leq r_1\leq r_2\leq
n-1,\eqno(4.80)$$
$$ w_{k,n}=y_{n,k}\qquad\for\;\;k\in\ol{1,n-1}\eqno(4.81)$$
and
$$w_{r_1,r_2}=\frac{y_{r_2,r_1}y_{r_2+1,r_1+1}\cdots
y_{n,n+r_1-r_2}}{y_{r_2,r_1+1}y_{r_2+1,r_1+2}\cdots
y_{n-1,n+r_1-r_2}}\qquad\for\;\;1\leq r_1<r_2\leq
n-1.\eqno(4.82)$$ Based on (4.76)-(4.78), we have
\begin{eqnarray*}& &\prod_{1\leq r_1\leq r_2\leq
n}y_{r_2,r_1}^{q^0_{r_2,r_1}+q^1_{r_2,r_1}+q_{r_2,r_1}+q_{r_1,r_2}}=\left[\prod_{r_1=1}^{n-1}
\prod_{r_2=2}^nw_{r_1,r_2}^{q_{r_1,r_2}}\right]\\ &&\times
\left[\prod_{r=1}^n(y_{r,r}y_{r+1,r+1}\cdots
y_{n,n})^{q_{r,r}^0}\right]\left[\prod_{r=1}^{n-1}(y_{n,r}y_{n,r+1}\cdots
y_{n,n})^{q_{n,r}^0}\right]\prod_{1\leq r_1<r_2\leq n-1}\\ &
&(y_{r_2,r_1}y_{r_2+1,r_1+1}\cdots
y_{n-1,n+r_1-r_2-1}y_{n,n+r_1-r_1}y_{n,n+r_1-r_1+1}\cdots
y_{n,n})^{q^0_{r_2,r_1}}.\hspace{2.4cm}(4.83)\end{eqnarray*}

Next we have
\begin{eqnarray*}\hspace{2cm}& &\sum_{q_{n-1,n}=0}^{\infty}
\left(\!\!\begin{array}{c}q^0_{n-1,n-1}+q_{n-1,n-1}^1+q_{n-1,n-1}\\
q_{n-1,n}\end{array}\!\!\right)w_{n-1,n}^{q_{n-1,n}}\\
&=&(1+w_{n-1,n})^{q^0_{n-1,n-1}+q_{n-1,n-1}^1+q_{n-1,n-1}}\\
&=&(1+w_{n-1,n})^{\sum_{s=1}^{n-1}(q_{s,s}^0+q_{s,s})-\sum_{s=1}^{n-2}q_{s,s+1}}
\hspace{5.5cm}(4.84)\end{eqnarray*} by (4.76). Set
$$ w^{(2n-1)}_{r,r+1}=\frac{w_{r,r+1}}{1+w_{n-1,n}},\qquad
w_{s,s}^{(2n-1)}=(1+w_{n-1,n})w_{s,s}\eqno(4.85)$$ for
$r\in\ol{1,n-2},\;s\in\ol{1,n-1}$ and
$$
w^{(2n-1)}_{r_1,r_2}=w_{r_1,r_2}\qquad\for\;\;r_1\in\ol{1,n-1},\;r_1,r_1+1\neq
r_2\in\ol{2,n}.\eqno(4.86)$$ Suppose that we have defined
$\{w_{r_1,r_2}^{(2k+1)}\mid r_1+r_2<2k+1\}$. Set
$$\ell=\max\{1,2k-n\}.\eqno(4.87)$$
Note
\begin{eqnarray*}&&
\sum_{q_{\ell,2k-\ell},q_{\ell+1,2k-\ell-1},\cdots,q_{2k-\ell-1,\ell+1}=0}^{\infty}\;
\prod_{r=\ell}^{k-1}\left(\!\!\begin{array}{c}q^0_{2k-1-r,r}+q^1_{2k-1-r,r}+q_{2k-1-r,r}+q_{r,2k-1-r}\\
q_{r,2k-r},q_{2k-r-1,r+1}\end{array}\!\!\right)\hspace{2cm}\end{eqnarray*}
\begin{eqnarray*} & &\times
w_{r,2k-r}^{(2k+1)q_{r,2k-r}}w_{2k-r-1,r+1}^{(2k+1)q_{2k-r-1,r+1}}\\
&=&\prod_{r=\ell}^{k-1}(1+w_{r,2k-r}^{(2k+1)}+w_{2k-r-1,r+1}^{(2k+1)})
^{q^0_{2k-1-r,r}+q^1_{2k-1-r,r}+q_{2k-1-r,r}+q_{r,2k-1-r}}\\
&=&\prod_{r=\ell}^{k-1}(1+w_{r,2k-r}^{(2k+1)}+w_{2k-r-1,r+1}^{(2k+1)})^{\sum_{s=1}^r
(q^0_{2k-1-2r+s,s}+q_{2k-1-2r+s,s}+q_{s,2k-1-2r+s})}\\ &
&\times(1+w_{r,2k-r}^{(2k+1)}+w_{2k-r-1,r+1}^{(2k+1)})^{-\sum_{s=1}^r
(q_{2k-1-2r+s,s+1}+q_{s,2k-2r+s})}\hspace{3.7cm}(4.88)\end{eqnarray*}
by (4.78). Set
$$w_{2k-1-2r+s,s}^{(2k)}=w_{2k-1-2r+s,s}^{(2k+1)}(1+w_{r,2k-r}^{(2k+1)}+w_{2k-r-1,r+1}^{(2k+1)}),
\eqno(4.89)$$
$$w_{s,2k-1-2r+s}^{(2k)}=w_{s,2k-1-2r+s}^{(2k+1)}(1+w_{r,2k-r}^{(2k+1)}+w_{2k-r-1,r+1}^{(2k+1)})
\eqno(4.90)$$ for $1\leq s\leq r\leq k-1$;
$$w_{2k-1-2r+s,s+1}^{(2k)}=\frac{w_{2k-1-2r+s,s+1}^{(2k+1)}}{1+w_{r,2k-r}^{(2k+1)}
+w_{2k-r-1,r+1}^{(2k+1)}}, \eqno(4.91)$$
$$w_{s,2k-2r+s}^{(2k)}=\frac{w_{s,2k-2r+s}^{(2k+1)}}{1+w_{r,2k-r}^{(2k+1)}+
w_{2k-r-1,r+1}^{(2k+1)}} \eqno(4.92)$$ for $1\leq s\leq r\leq k-1$
and
$$w^{(2k)}_{r_1,r_2}=w^{(2k+1)}_{r_1,r_2}\;\;\mbox{for the other
pairs}\;\;(r_1,r_2)\;\;\mbox{such that}\;\;r_1+r_2\leq
2k.\eqno(4.93)$$

Let
$$\iota=\max\{1,2k-1-n\}.\eqno(4.94)$$
Moreover,
\begin{eqnarray*}&&\sum_{q_{\iota,2k-\iota-1},q_{\iota+1,2k-\iota-2},\cdots,
q_{2k-\iota-2,\iota+1}=0}
^{\infty}[\prod_{r=\iota}^{k-2}\left(\!\!\begin{array}{c}q^0_{2k-2-r,r}+q^1_{2k-2-r,r}+q_{2k-2-r,r}+q_{r,2k-2-r}
\\ q_{r,2k-r-1},q_{2k-r-2,r+1}\end{array}\!\!\right)\\ &&\times
w_{r,2k-r-1}^{(2k)q_{r,2k-r-1}}w_{2k-r-2,r+1}^{(2k)q_{2k-r-2,r+1}}]\left(\!\!\begin{array}{c}
q^0_{k-1,k-1}+q^1_{k-1,k-1}
+q_{k-1,k-1}\\ q_{k-1,k}\end{array}\!\!\right)w_{k-1,k}^{(2k)q_{k-1,k}}\\
&=&\left[\prod_{r=\iota}^{k-2}(1+w_{r,2k-r-1}^{(2k)}+w^{(2k)}_{2k-r-2,r+1})
^{q^0_{2k-2-r,r}+q^1_{2k-2-r,r}+q_{2k-2-r,r}+q_{r,2k-2-r}}\right]\\
&&\times (1+w_{k-1,k}^{(2k)})^{q^0_{k-1,k-1}+q^1_{k-1,k-1}
+q_{k-1,k-1}}\\
&=&[\prod_{r=\iota}^{k-2}(1+w_{r,2k-r-1}^{(2k)}+w^{(2k)}_{2k-r-2,r+1})
^{\sum_{s=1}^r(q^0_{2k-2-2r+s,s}+q_{2k-2-2r+s,s}+q_{s,2k-2-2r+s})}\\
& &\times(1+w_{r,2k-r-1}^{(2k)}+w^{(2k)}_{2k-r-2,r+1})
^{-\sum_{s=1}^{r-1}(q_{2k-2-2r+s,s+1}+q_{s,2k-1-2r+s})}\\ &&\times
(1+w_{k-1,k}^{(2k)})^{\sum_{r=1}^{k-1}(q^0_{r,r}+q_{r,r})-\sum_{r=1}^{k-1}q_{r,r+1}}
\hspace{7cm}(4.95)\end{eqnarray*} by (4.76) and (4.78). Set
$$w_{2k-2-2r+s,s}^{(2k-1)}=w_{2k-2-2r+s,s}^{(2k)}(1+w_{r,2k-r-1}^{(2k)}+w_{2k-r-2,r+1}^{(2k)}),
\eqno(4.96)$$
$$w_{s,2k-2-2r+s}^{(2k-1)}=w_{s,2k-2-2r+s}^{(2k)}(1+w_{r,2k-r-1}^{(2k)}+w_{2k-r-2,r+1}^{(2k)})
\eqno(4.97)$$ for $1\leq s\leq r\leq k-2$;
$$w_{2k-2-2r+s,s+1}^{(2k-1)}=\frac{w_{2k-2-2r+s,s+1}^{(2k)}}{1+w_{r,2k-r-1}^{(2k)}
+w_{2k-r-2,r+1}^{(2k)}}, \eqno(4.98)$$
$$w_{s,2k-1-2r+s}^{(2k-1)}=\frac{w_{s,2k-1-2r+s}^{(2k)}}{1+w_{r,2k-r-1}^{(2k)}+
w_{2k-r-2,r+1}^{(2k)}} \eqno(4.99)$$ for $1\leq s\leq r\leq k-2$;
$$w_{r,r}^{(2k-1)}=w_{r,r}^{(2k)}(1+w_{k-1,k}^{(2k)}),\qquad
w_{s,s+1}^{(2k-1)}=\frac{w_{s,s+1}^{(2k)}}{1+w_{k-1,k}^{(2k)}}\eqno(4.100)$$
for $r\in\ol{1,k-1},\;s\in\ol{1,k-2}$; and
$$w^{(2k-1)}_{r_1,r_2}=w^{(2k)}_{r_1,r_2}\;\;\mbox{for the other
pairs}\;\;(r_1,r_2)\;\;\mbox{such that}\;\;r_1+r_2\leq
2k-1.\eqno(4.101)$$ By induction, we have defined
$$\{w^{(k)}_{r_1,r_2}\mid
k\in\ol{1,2n-1};\;r_1\in\ol{1,n-1},\;r_2\in\ol{1,n},\;r_1+r_2\leq
k\}.\eqno(4.102)$$ For convenience, we denote
$$
w_{r_1,r_2}^{(2n)}=w_{r_1,r_2}\qquad\for\;\;r_1\in\ol{1,n-1},\;r_2\in\ol{2,n}.\eqno(4.103)$$
Based on (4.84), (4.88) and (4.95), we have
\begin{eqnarray*}&&\sum_{r_1=1}^{n-1}\;\sum_{r_2=2}^n\;\sum_{q_{r_1,r_2}=0}^{\infty}
[\prod_{1\leq l_1\leq l_2\leq n-1}
\left(\!\!\begin{array}{c}q^0_{l_2,l_1}+q_{l_2,l_1}^1+q_{l_1,l_2}+q_{l_2,l_1}
\\ q_{l_1,l_2+1},q_{l_2,l_1+1}\end{array}\!\!\right)\\ & &\times
w_{l_1,l_2+1}^{q_{l_1,l_2+1}}w_{l_2,l_1+1}^{q_{l_2,l_1+1}}]\left[\prod_{r=1}^{n-1}
\left(\!\!\begin{array}{c}q^0_{r,r}+q_{r,r}^1+q_{r,r}\\
q_{r,r+1}\end{array}\!\!\right)w_{r,r+1}^{q_{r,r+1}}\right]\\
&=&\left[\prod_{k=2}^n(1+w_{k-1,k}^{(2k)})^{\sum_{r=1}^{k-1}q_{r,r}^0}\right]\left[\prod_{k=2}
^{n-1}\prod_{r=\ell}^{k-1}(1+w_{r,2k-r}^{(2k+1)}+w_{2k-r-1,r+1}^{(2k+1)})^{\sum_{s=1}^r
q_{2k-1-2r+s,s}^0}\right]\\ &&\times
\left[\prod_{k=3}^{n-1}\prod_{r=\iota}^{k-2}(1+w_{r,2k-r-1}^{(2k)}+w_{2k-r-2,r+1}^{(2k)})
^{\sum_{s=1}^rq_{2k-2-2r+s,s}^0}\right]\\ &=&[\prod_{1\leq
r_1<r_2\leq n-1;r_2-r_1\;\mbox{\small is
odd}}\;(\prod_{k=(r_2-r_1+3)/2}^{n-1}(1+w^{(2k+1)}_{(2k+r_1-r_2-1)/2,(2k-r_1+r_2+1)/2}\\
&&+w^{(2k+1)}_{(2k-r_1+r_2-1)/2,(2k+r_1-r_2+1)/2})^{q^0_{r_2,r_1}})]
[\prod_{1\leq r_1<r_2\leq n-1;r_2-r_1\;\mbox{\small is
even}}\hspace{4cm}\end{eqnarray*}
\begin{eqnarray*}& &
(\prod_{k=(r_2-r_1+2)/2}^{n-1}(1+w^{(2k+2)}_{(2k+r_1-r_2)/2,(2k-r_1+r_2+2)/2}+w^{(2k+2)}
_{(2k-r_1+r_2)/2,(2k+r_1-r_2+2)/2})^{q^0_{r_2,r_1}})]\\ & &\times
\left[\prod_{r=1}^{n-1}\left(\prod_{k=r}^{n-1}(1+w_{k,k+1}^{(2k)})\right)^{q_{r,r}^0}\right].
\hspace{8.4cm}(4.104)\end{eqnarray*}

According to (4.74), (4.83) and (4.104), we define
$$ \xi_{n,i}^C=(-1)^{n+i+1}\frac{y_{n,i}y_{n,i+1}\cdots
y_{n,n}}{2^{1+\dlt_{n,i}}}\qquad\for\;\;1\leq i\leq
n,\eqno(4.105)$$
$$\xi_{i,i}^C=-\frac{y_{i,i}\cdots
y_{n,n}\prod_{k=i}^{n-1}(1+w_{k,k+1}^{(2k)})}{4}\eqno(4.106)$$ for
$1\leq i\leq n-1$,
\begin{eqnarray*}&
&\xi_{r_2,r_1}^C=\frac{(-1)^{r_1+r_2+1}}{2}y_{r_2,r_1}y_{r_2+1,r_1+1}\cdots
y_{n-1,n+r_1-r_2-1}y_{n,n+r_1-r_2}y_{n,n+r_1-r_2+1}\cdots
y_{n,n}\\ &
&\prod_{k=(r_2-r_1+3)/2}^{n-1}(1+w^{(2k+1)}_{(2k+r_1-r_2-1)/2,(2k-r_1+r_2+1)/2}
+w^{(2k+1)}_{(2k-r_1+r_2-1)/2,(2k+r_1-r_2+1)/2})\hspace{0.5cm}(4.107)\end{eqnarray*}
for $1\leq r_1<r_2\leq n-1$ such that $r_2-r_1$ is odd, and
\begin{eqnarray*}&
&\xi_{r_2,r_1}^C=\frac{(-1)^{r_1+r_2+1}}{2}y_{r_2,r_1}y_{r_2+1,r_1+1}\cdots
y_{n-1,n+r_1-r_2-1}y_{n,n+r_1-r_2}y_{n,n+r_1-r_2+1}\cdots
y_{n,n}\\ &
&\prod_{k=(r_2-r_1+2)/2}^{n-1}(1+w^{(2k+2)}_{(2k+r_1-r_2)/2,(2k-r_1+r_2+2)/2}+w^{(2k+2)}
_{(2k-r_1+r_2)/2,(2k+r_1-r_2+2)/2})\hspace{1.2cm}(4.108)\end{eqnarray*}
for $1\leq r_1<r_2\leq n-1$ such that $r_2-r_1$ is even.

 For
$\al\in\G_C$, we set
$$\al^c_i=\sum_{r=1}^i\al_{n+r,i}+\sum_{s=i}^n\al_{n+s,i}\qquad\for\;\;i\in\ol{1,n}\eqno(4.109)$$
and
$$\al^c=\sum_{1\leq r_1\leq r_2\leq n}\al_{n+r_2,r_1}.\eqno(4.110)$$

By (3.40), (4.74), (4.75), (4.79)-(4.83) and (4.104)-(4.110),  we
have
\begin{eqnarray*}& &E_C(z_1,z_2,...,z_n)\\ &=&\frac{z_n^{n+\lmd_n+1}\prod_{i=1}^{n-1}
z_i^{2n-i+1/2}}{\left[\prod_{1\leq r_1<r_2\leq
n}(z_{r_1}-z_{r_2})\right]\left[\prod_{1\leq s_1 \leq s_2\leq
n}(z_{s_1}z_{s_2}-1)\right]}\sum_{\gm\in\G_C}\frac{1}{\gm!\la\lmd_n\ra_{_{\gm^c}}}
\\& &\times\prod_{r=1}^n\la \gm_{\ol{r}}-\gm_{\ul{r}}+\gm_r^c-1/2\ra_{_{\gm_{\ol{r}}+\gm_r^c}}
\xi^\gm,\hspace{7.6cm}(4.111)\end{eqnarray*} where
$$\gm !=\left[\prod_{1\leq r_1<r_2\leq n}\gm_{r_2,r_1}!\right]\left[\prod_{1\leq s_1\leq s_2
\leq n}\gm_{n+s_2,s_1}!\right]\eqno(4.112)$$ and
$$\xi^\gm=
\left[\prod_{1\leq r_1<r_2\leq
n}(\xi^A_{r_2,r_1})^{\gm_{r_2,r_1}}\right] \left[\prod_{1\leq
s_1\leq s_2\leq
n}(\xi^C_{n+s_2,s_1})^{\gm_{n+s_2,s_1}}\right].\eqno(4.112)$$
Recall the notations in (1.16) and (3.38). We define our
hypergeometric function of type C by
$${\cal X}_C(\tau_1,...,\tau_n;\vt)\{z_{r_2,r_1},z_{n+s_2,s_1}\}=
\sum_{\al\in\G_C}\frac{\prod_{r=1}^n(\tau_r-\al_{\ul{r}})_{_{\al_{\ol{r}}+\al^c_r}}}
{\al!(\vt)_{\al^c}}z^\al.\eqno(4.114)$$
 By (2.18), (4.110) and (4.114), we obtain the following main theorem in this section:
\psp

{\bf Theorem 4.1}. {\it The Etingof trace function}
\begin{eqnarray*}E_C(z_1,z_2,...,z_n)&=&
\frac{z_n^{n+\lmd_n+1}\prod_{i=1}^{n-1}z_i^{2n-i+1/2}}
{\left[\prod_{1\leq r_1<r_2\leq
n}(z_{r_1}-z_{r_2})\right]\left[\prod_{1\leq s_1\leq s_2\leq
n}(z_{s_1}z_{s_2}-1)\right]}\\ & &\times {\cal
X}_C(1/2,...,1/2;-\lmd_n)\{\xi^A_{r_2,r_1},\xi^C_{n+s_2,s_1}\}.
\hspace{4cm}(4.115)\end{eqnarray*}

\section{Path Hypergeometric Functions}

In this section, we find the differential properties, integral
representations and differential equations  for the hypergeometric
functions of type A in (3.38),  and the differential properties
and differential equations  for the hypergeometric functions of
type C in (4.114). Moreover, we define our hypergeometric
functions of type B and D analogously as those of type C.

For two positive integers $k_1$ and $k_2$ such that $k_1<k_2$, a {\it path} from  $k_1$ to
 $k_2$ is a sequence $(m_0,m_1....,m_r)$ of positive integers such that
$$k_1=m_0<m_1<m_2<\cdots <m_{r-1}<m_r=k_2.\eqno(5.1)$$
One can imagine a path from $k_1$ to $k_2$ is a way of a super man going from $k_1$th floor to
$k_2$th floor through a stairway. Let
$${\cal P}_{k_1}^{k_2}=\mbox{the set of all paths from}\;k_1\;\mbox{to}\;k_2.\eqno(5.2)$$
The {\it path polynomial} from $k_1$ to $k_2$ is defined as
$$P_{[k_1,k_2]}=\sum_{(m_0,m_1,...,m_r)\in {\cal P}_{k_1}^{k_2}}(-1)^rz_{m_1,m_0}z_{m_2,m_1}
\cdots z_{m_{r-1},m_{r-2}}z_{m_r,m_{r-1}}.\eqno(5.3)$$
Moreover, we set
$$P_{[k,k]}=1\qquad\for\;\;0<k\in\mbb{N}.\eqno(5.4)$$
For convenience, we simply denote
$${\cal X}_A={\cal X}_A(\tau_1,..,\tau_n;\vt)\{z_{j,k}\},\eqno(5.5)$$
$${\cal X}_A[i,j]={\cal X}_A(\tau_1,...,\tau_i+1,...,\tau_j-1,...\tau_n;\vt)\{z_{r_2,r_1}\}
\eqno(5.6)$$
obtained from ${\cal X}_A$ by changing $\tau_i$ to $\tau_i+1$ and $\tau_j$ to $\tau_j-1$ for
$1\leq i<j\leq n-1$ and
$${\cal X}_A[k,n]={\cal X}_A(\tau_1,..,\tau_k+1,...,\tau_n+1;\vt+1)\{z_{r_2,r_1}\}\eqno(5.7)$$
obtained from ${\cal X}_A$ by changing $\tau_i$ to $\tau_i+1$, $\tau_n$ to $\tau_n+1$ and
$\vt$ to $\vt+1$ for $k\in\ol{1,n-1}$.
\psp

{\bf Theorem 5.1}. {\it For $1\leq r_1<r_2\leq n-1$ and $r\in \ol{1,n-1}$, we have
$$\ptl_{z_{r_2,r_1}}({\cal X}_A)=\sum_{s=1}^{r_1}\tau_sP_{[s,r_1]}{\cal X}_A[s,r_2],\eqno(5.8)$$
$$\ptl_{z_{n,r}}({\cal X}_A)=\frac{\tau_n}{\vt}\sum_{s=1}^r\tau_sP_{[s,r]}{\cal X}_A[s,n].
\eqno(5.9)$$}

{\it Proof}. Note that
\begin{eqnarray*}& & \ptl_{z_{r_2,r_1}}({\cal X}_A)\\ &=& \sum_{\be\in\G_A}
\frac{\left[\prod_{s=1}^{n-1}(\tau_s-\be_{\ul{s}})_{_{\be_{\ol{s}}}}\right]
(\tau_n)_{_{\be_{\ul{n}}}}}{\be!(\vt)_{_{\be_{\ul{n}}}}}\be_{r_2,r_1}z^{\be-\es_{r_2,r_1}}\\
 &=&\sum_{\be\in\G_A}\frac{\left[(\tau_{r_1}-\be_{\ul{r_1}})(\tau_{r_1}+1
 -\be_{\ul{r_1}})_{\be_{\ol{r_1}}}(\tau_{r_2}-1-\be_{\ul{r_2}})_{\be_{\ol{r_2}}}
 \prod_{s\neq r_1,r_2}(\tau_s-\be_{\ul{s}})_{_{\be_{\ol{s}}}}\right](\tau_n)_{_{\be_{\ul{n}}}}}
 {\be!(\vt)_{_{\be_{\ul{n}}}}}z^\be\\ &=&(\tau_{r_1}-\sum_{s=1}^{r_1-1}z_{r_1,s}\ptl_{z_{r_1,s}})
 ({\cal X}_A[r_1,r_2]).\hspace{8.2cm}(5.10)\end{eqnarray*}
In particular,
$$ \ptl_{z_{r,1}}({\cal X}_A)=\tau_1{\cal X}_A[1,r]\qquad\for\;\;r\in\ol{1,n-1}.\eqno(5.11)$$
By (5.10), (5.11) and  induction, we get (5.8). Moreover,
\begin{eqnarray*}\qquad \ptl_{z_{n,r}}({\cal X}_A)&=& \sum_{\be\in\G_A}\frac{\left[
\prod_{s=1}^{n-1}(\tau_s-\be_{\ul{s}})_{_{\be_{\ol{s}}}}\right](\tau_n)_{_{\be_{\ul{n}}}}}
{\be!(\vt)_{_{\be_{\ul{n}}}}}\be_{n,r}z^{\be-\es_{n,r}}\\ &=&
\sum_{\be\in\G_A}\frac{\left[(\tau_r-\be_{\ul{r}})(\tau_r+1-\be_{\ul{r}})_{\be_{\ol{r}}}
\prod_{s\neq
r}(\tau_s-\be_{\ul{s}})_{_{\be_{\ol{s}}}}\right]\tau_n(\tau_n+1)_{_{\be_{\ul{n}}}}}
{\be!\vt(\vt+1)_{_{\be_{\ul{n}}}}}z^\be\\
&=&\frac{\tau_n}{\vt}(\tau_r-\sum_{s=1}^rz_{r,s}
\ptl_{z_{r,s}})({\cal
X}_A[r,n]),\hspace{5.9cm}(5.12)\end{eqnarray*} by which we obtain
(5.9).$\qquad\Box$ \psp

For any $z\in \mbb{C}\setminus (-\mbb{N})$, the gamma function
$$\G(z)=\left[z e^{cz}\prod_{m=1}^{\infty}\left\{\left(1+\frac{z}{m}\right)e^{-z/m}\right\}
\right]^{-1},\eqno(5.13)$$ where $c$ is Euler's constant given by
$$c=\lim_{m\rta\infty}\left(\sum_{k=1}^m\frac{1}{k}-\ln m\right).\eqno(5.14)$$
When $\mbox{Re}\:z>0$, we have
$$\G(z)=\int_0^{\infty}t^{z-1}e^{-t}dt.\eqno(5.15)$$
Now we have the following theorem of integral representation:
\psp

{\bf Theorem 5.2}. {\it Suppose $\mbox{\it Re}\:\tau_n>0$ and $\mbox{\it Re}\:(\vt-\tau_n)>0$.
 We have}
$${\cal X}_A=\frac{\G(\vt)}{\G(\vt-\tau_n)\G(\tau_n)}\int_0^1\left[\prod_{r=1}^{n-1}
(\sum_{s=r}^{n-1}P_{[r,s]}+tP_{[r,n]})^{-\tau_r}\right]t^{\tau_n-1}(1-t)^{\vt-\tau_n-1}dt.
\eqno(5.16)$$

{\it Proof}.  For $\mu_1,\mu_2\in\mbb{C}$ with $\mbox{Re}\:\mu_1,\mbox{Re}\:\mu_1>0$, we have the
following Euler integral
$$\int_0^1t^{\mu_1-1}(1-t)^{\mu_2-1}dt=\frac{\G(\mu_1)\G(\mu_2)}{\G(\mu_1+\mu_2)}.\eqno(5.17)$$
Moreover,
$$(\mu)_m=\frac{\G(\mu+m)}{\G(\mu)}\qquad\for\;\;m\in\mbb{N},\;\mu\in\mbb{C}.\eqno(5.18)$$
Recall
 $${\cal X}_A=\sum_{\be\in\G_A}\frac{\left[\prod_{s=1}^{n-1}(\tau_s-\be_{\ul{s}})_{_{\be_{\ol{s}}}}
 \right](\tau_n)_{_{\be_{\ul{n}}}}}{\be!(\vt)_{_{\be_{\ul{n}}}}}z^\be.\eqno(5.19)$$
We have
\begin{eqnarray*}\qquad & &\frac{(\tau_n)_{_{\be_{\ul{n}}}}}{(\vt)_{_{\be_{\ul{n}}}}}=
\frac{\G(\vt)}{\G(\vt-\tau_n)\G(\tau_n)}\frac{\G(\tau_n+\be_{\ul{n}})\G(\vt-\tau_n)}{\G(\vt
+\be_{\ul{n}})}\\
&=&\frac{\G(\vt)}{\G(\vt-\tau_n)\G(\tau_n)}\int_0^1t^{\tau_n+\be_{\ul{n}}-1}
(1-t)^{\vt-\tau_n-1}dt\hspace{5.5cm}(5.20)\end{eqnarray*} by
(5.17) and (5.19).

Denote
$$\G_{(s)}=\sum_{1\leq r_1<r_2\leq s}\mbb{N}\es_{r_2,r_1},\qquad\be_{\ol{r}}^s=\sum_{p=r+1}^s
\be_{p,r}\qquad\for\;\;\be\in\G_A,\;1\leq r<s\leq
n-1.\eqno(5.21)$$ Note that
\begin{eqnarray*}& &\sum_{\be\in\G_A}\frac{\prod_{s=1}^{n-1}(\tau_s-\be_{\ul{s}})_{_{\be_{\ol{s}}}}}
{\be!}z^\be
\\&=&\sum_{\be\in\G_{(n-1)}}\frac{\prod_{s=1}^{n-2}(\tau_s-\be_{\ul{s}})
_{_{\be^{n-1}_{\ol{s}}}}}{\be!}z^\be\sum_{\be_{n,1},...,\be_{n,n-1}=0}^n\frac{(\tau_{n-1}
-\be_{\ul{n-1}})_{\be_{n,n-1}}z_{n,n-1}^{\be_{n,n-1}}}{\be_{n,n-1}!}\\
& &\times\prod_{r=1}^{n-2}
\frac{(\tau_r-\be_{\ul{r}}+\be^{n-1}_{\ol{r}})_{\be_{n,r}}z_{n,r}^{\be_{n,r}}}{\be_{n,r}!}
\hspace{12cm}\end{eqnarray*}
 \begin{eqnarray*} &=&
\sum_{\be\in\G_{(n-1)}}\frac{\prod_{s=1}^{n-2}(\tau_s-\be_{\ul{s}})_{_{\be^{n-1}_{\ol{s}}}}}
{\be!}z^\be(1-z_{n,n-1})^{\be_{\ul{n-1}}-\tau_{n-1}}\prod_{r=1}^{n-2}(1-z_{n,r})^{\be_{\ul{r}}
-\be_{\ol{r}}^{n-1}-\tau_r}\\
&=& \left[\prod_{s=1}^{n-1}(1-z_{n,s})^{-\tau_s}\right]
\sum_{\be\in\G_{(n-1)}}\frac{\prod_{s=1}^{n-2}(\tau_s-\be_{\ul{s}})_{_{\be^{n-1}_{\ol{s}}}}}
{\be!}\\
& &\times\prod_{1\leq r_1<r_2\leq
n-1}\left[\frac{(1-z_{n.r_2})z_{r_2,r_1}}{1-z_{n,r_1}}\right]^{\be_{r_2,r_1}}.
\hspace{7.3cm}(5.22)\end{eqnarray*} Observe that
$$1-\frac{(1-z_{n.n-1})z_{n-1,r}}{1-z_{n,r}}=\frac{1-z_{n-1,r}-z_{n,r}+z_{n.n-1}z_{n-1,r}}{1
-z_{n,r}}\qquad\for\;\;r\in\ol{1,n-1}.\eqno(5.23)$$ By (5.22) and
(5.23), we have
\begin{eqnarray*}& &\sum_{\be\in\G_A}\frac{\prod_{s=1}^{n-1}(\tau_s-\be_{\ul{s}})_{_{\be_{\ol{s}}}}}
{\be!}z^\be=(1-z_{n,n-1})^{-\tau_{n-1}}(1-z_{n-1,n-2}+P_{[n-2,n]})^{-\tau_{n-2}}\\
& &\times
\left[\prod_{s=1}^{n-2}(1-z_{n-1,s}-z_{n,s}+z_{n.n-1}z_{n-1,s})^{-\tau_s}\right]
\sum_{\be\in\G_{(n-2)}}\frac{\prod_{s=1}^{n-3}(\tau_s-\be_{\ul{s}})_{_{\be^{n-2}_{\ol{s}}}}}
{\be !}\\ & &\times \prod_{1\leq r_1<r_2\leq
n-2}\left[\frac{(1-z_{n-1,r_2}-z_{n,r_2}
+z_{n.n-1}z_{n-1,r_2})z_{r_2,r_1}}{1-z_{n-1,r_1}-z_{n,r_1}+z_{n.n-1}z_{n-1,r_1}}
\right]^{\be_{r_2,r_1}}.\hspace{3.3cm}(5.24)\end{eqnarray*} Now
\begin{eqnarray*} & &1-\frac{(1-z_{n-1,n-2}-z_{n,r-2}+z_{n.n-1}z_{n-1,n-2})z_{n-2,r}}
{1-z_{n-1,r}-z_{n,r}+z_{n.n-1}z_{n-1,r}}\\
&=&\frac{1}{1-z_{n-1,r}-z_{n,r}+z_{n.n-1}z_{n-1,r}}
[1-z_{n-2,r}-z_{n-1,r}+z_{n-1,n-2}z_{n-2,r}\\ &
&-z_{n,r}+z_{n.n-1}z_{n-1,r}+z_{n,r-2}z_{n-2,r}
-z_{n.n-1}z_{n-1,n-2}z_{n-2,r}].\hspace{3.7cm}(5.25)\end{eqnarray*}
By induction, we can prove that
$$\sum_{\be\in\G_A}\frac{\prod_{s=1}^{n-1}(\tau_s-\be_{\ul{s}})_{_{\be_{\ol{s}}}}}{\be!}z^\be
=\prod_{r=1}^{n-1}(\sum_{s=r}^nP_{[r,s]})^{-\tau_r}.\eqno(5.26)$$
Thus
$$\sum_{\be\in\G_A}\frac{\prod_{s=1}^{n-1}(\tau_s-\be_{\ul{s}})_{_{\be_{\ol{s}}}}}{\be!}
z^\be
t^{\be_{\ul{n}}}=\prod_{r=1}^{n-1}(\sum_{s=r}^{n-1}P_{[r,s]}+tP_{[r,n]})^{-\tau_r}.
\eqno(5.27)$$ Hence we obtain (5.16) by (5.20) and
(5.27).$\qquad\Box$ \psp

{\bf Remark 5.3}. According to [M],
$$\left(\begin{array}{ccccc}1&0&0&\cdots&0\\ z_{2,1}&
1&0&\cdots&0\\ z_{3,1}&z_{3,2}&1&\ddots&\vdots\\
\vdots&\vdots&\ddots&\ddots&0\\ z_{n,1}&z_{n,2}&\cdots&z_{n,n-1}&
1\end{array}\right)^{-1}=\left(\begin{array}{ccccc}1&0&0&\cdots&0\\
P_{1,2}&
1&0&\cdots&0\\ P_{1,3}&P_{2,3}&1&\ddots&\vdots\\
\vdots&\vdots&\ddots&\ddots&0\\ P_{1,n}&P_{2,n}&\cdots&P_{n-1,n}&
1\end{array}\right).\eqno(5.28)$$

Denote
$${\cal D}_{\ul{i}}^A=\sum_{r=1}^{i-1}z_{i,r}\ptl_{z_{i,r}},\qquad {\cal D}_{\ol{i}}^A=
\sum_{s=i+1}^nz_{s,i}\ptl_{z_{s,i}}\qquad\for\;\;i\in\ol{1,n}.\eqno(5.29)$$
\pse

{\bf Theorem 5.4}. {\it We have:
$$(\tau_{r_2}-1-{\cal D}_{\ul{r_2}}^A+{\cal D}_{\ol{r_2}}^A)\ptl_{z_{r_2,r_1}}({\cal X}_A)
=(\tau_{r_2}-1-{\cal D}_{\ul{r_2}}^A)(\tau_{r_1}-{\cal
D}_{\ul{r_1}}^A+{\cal D}_{\ol{r_1}}^A) ({\cal X}_A)\eqno(5.30)$$
for $1\leq r_1<r_2\leq n-1$ and
$$(\vt+{\cal D}_{\ul{n}}^A)\ptl_{z_{n,r}}({\cal X}_A)=(\tau_n+{\cal D}_{\ul{n}}^A)
(\tau_r-{\cal D}_{\ul{r}}^A+{\cal D}_{\ol{r}}^A)({\cal
X}_A)\eqno(5.31)$$ for $r\in\ol{1,n-1}$.} \pse

{\it Proof}. Let
$$a_\be=\frac{\left[\prod_{s=1}^{n-1}(\tau_s-\be_{\ul{s}})_{_{\be_{\ol{s}}}}\right]
(\tau_n)_{_{\be_{\ul{n}}}}}{\be!(\vt)_{_{\be_{\ul{n}}}}}\qquad\for\;\;\be\in\G_A.\eqno(5.32)$$
Then
$$\frac{a_{\be+\es_{r_2,r_1}}}{a_\be}=\frac{(\tau_{r_2}-\be_{\ul{r_2}}-1)(\tau_{r_1}-
\be_{\ul{r_1}}+\be_{\ol{r_1}})}{(\tau_{r_2}-1-\be_{\ul{r_2}}+\be_{\ol{r_2}})(\be_{r_2,r_1}+1)},
\eqno(5.33)$$ equivalently,
$$(\tau_{r_2}-1-\be_{\ul{r_2}}+\be_{\ol{r_2}})(\be_{r_2,r_1}+1)a_{\be+\es_{r_2,r_1}}=(\tau_{r_2}
-\be_{\ul{r_2}}-1)(\tau_{r_1}-\be_{\ul{r_1}}+\be_{\ol{r_1}})a_\be.\eqno(2.34)$$
Thus
$$\sum_{\be\in\G_A}(\tau_{r_2}-1-\be_{\ul{r_2}}+\be_{\ol{r_2}})(\be_{r_2,r_1}+1)
a_{\be+\es_{r_2,r_1}}z^\be=\sum_{\be\in\G_A}(\tau_{r_2}-\be_{\ul{r_2}}-1)(\tau_{r_1}
-\be_{\ul{r_1}}+\be_{\ol{r_1}})a_\be z^\be,\eqno(2.35)$$ which is
equivalent to (5.30). Similarly, (5.31) follows from
$$\frac{a_{\be+\es_{n,r}}}{a_\be}=\frac{(\tau_n-\be_{\ul{n}})(\tau_r-\be_{\ul{r}}+\be_{\ol{r}})}
{(\vt+\be_{\ul{n}})(\be_{n,r}+1)}.\qquad\Box\eqno(5.36)$$ \psp

Next we study our hypergeometric functions of type C. For convenience, we simply denote
$${\cal X}_C={\cal X}_C(\tau_1,...,\tau_n;\vt)\{z_{r_2,r_1},z_{n+s_2,s_1}\}\eqno(5.37)$$
(cf. (4.114)). First we have
\begin{eqnarray*}& &\ptl_{z_{r_2,r_1}}({\cal X}_C)\\ &=&\sum_{\al\in\G_C}\frac{\prod_{r=1}^n
(\tau_r-\al_{\ul{r}})_{_{\al_{\ol{r}}+\al^c_r}}}{\al!(\vt)_{\al^c}}\al_{r_2,r_1}z^{\al
-\es_{r_2,r_1}}\hspace{9cm}\end{eqnarray*}
\begin{eqnarray*}
 &=& \sum_{\al\in\G_C}\frac{\prod_{r\neq
r_1,r_2}(\tau_r-\al_{\ul{r}})_{_{\al_{\ol{r}}+\al^c_r}}}{\al!(\vt)_{\al^c}}(\tau_{r_1}
-\al_{\ul{r_1}})(\tau_{r_1}+1-\al_{\ul{r_1}})_{_{\al_{\ol{r_1}}+\al^c_{r_1}}}\\
&
&\times(\tau_{r_2}-1-\al_{\ul{r_2}})_{_{\al_{\ol{r_2}}+\al^c_{r_2}}}z^\al.
\hspace{9.3cm}(5.38)\end{eqnarray*} or $1\leq r_1<r_2\leq n$.
Moreover,
\begin{eqnarray*}& &\ptl_{z_{n+r_2,r_1}}({\cal X}_C)\\ &=&\sum_{\al\in\G_C}
\frac{\prod_{r=1}^{n-1}(\tau_r-\al_{\ul{r}})_{_{\al_{\ol{r}}+\al^c_r}}}{\al!(\vt)_{\al^c}}
\al_{n+r_2,r_1}z^{\al-\es_{n+r_2,r_1}}\\ &=&
\frac{1}{\vt}\sum_{\al\in\G_C} \frac{\prod_{r\neq
r_1,r_2}(\tau_r-\al_{\ul{r}})_{_{\al_{\ol{r}}+\al^c_r}}}
{\al!(\vt+1)_{\al^c}}(\tau_{r_1}-\al_{\ul{r_1}})(\tau_{r_2}-\al_{\ul{r_2}})(\tau_{r_1}
+1-\al_{\ul{r_1}})_{_{\al_{\ol{r_1}}+\al^c_{r_1}}}\\ & &\times
(\tau_{r_2}+1-
\al_{\ul{r_2}})_{_{\al_{\ol{r_2}}+\al^c_{r_2}}}z^{\al}\hspace{9.3cm}(5.39)\end{eqnarray*}
for $1\leq r_1<r_2\leq n$ and
\begin{eqnarray*}\ptl_{z_{n+s,s}}({\cal X}_C)&=&\sum_{\al\in\G_C}\frac{\prod_{r=1}^n(\tau_r
-\al_{\ul{r}})_{_{\al_{\ol{r}}+\al^c_r}}}{\al!(\vt)_{\al^c}}\al_{n+s,s}z^{\al-\es_{n+s,s}}\\
&=&\frac{1}{\vt}\sum_{\al\in\G_C}\frac{\prod_{r\neq s}(\tau_r-\al_{\ul{r}})_{_{\al_{\ol{r}}
+\al^c_r}}}{\al!(\vt+1)_{\al^c}}(\tau_{s}-\al_{\ul{s}})(\tau_{s}+1-\al_{\ul{s}})\\
& &\times
(\tau_s+2-\al_{\ul{s}})_{_{\al_{\ol{s}}+\al^c_{s}}}z^{\al}\hspace{8cm}(5.40)\end{eqnarray*}
for $s\in\ol{1,n}$.

Expressions (5.38)-(5.40) motivate us to define
$${\cal X}_C[i_1,i_2]={\cal X}_C(\tau_1,...,\tau_{i_1}+1,...,\tau_{i_2}-1,...\tau_n;\vt)
\{z_{r_2,r_1},z_{n+s_2,s_1}\}\eqno(5.41)$$ obtained from ${\cal
X}_C$  by changing $\tau_{i_1}$ to $\tau_{i_1}+1$ and $\tau_{i_2}$
to $\tau_{i_2}-1$ for $1\leq i_1<i_2\leq n$,
$${\cal X}_C[j_1,j_2;1]={\cal X}_C(\tau_1,...,\tau_{j_1}+1,...,\tau_{j_2}+1,...\tau_n;\vt+1)
\{z_{r_2,r_1},z_{n+s_2,s_1}\}\eqno(5.42)$$ obtained from ${\cal
X}_C$  by changing $\tau_{j_1}$ to $\tau_{j_1}+1$, $\tau_{j_2}$ to
$\tau_{j_2}+1$ and $\vt$ to $\vt+1$ for $1\leq j_1<j_2\leq n$ and
$${\cal X}_C[k(2)]={\cal X}_C(\tau_1,...,\tau_k+2,...,\tau_n;\vt+1)\{z_{r_2,r_1},z_{n+s_2,s_1}\}
\eqno(5.43)$$ obtained from ${\cal X}_C$  by changing $\tau_k$ to
$\tau_k+2$ and $\vt$ to $\vt+1$ for $k\in\ol{1,n}$. Now (5.38) can
be written as
$$\ptl_{z_{r_2,r_1}}({\cal X}_C)=(\tau_{r_1}-\sum_{s=1}^{r_1-1}z_{r_1,s}\ptl_{z_{r_1,s}})
({\cal X}_C[i_1,i_2])\eqno(5.44)$$ and (5.39) becomes
$$\ptl_{z_{n+r_2,r_1}}({\cal X}_C)=\frac{1}{\vt}(\tau_{r_1}-\sum_{s_1=1}^{r_1-1}z_{r_1,s_1}
\ptl_{z_{r_1,s_1}})(\tau_{r_2}-\sum_{s_2=1}^{r_2-1}z_{r_2,s_2}\ptl_{z_{r_2,s_2}})
({\cal X}_C[i_1,i_2;1])\eqno(5.45)$$ for $1\leq r_1<r_2\leq n.$
Moreover, (5.40) is equivalent to
$$\ptl_{z_{n+s,s}}({\cal X}_C)=\frac{1}{\vt}(\tau_s-\sum_{r=1}^{s-1}z_{s,r}\ptl_{z_{s,r}})
(\tau_s+1-\sum_{r=1}^{s-1}z_{s,r}\ptl_{z_{s,r}})({\cal
X}_C[s(2)])\eqno(5.46)$$ for $s\in\ol{1,n}$. By (5.44)-(5.46) and
induction, we obtain: \psp

{\bf Theorem 5.5}. {\it The following equations hold for ${\cal
X}_C$:
$$\ptl_{z_{r_2,r_1}}({\cal X}_C)=\sum_{s=1}^{r_1}\tau_sP_{[s,r_1]}{\cal X}_C[s,r_2],
\eqno(5.47)$$
\begin{eqnarray*}\ptl_{z_{n+r_2,r_1}}({\cal X}_C)&=&\frac{1}{\vt}[\sum_{i=1}^{r_1}\tau_i^2
P_{[i,r_1]}P_{[i,r_2]}{\cal
X}_C[i(2)]+\sum_{s_1=1}^{r_1}\sum_{s_2=r_1+1}^{r_2}\tau_{s_1}
\tau_{s_2}P_{[s_1,r_1]}P_{[s_2,r_2]}{\cal X}_C[s_1,s_2;1]\\ &
&+\sum_{1\leq s_1<s_2\leq r_1}
\tau_{s_1}\tau_{s_2}(P_{[s_1,r_1]}P_{[s_2,r_2]}+P_{[s_2,r_1]}P_{[s_1,r_2]})
{\cal X}_C[s_1,s_2;1]]\hspace{1cm}(5.48)\end{eqnarray*} for $1\leq
r_1<r_2\leq n$ and
\begin{eqnarray*}\ptl_{z_{n+s,s}}({\cal X}_C)&=&\frac{1}{\vt}[\sum_{i=1}^s\tau_i^2P_{[i,s]}^2
{\cal X}_C[i(2)]+\tau_s{\cal X}_C[s(2)]+
\sum_{i=1}^{s-1}\tau_iP_{[i,s]}{\cal X}_C[i,s;1]\\ & &
+2\sum_{1\leq r_1<r_2\leq
s}\tau_{r_1}\tau_{r_2}P_{[r_1,s]}P_{[r_2,s]}{\cal X}_C[r_1,r_2;1]]
\hspace{4.5cm}(5.49)\end{eqnarray*} for $s\in\ol{1,n}$.} \psp

Up to this stage, we have not found a nice integral representation for ${\cal X}_C$. In fact,
there is no simple integral formula of Euler type with an elementary integrand for Lauricella
third multiple hypergeometric function (e.g., cf. [Eh]). It might also be the case for our
hypergeometric function ${\cal X}_C$.

By (4.114), we set
$$c_\al=\frac{\prod_{r=1}^n(\tau_r-\al_{\ul{r}})_{_{\al_{\ol{r}}+\al^c_r}}}{\al!(\vt)_{\al^c}}
z^\al\qquad\for\;\;\al\in\G_C.\eqno(5.50)$$ Note that
$$\frac{c_{\al+\es_{r_2,r_1}}}{c_\al}=\frac{(\tau_{r_1}-\al_{_{\ul{r_1}}}+\al_{_{\ol{r_1}}}
+\al^c_{r_1})(\tau_{r_2}-1-\al_{_{\ul{r_2}}})}{(\al_{r_2,r_1}+1)(\tau_{r_2}-1-\al_{_{\ul{r_2}}}
+\al_{_{\ol{r_2}}}+\al^c_{r_2})},\eqno(5.51)$$
$$\frac{c_{\al+\es_{n+r_2,r_1}}}{c_\al}=\frac{(\tau_{r_1}-\al_{_{\ul{r_1}}}+\al_{_{\ol{r_1}}}
+\al^c_{r_1})(\tau_{r_2}-\al_{_{\ul{r_2}}}+\al_{_{\ol{r_2}}}+\al^c_{r_2})}{(\al_{n+r_2,r_1}+1)
(\vt+\al^c)}\eqno(5.52)$$ for $1\leq r_1<r_2\leq n$ and
$$\frac{c_{\al+\es_{n+s,s}}}{c_\al}=\frac{(\tau_s-\al_{_{\ul{s}}}+\al_{_{\ol{s}}}+\al^c_{s})
(\tau_s+1-\al_{_{\ul{s}}}+\al_{_{\ol{s}}}+\al^c_{s})}{(\al_{n+s,s}+1)(\vt+\al^c)}\eqno(5.53)$$
for $s\in\ol{1,n}$. Let
$${\cal D}^C_r=\sum_{i=1}^rz_{n+r,i}\ptl_{z_{n+r,i}}+\sum_{s=r}^nz_{n+s,r}\ptl_{z_{n+s,r}}
\qquad\for\;\;r\in\ol{1,n}\eqno(5.54)$$ and
$${\cal D}^C=\sum_{1\leq r_1\leq r_2\leq n}z_{n+r_2,r_1}\ptl_{z_{n+r_2,r_1}}.\eqno(5.55)$$
By the proof of Theorem 5.4, we have: \psp

{\bf Theorem 5.6}. {\it The function ${\cal X}_C$ satisfies:
$$(\tau_{r_2}-1-{\cal D}_{\ul{r_2}}^A+{\cal D}_{\ol{r_2}}^A+{\cal D}^C_{r_2})\ptl_{z_{r_2,r_1}}
({\cal X}_C)=(\tau_{r_2}-1-{\cal D}_{\ul{r_2}}^A)(\tau_{r_1}-{\cal
D}_{\ul{r_1}}^A +{\cal D}_{\ol{r_1}}^A+{\cal D}_{r_1}^C)({\cal
X}_C),\eqno(5.56)$$
$$(\vt+{\cal D}^C)\ptl_{z_{n+r_2,r_1}}({\cal X}_C)=(\tau_{r_2}-{\cal D}_{\ul{r_2}}^A
+{\cal D}_{\ol{r_2}}^A+{\cal D}^C_{r_2})(\tau_{r_1}-{\cal
D}_{\ul{r_1}}^A+{\cal D}_{\ol{r_1}}^A +{\cal D}_{r_1}^C)({\cal
X}_C)\eqno(5.57)$$ for $1\leq r_1<r_2\leq n$ and
$$(\vt+{\cal D}^C)\ptl_{z_{n+s,s}}({\cal X}_C)=(\tau_s-{\cal D}_{\ul{s}}^A+{\cal D}_{\ol{s}}^A
+{\cal D}^C_s)(\tau_s+1-{\cal D}_{\ul{s}}^A+{\cal
D}_{\ol{s}}^A+{\cal D}_{s}^C)({\cal X}_C) \eqno(5.58)$$ for
$s\in\ol{1,n}$.} \psp

Set
$$\G_D=\sum_{1\leq r_1<r_2\leq n}(\mbb{N}\es_{r_2,r_1}+\mbb{N}\es_{n+r_2,r_1})\subset \G_C.
\eqno(5.59)$$
Moreover, we let
$$\be_r^D=\sum_{i=1}^{r-1}\be_{n+r,i}+\sum_{s=r+1}^n\be_{n+s,r},\;\;\be^D=\sum_{1\leq r_1<r_2
\leq n}\be_{n+r_2,r_1}\qquad\for\;\;\be\in\G_D\eqno(5.60)$$ and
$$\be_r^B=\sum_{i=1}^r\be_{n+r,i}+\sum_{s=r+1}^n\be_{n+s,r},
\;\;\be^B=\be^c=\sum_{1\leq r_1\leq r_2\leq
n}\be_{n+r_2,r_1}\qquad\for\;\;\be\in\G_C. \eqno(5.61)$$ We define
the following hypergeometric functions:
$${\cal X}_D(\tau_1,...,\tau_n;\vt)\{z_{r_2,r_1},z_{n+r_2,r_1}\}=\sum_{\al\in\G_D}
\frac{\prod_{r=1}^n(\tau_r-\al_{\ul{r}})_{_{\al_{\ol{r}}+\al^D_r}}}{\al!(\vt)_{\al^D}}z^\al.
\eqno(5.62)$$
$${\cal X}_B(\tau_1,...,\tau_n;\vt)\{z_{r_2,r_1},z_{n+s_2,s_1}\}=\sum_{\al\in\G_C}
\frac{\prod_{r=1}^n(\tau_r-\al_{\ul{r}})_{_{\al_{\ol{r}}+\al^B_r}}}{\al!(\vt)_{\al^B}}z^\al.
\eqno(5.63)$$ Furthermore, we define ${\cal X}_B,\;{\cal
X}_B[i_1,i_2],\;{\cal X}_B[j_1,j_2;1],\;{\cal X}_B[k(2)]$ and
${\cal X}_D,\;{\cal X}_D[i_1,i_2],\;{\cal X}_D[j_1,j_2;1],$ ${\cal
X}_D[k(2)]$ as those of type C in (5.37) and (5.41)-(5.43). We let
$${\cal X}_B[k]={\cal X}_N(\tau_1,...,\tau_k+2,...,\tau_n;\vt+1)\{z_{r_2,r_1},z_{n+s_2,s_1}\}
\eqno(5.64)$$ obtained from ${\cal X}_B$  by changing $\tau_k$ to
$\tau_k+1$ and $\vt$ to $\vt+1$ for $k\in\ol{1,n}$. \psp

{\bf Theorem 5.7}. {\it The following equations hold for ${\cal
X}_C$:
$$\ptl_{z_{r_2,r_1}}({\cal X}_B)=\sum_{s=1}^{r_1}\tau_sP_{[s,r_1]}{\cal X}_B[s,r_2],
\eqno(5.65)$$
\begin{eqnarray*}\ptl_{z_{n+r_2,r_1}}({\cal X}_B)&=&\frac{1}{\vt}[\sum_{i=1}^{r_1}
\tau_i^2P_{[i,r_1]}P_{[i,r_2]}{\cal
X}_B[i(2)]+\sum_{s_1=1}^{r_1}\sum_{s_2=r_1+1}^{r_2}
\tau_{s_1}\tau_{s_2}P_{[s_1,r_1]}P_{[s_2,r_2]}{\cal
X}_B[s_1,s_2;1]\\ & &+\sum_{1\leq s_1<s_2 \leq
r_1}\tau_{s_1}\tau_{s_2}(P_{[s_1,r_1]}P_{[s_2,r_2]}+P_{[s_2,r_1]}P_{[s_1,r_2]})
{\cal X}_B[s_1,s_2;1]],\hspace{0.9cm}(5.66)\end{eqnarray*}
$$\ptl_{z_{r_2,r_1}}({\cal X}_D)=\sum_{s=1}^{r_1}\tau_sP_{[s,r_1]}{\cal X}_D[s,r_2],\eqno(5.67)$$
\begin{eqnarray*}\ptl_{z_{n+r_2,r_1}}({\cal X}_D)&=&\frac{1}{\vt}[\sum_{i=1}^{r_1}
\tau_i^2P_{[i,r_1]}P_{[i,r_2]}{\cal
X}_D[i(2)]+\sum_{s_1=1}^{r_1}\sum_{s_2=r_1+1}^{r_2}
\tau_{s_1}\tau_{s_2}P_{[s_1,r_1]}P_{[s_2,r_2]}{\cal
X}_D[s_1,s_2;1]\\ & &+ \sum_{1\leq s_1<s_2\leq
r_1}\tau_{s_1}\tau_{s_2}(P_{[s_1,r_1]}P_{[s_2,r_2]}+
P_{[s_2,r_1]}P_{[s_1,r_2]}){\cal
X}_D[s_1,s_2;1]],\hspace{0.9cm}(5.68)\end{eqnarray*} for $1\leq
r_1<r_2\leq n$ and
$$\ptl_{z_{n+s,s}}({\cal X}_B)=\sum_{r=1}^s\tau_rP_{[s,r]}X_B[s]\eqno(5.69)$$
for $s\in\ol{1,n}$.}\psp

 Let
$${\cal D}^B_r=\sum_{i=1}^rz_{n+r,i}\ptl_{z_{n+r,i}}+\sum_{s=r+1}^nz_{n+s,r}\ptl_{z_{n+s,r}},
\eqno(5.70)$$
$${\cal D}^D_r=\sum_{i=1}^{r-1}z_{n+r,i}\ptl_{z_{n+r,i}}+\sum_{s=r+1}^nz_{n+s,r}
\ptl_{z_{n+s,r}}\eqno(5.71)$$ for $r\in\ol{1,n}$ and
$${\cal D}^B={\cal D}^C,\;\;{\cal D}^D=\sum_{1\leq r_1< r_2\leq n}z_{n+r_2,r_1}
\ptl_{z_{n+r_2,r_1}}.\eqno(5.72)$$ As Theorem 5.4, we have: \psp

{\bf Theorem 5.7}. {\it The functions ${\cal X}_B$ and ${\cal X}_C$ satisfy:
$$(\tau_{r_2}-1-{\cal D}_{\ul{r_2}}^A+{\cal D}_{\ol{r_2}}^A+{\cal D}^B_{r_2})\ptl_{z_{r_2,r_1}}
({\cal X}_B)=(\tau_{r_2}-1-{\cal D}_{\ul{r_2}}^A)(\tau_{r_1}-{\cal
D}_{\ul{r_1}}^A+ {\cal D}_{\ol{r_1}}^A+{\cal D}_{r_1}^B)({\cal
X}_B),\eqno(5.73)$$
$$(\tau_{r_2}-1-{\cal D}_{\ul{r_2}}^A+{\cal D}_{\ol{r_2}}^A+{\cal D}^D_{r_2})\ptl_{z_{r_2,r_1}}
({\cal X}_D)=(\tau_{r_2}-1-{\cal D}_{\ul{r_2}}^A)(\tau_{r_1}-{\cal
D}_{\ul{r_1}}^A+ {\cal D}_{\ol{r_1}}^A+{\cal D}_{r_1}^D)({\cal
X}_D),\eqno(5.74)$$
$$(\vt+{\cal D}^B)\ptl_{z_{n+r_2,r_1}}({\cal X}_B)=(\tau_{r_2}-{\cal D}_{\ul{r_2}}^A
+{\cal D}_{\ol{r_2}}^A+{\cal D}^B_{r_2})(\tau_{r_1}-{\cal
D}_{\ul{r_1}}^A+{\cal D}_{\ol{r_1}}^A +{\cal D}_{r_1}^B)({\cal
X}_B),\eqno(5.75)$$
$$(\vt+{\cal D}^D)\ptl_{z_{n+r_2,r_1}}({\cal X}_D)=(\tau_{r_2}-{\cal D}_{\ul{r_2}}^A
+{\cal D}_{\ol{r_2}}^A+{\cal D}^D_{r_2})(\tau_{r_1}-{\cal
D}_{\ul{r_1}}^A+{\cal D}_{\ol{r_1}}^A +{\cal D}_{r_1}^D)({\cal
X}_D)\eqno(5.76)$$ for $1\leq r_1<r_2\leq n$ and
$$(\vt+{\cal D}^B)\ptl_{z_{n+s,s}}({\cal X}_B)=(\tau_s-{\cal D}_{\ul{s}}^A+{\cal D}_{\ol{s}}^A
+{\cal D}^B_s)({\cal X}_B)\eqno(5.77)$$ for $s\in\ol{1,n}$.}

\section{Weyl Functions}

In this section, we show that certain variations of the  Weyl functions of the
finite-dimensional simple Lie algebras of types A, B, C and  D  give rise to special solutions
 of the Olshanesky-Perelomov equation (1.13) with different constants $K_1$-$K_4$ and $\nu$.
 Part of these results are already known. Nevertheless, we gave a self-contained complete
 exposition for reader's convenience. Like Theorem 2.2, the results in this section and further
  study on the systems of partial differential equations in last section would yield more
  solutions of the Olshanesky-Perelomov equation (1.13).

Let $f_{i,j}(z)\mid i,j\in\ol{1,n}\}$ be a set of one-variable differentiable functions and
 let $d_i$ be a one-variable differential operator in $z_i$ for $i\in\ol{1,n}$. It is easy to
 verified the following lemma:
\psp

{\bf Lemma 6.1}. {\it We have the following equation on differentiation of determinants}:
\begin{eqnarray*}\qquad & &(\sum_{i=1}^nd_i)\left(\left|\begin{array}{cccc}f_{1,1}(z_1)&
f_{1,2}(z_2)&\cdots & f_{1,n}(z_n)\\ f_{2,1}(z_1)& f_{2,2}(z_2)&\cdots & f_{2,n}(z_n)\\
 \vdots&\vdots&\vdots&\vdots\\ f_{n,1}(z_1)& f_{n,2}(z_2)&\cdots &f_{n,n}(z_n)\end{array}
 \right|\right)\hspace{6cm}\end{eqnarray*}
\begin{eqnarray*}\qquad &=&\sum_{i=1}^n\left|\begin{array}{cccc}f_{1,1}(z_1)& f_{1,2}(z_2)&
\cdots & f_{1,n}(z_n)\\ \vdots&\vdots&\vdots&\vdots\\ f_{i-1,1}(z_1)& f_{i-1,2}(z_2)&\cdots &
f_{i-1,n}(z_n)\\ d_1(f_{i,1}(z_1))& d_2(f_{i,2}(z_2))&\cdots & d_n(f_{i,n}(z_n))
\\ f_{i+1,1}(z_1)& f_{i+1,2}(z_2)&\cdots & f_{i+1,n}(z_n)\\ \vdots&\vdots&\vdots&\vdots\\
 f_{n,1}(z_1)& f_{n,2}(z_2)&\cdots & f_{n,n}(z_n)\end{array}\right|.\hspace{4.3cm}(6.1)
 \end{eqnarray*}

Suppose that $R^+_X\subset \mbb{R}^n$ is the set of positive roots of the finite-dimensional
 simple Lie algebra of type X (e.g. cf. Section 12.1 in [Hj]). We denote
$$z^\al=\prod_{i=1}^nz_i^{\al_i}\qquad\for\;\;\al\in\mbb{R}^n.\eqno(6.2)$$
The Weyl function
$${\cal W}_X=\prod_{\al\in R^+_X}(z^{\al/2}-z^{-\al/2}).\eqno(6.3)$$

First  the Weyl function of $sl(n)$ is
$${\cal W}_{A_{n-1}}=\prod_{1\leq i<j\leq n}(z_i^{1/2}z_j^{-1/2}-z_i^{-1/2}z_j^{1/2})=
(z_1z_2\cdots z_n)^{(1-n)/2}\prod_{1\leq i<j\leq n}(z_i-z_j),\eqno(6.4)$$
Denote the Vandermonde determinant
$$W(z_1,z_2,...,z_n)=\left|\begin{array}{cccc}1&1&\cdots &1\\ z_1&z_2&\cdots &z_n\\ z_1^2&z_2^2
&\cdots &z_n^2 \\ \vdots&\vdots&\vdots&\vdots\\ z_1^{n-1}&z_2^{n-1}&\cdots &z_n^{n-1}
\end{array}\right|=\prod_{1\leq i<j\leq n}(z_i-z_j).\eqno(6.5)$$
Then $W(z_1,z_2,...,z_n)$ is the fundamental part of ${\cal W}_{A_{n-1}}$. By Lemma 6.1,
 we have
\begin{eqnarray*} (\sum_{i=1}^n(z_i\ptl_{z_i})^2)(W(z_1,z_2,...,z_n))&=&(\sum_{i=1}^{n-1}i^2)
W(z_1,z_2,...,z_n)\\ &=&\frac{(n-1)n(2n-1)}{6}W(z_1,z_2,...,z_n).\hspace{2.9cm}(6.6)
\end{eqnarray*}
On the other hand,
\begin{eqnarray*} & &(\sum_{i=1}^n(z_i\ptl_{z_i})^2)(W(z_1,z_2,...,z_n))\\ &=&(\sum_{r=1}^n
(z_r\ptl_{z_r})^2)(\prod_{1\leq i<j\leq n}(z_i-z_j))\\
&=&\sum_{r=1}^nz_r\left[-\sum_{s=1}^{r-1} \left(\prod_{1\leq
i<j\leq n;\:(i,j)\neq (s,r)}(z_i-z_j)\right)+\sum_{s=r+1}^n
\left(\prod_{1\leq i<j\leq n;\:(i,j)\neq
(r,s)}(z_i-z_j)\right)\right]\\ & &+2\sum_{r=1}^n
z_r^2(\sum_{1\leq s_1<s_2<r}\left[\prod_{1\leq i<j\leq
n;\:(i,j)\neq (s_1,r),(s_2,r)}(z_i-z_j) \right]\\ & &+
\sum_{r<s_1<s_2\leq n}\left[\prod_{1\leq i<j\leq n;\:(i,j)\neq
(r,s_1),(r,s_2)}(z_i-z_j)\right]
\\ & &-\sum_{1\leq s_1<r<s_2\leq n}\left[\prod_{1\leq i<j\leq n;\:(i,j)\neq (s_1,r),(r,s_2)}
(z_i-z_j)\right])\\ &=&\sum_{1\leq s<r\leq n}(z_s-z_r)\prod_{1\leq i<j\leq n;\:(i,j)\neq (s,r)}
(z_i-z_j)\\ & &+2W(z_1,z_2,...,z_n)\sum_{r=1}^n\sum_{1\leq s_1<s_2\leq n;\:s_1,s_2\neq r}
\frac{z_r^2}{(z_{s_1}-z_r)(z_{s_2}-z_r)}\\ &=&\left(\frac{n(n-1)}{2}+2
\sum_{r=1}^n\sum_{1\leq s_1<s_2\leq n;\:s_1,s_2\neq r}\frac{z_r^2}{(z_{s_1}-z_r)(z_{s_2}-z_r)}
\right)W(z_1,z_2,...,z_n).\hspace{1cm}(6.7)\end{eqnarray*}
Thus (6.6) and (6.7) yield
\begin{eqnarray*}\qquad\qquad& &\sum_{r=1}^n\sum_{1\leq s_1<s_2\leq n;\:s_1,s_2\neq r}
\frac{z_r^2}{(z_{s_1}-z_r)(z_{s_2}-z_r)}\\
&=&\frac{1}{2}\left[\frac{(n-1)n(2n-1)}{6}
-\frac{n(n-1)}{2}\right]\\
&=&\frac{(n-1)n(n-2)}{6}=\left(\!\!\begin{array}{c}n\\
3\end{array} \!\!\right).\hspace{7.9cm}(6.8)\end{eqnarray*}

Let
$$\phi^A_{\mu_1,\mu_2}=(z_1z_2\cdots z_n)^{\mu_1}W^{\mu_2}(z_1,z_2,...,z_n)\qquad\for\;\;
\mu_1,\mu_2\in\mbb{C}.\eqno(6.9)$$
Then
$$z_r\ptl_{z_r}(\phi^A_{\mu_1,\mu_2})=\left(\mu_1-\mu_2\sum_{s=1}^{r-1}\frac{z_r}{z_s-z_r}
+\mu_2\sum_{s=r+1}^n\frac{z_r}{z_r-z_s}\right)\phi^A_{\mu_1,\mu_2}\eqno(6.10)$$
for $r\in\ol{1,n}$. Thus
\begin{eqnarray*}& &\sum_{r=1}^n(z_r\ptl_{z_r})^2(\phi^A_{\mu_1,\mu_2})\\ &=&\sum_{r=1}^n
[\mu_1^2-2\mu_1\mu_2\sum_{s=1}^{r-1}\frac{z_r}{z_s-z_r}+2\mu_1\mu_2\sum_{s=r+1}^n\frac{z_r}
{z_r-z_s}-\mu_2\sum_{r\neq s\in\ol{1,n}}\frac{z_sz_r}{(z_s-z_r)^2}\\
&&+\mu_2^2\sum_{r\neq
s\in\ol{1,n}}\frac{z_r^2}{(z_s-z_r)^2}+2\mu_2^2 \sum_{1\leq
s_1<s_2\leq n;\:s_1,s_2\neq
r}\frac{z_r^2}{(z_{s_1}-z_r)(z_{s_2}-z_r)}] \phi^A_{\mu_1,\mu_2}\\
&=&[n\mu_1^2+n(n-1)\mu_1\mu_2+\mu_2^2
\left(\!\!\begin{array}{c}n\\
3\end{array}\!\!\right)-2\mu_2\sum_{1\leq r<s\leq n}
\frac{z_sz_r}{(z_s-z_r)^2}\\ & &+2\mu_2^2\sum_{1\leq r<s\leq
n}\frac{z_r^2+z_s^2}
{(z_s-z_r)^2}]\phi^A_{\mu_1,\mu_2}\\&=&[n\mu_1^2+n(n-1)\mu_1\mu_2+2\mu_2^2
\left(\!\!\begin{array}{c}n\\
3\end{array}\!\!\right)-2\mu_2\sum_{1\leq r<s\leq n}
\frac{z_sz_r}{(z_s-z_r)^2}\\ & &+\mu_2^2\sum_{1\leq r<s\leq
n}\frac{z_r^2+z_s^2-2z_rz_s+2z_rz_s}{(z_s-z_r)^2}]\phi^A_{\mu_1,\mu_2}\\
&=&[n\mu_1^2+n(n-1)(\mu_1+\mu_2/2)\mu_2+2\mu_2^2\left(\!\!\begin{array}{c}n\\
3\end{array}\!\!\right)\\ & &+2\mu_2(\mu_2-1)\sum_{1\leq r<s\leq
n}\frac{z_sz_r}{(z_s-z_r)^2}]\phi^A_{\mu_1,\mu_2}\hspace{7.2cm}(6.11)\end{eqnarray*}
by (6.8) and (6.10). Therefore, we have: \psp

{\bf Theorem 6.2}. {\it The function $\phi^A_{\mu_1,\mu_2}$ satisfies:
\begin{eqnarray*}\qquad& &\sum_{r=1}^n(z_r\ptl_{z_r})^2(\phi^A_{\mu_1,\mu_2})+2\mu_2(1-\mu_2)
\left(\sum_{1\leq i<j\leq n}\frac{z_iz_j}{(z_i-z_j)^2}\right)\phi^A_{\mu_1,\mu_2}\\
&=&\left[n\mu_1^2+n(n-1)(\mu_1+\mu_2/2)\mu_2+2\left(\!\!\begin{array}{c}n\\ 3\end{array}
\!\!\right)\mu_2^2\right]\phi^A_{\mu_1,\mu_2},\hspace{4cm}(6.12)\end{eqnarray*}
which is a Calogero-Sutherland equation.}
\psp

We remark that above result was  known when $\mu_1=\mu_2$ or $\mu_1=0$. Moreover,
${\cal W}_{A_n}=\phi^A_{(1-n)/2,1}$.

 According to simplicity, we secondly consider the Weyl function of  $so(2n)$. By the root
 structure,
\begin{eqnarray*}\qquad{\cal W}_{D_n}&=&\prod_{1\leq i<j\leq n}(z_i^{1/2}z_j^{-1/2}-z_i^{-1/2}
z_j^{1/2})(z_i^{1/2}z_j^{1/2}-z_i^{-1/2}z_j^{-1/2})\\ &=&
\prod_{1\leq i<j\leq n}(z_i+z_i^{-1}-(z_j+z_j^{-1}))\\
&=&W(z_1+z_1^{-1},z_2+z_2^{-1},...,z_n
+z_n^{-1})\hspace{5.8cm}(6.13)\end{eqnarray*}(cf. (6.5)). Note
that
\begin{eqnarray*}(z_i\ptl_{z_i})^2[(z_i+z_i^{-1})^k]&=&k(z_i+z_i^{-1})^k+k(k-1)(z_i-z_i^{-1})^2
(z_i+z_i^{-1})^{k-2}\\ &=& k(z_i+z_i^{-1})^k+k(k-1)(z_i^2-2+z_i^{-2})(z_i+z_i^{-1})^{k-2}\\
 &=& k^2(z_i+z_i^{-1})^k-4k(k-1)(z_i+z_i^{-1})^{k-2}.\hspace{3.4cm}(6.14)\end{eqnarray*}
By Lemma 6.1 and a similar calculation as (6.5),
$$(\sum_{i=1}^n(z_i\ptl_{z_i})^2)({\cal W}_{D_n})=\frac{(n-1)n(2n-1)}{6}{\cal W}_{D_n}.
\eqno(6.15)$$
On the other hand, by (6.7) and (6.14), we have
\begin{eqnarray*} & &\sum_{r=1}^n\sum_{1\leq s_1<s_2\leq n;\:s_1,s_2\neq r}
\frac{(z_r-z_r^{-1})^2}{(z_{s_1}+z_{s_1}^{-1}-z_r-z_r^{-1})(z_{s_2}+z_{s_2}^{-1}-z_r-z_r^{-1})}
\\ &=&\frac{(n-1)n(n-2)}{6}=\left(\!\!\begin{array}{c}n\\ 3\end{array}\!\!\right).
\hspace{9.4cm}(6.16)\end{eqnarray*} Moreover,
\begin{eqnarray*} & &\frac{(z_r-z_r^{-1})^2+(z_s-z_s^{-1})^2}{(z_s+z_s^{-1}-z_r-z_r^{-1})^2}\\
 &=& \frac{z_r^2+z_r^{-2}+z_s^2+z_s^{-2}-4}{(z_rz_s)^{-2}(z_s-z_r)^2(z_rz_s-1)^2}\\ &=&
 \frac{z_r^4z_s^2+z_s^2+z_r^2z_s^4+z_s^2-4z_r^2z_s^2}{(z_s-z_r)^2(z_rz_s-1)^2}\\ &=&
\frac{1}{(z_s-z_r)^2(z_rz_s-1)^2}
 [z_r^4z_s^2+z_s^2+z_r^2z_s^4+z_s^2-2z_rz_s(z_r^2+z_s^2+(z_rz_s)^2-2z_rz_s+1)\\ & &
 +2z_rz_s(z_r^2+z_s^2+(z_rz_s)^2-4z_rz_s+1)]\\ &=&\frac{(z_s-z_r)^2(z_rz_s-1)^2
 +2z_rz_s[(z_r-z_s)^2+(z_rz_s-1)^2]}{(z_s-z_r)^2(z_rz_s-1)^2}\\ &=&1 +2\frac{z_rz_s}
 {(z_r-z_s)^2}+2\frac{z_rz_s}{(z_rz_s-1)^2}\hspace{8cm}(6.17)\end{eqnarray*}
and
\begin{eqnarray*}\qquad& &z_r\ptl_{z_r}\left(\frac{z_r-z_r^{-1}}{z_r+z_r^{-1}-z_s-z_s^{-1}}
\right)\\ &=&\frac{(z_r+z_r^{-1})(z_r+z_r^{-1}-z_s-z_s^{-1})-(z_r-z_r^{-1})^2}{(z_r+z_r^{-1}
-z_s-z_s^{-1})^2}\\ &=& \frac{(z_r+z_r^{-1})^2-(z_r+z_r^{-1})(z_s+z_s^{-1})-(z_r-z_r^{-1})^2}
{(z_r+z_r^{-1}-z_s-z_s^{-1})^2}\\ &=&
 \frac{4-z_rz_s-z_rz_s^{-1}-z_r^{-1}z_s-z_r^{-1}z_s^{-1}}{(z_rz_s)^{-2}(z_s-z_r)^2(z_rz_s-1)^2}
 \\ &=& -z_rz_s\frac{(z_rz_s)^2+z_r^2+z_s^2+1-4z_rz_s}{(z_s-z_r)^2(z_rz_s-1)^2}\\ &=&
 -\frac{z_rz_s}{(z_s-z_r)^2}-\frac{z_rz_s}{(z_rz_s-1)^2}\hspace{8.4cm}(6.18)\end{eqnarray*}
for $r,s\in\ol{1,n}$ such that $r\neq s$.

Set
$$\phi_\mu^D=({\cal W}_{D_n})^\mu=\prod_{1\leq i<j\leq n}(z_i+z_i^{-1}-z_j-z_j^{-1})^{\mu}
\qquad\for\;\;\mu\in\mbb{C}.\eqno(6.19)$$
Then
\begin{eqnarray*} & &(\sum_{r=1}^n(z_r\ptl_{z_r})^2)(\phi_\mu^D)\\ &=&\sum_{r=1}^n
[\sum_{r\neq s\in\ol{1,n}}\left[\mu
z_r\ptl_{z_r}\left(\frac{z_r-z_r^{-1}}
{z_r+z_r^{-1}-z_s-z_s^{-1}}\right)+\mu^2\frac{(z_r-z_r^{-1})^2}{(z_s+z_s^{-1}-z_r-z_r^{-1})^2}
\right]\\ & &+2\mu^2\sum_{1\leq s_1<s_2\leq n;\:s_1,s_2\neq
r}\frac{(z_r-z_r^{-1})^2}{(z_{s_1}
+z_{s_1}^{-1}-z_r-z_r^{-1})(z_{s_2}+z_{s_2}^{-1}-z_r-z_r^{-1})}]\phi_\mu^D\\
&=& [\sum_{1\leq r<s\leq
n}\left[-2\mu\left(\frac{z_rz_s}{(z_s-z_r)^2}+\frac{z_rz_s}{(z_rz_s-1)^2}
\right)+\mu^2\frac{(z_r-z_r^{-1})^2+(z_s-z_s^{-1})^2}{(z_s+z_s^{-1}-z_r-z_r^{-1})^2}\right]
\phi_\mu^D\\&&+2\mu^2\sum_{r=1}^n\;\sum_{1\leq s_1<s_2\leq
n;\:s_1,s_2\neq r}
\frac{(z_r-z_r^{-1})^2}{(z_{s_1}+z_{s_1}^{-1}-z_r-z_r^{-1})(z_{s_2}+z_{s_2}^{-1}-z_r-z_r^{-1})}]
\\ &=& \left[2\mu(\mu-1)\sum_{1\leq r<s\leq n}\left(\frac{z_rz_s}{(z_s-z_r)^2}+\frac{z_rz_s}
{(z_rz_s-1)^2}\right)+\mu^2\left(\frac{n(n-1)}{2}+2\left(\!\!\begin{array}{c}n\\
3\end{array} \!\!\right)\right)\right]\phi_\mu^D\\
&=&\left[2\mu(\mu-1)\sum_{1\leq r<s\leq
n}\left(\frac{z_rz_s}{(z_s-z_r)^2}+\frac{z_rz_s}{(z_rz_s-1)^2}\right)
+\frac{n(n-1)(2n-1)}{6}\mu^2\right]\phi_\mu^D.\hspace{0.5cm}(6.20)\end{eqnarray*}
Thus we obtain: \psp

{\bf Theorem 6.3}. {\it The function $\phi_\mu^D$ the equation:
\begin{eqnarray*}& &\sum_{r=1}^n(z_r\ptl_{z_r})^2(\phi^D_\mu)+2\mu(1-\mu)
\left(\sum_{1\leq i<j\leq n}\frac{z_iz_j}{(z_i-z_j)^2}+\sum_{1\leq i<j\leq n}\frac{z_iz_j}
{(z_iz_j-1)^2}\right)\phi^D_\mu\\ &=&\frac{n(n-1)(2n-1)\mu^2}{6}\phi_\mu^D,
\hspace{9.7cm}(6.21)\end{eqnarray*}
which is an Olshanesky-Perelomov equation}.
\psp

We remark that the above result with $n=4$ was given in [NFP].

The Weyl function of $sp(2n)$ is
$${\cal W}_{C_n}=\left[\prod_{r=1}^n(z_r-z_r^{-1})\right]\prod_{1\leq i<j\leq n}
(z_i+z_i^{-1}-(z_j+z_j^{-1})),\eqno(6.22)$$
In terms of determinant,
$${\cal W}_{C_n}=\left|\begin{array}{ccc}z_1-z_1^{-1}&\cdots& z_n-z_n^{-1}\\ (z_1-z_1^{-1})
(z_1+z_1^{-1})&\cdots& (z_n-z_n^{-1})(z_n+z_n^{-1})\\ (z_1-z_1^{-1})(z_1+z_1^{-1})^2&
\cdots& (z_n-z_n^{-1})(z_n+z_n^{-1})^2\\ \vdots&\vdots&\vdots\\ (z_1-z_1^{-1})
(z_1+z_1^{-1})^{n-1}&\cdots& (z_n-z_n^{-1})(z_n+z_n^{-1})^{n-1}\end{array}\right|.
\eqno(6.23)$$
Moreover,
\begin{eqnarray*} (z_r\ptl{z_r})^2[(z_r-z_r^{-1})(z_r+z_r^{-1})^k]&=&(k+1)^2(z_r-z_r^{-1})
(z_r+z_r^{-1})^k\\ &
&+4k(k-1)(z_r-z_r^{-1})(z_r+z_r^{-1})^{k-2}\hspace{2.3cm}(6.24)
\end{eqnarray*}
by (6.14). Thus
$$(\sum_{i=1}^n(z_i\ptl_{z_i})^2)({\cal W}_{C_n})=\frac{n(n+1)(2n+1)}{6}{\cal W}_{C_n}
\eqno(6.25)$$
by Lemma 6.1 and (6.24).

Set
$$\psi^C_\mu=\prod_{r=1}^n(z_r-z_r^{-1})^{\mu}\qquad\for\;\;\mu\in\mbb{C}.\eqno(6.26)$$
Then
\begin{eqnarray*}\sum_{r=1}^n(z_r\ptl{z_r})^2(\psi^C_\mu)&=&\mu\psi^C_\mu\sum_{r=1}^n
\left[z_r\ptl_{z_r}\left(\frac{z_r+z_r^{-1}}{z_r-z_r^{-1}}\right)+\mu\left(\frac{z_r+z_r^{-1}}
{z_r-z_r^{-1}}\right)^2\right]\\ &=&
\mu\psi^C_\mu\sum_{r=1}^n\left[\frac{(z_r-z_r^{-1})^2-(z_r+z_r^{-1})^2}{(z_r-z_r^{-1})^2}
+\mu\frac{(z_r+z_r^{-1})^2}{(z_r-z_r^{-1})^2}\right]\\
&=&\mu\sum_{r=1}^n\left[\mu+4(\mu-1)
\frac{1}{(z_r-z_r^{-1})^2}\right]\psi^C_\mu\\
&=&\left[n\mu^2+4\mu(\mu-1)\sum_{r=1}^n
\frac{z_r^2}{(z_r^2-1)^2}\right]\psi^C_\mu.\hspace{4.5cm}(6.27)\end{eqnarray*}
Moreover,
\begin{eqnarray*}\qquad& &\sum_{r=1}^n(z_r\ptl{z_r})(\psi^C_{\mu_1})(z_r\ptl{z_r})(\phi_{\mu_2}
^D)\\
&=&\mu_1\mu_2\sum_{r=1}^n\frac{z_r+z_r^{-1}}{z_r-z_r^{-1}}\psi^C_{\mu_1}
\sum_{r\neq
s\in\ol{1,n}}\frac{z_r-z_r^{-1}}{z_r+z_r^{-1}-z_s-z_s^{-1}}\phi_{\mu_2}^D
\hspace{6cm}\end{eqnarray*}
\begin{eqnarray*}
&=&\mu_1\mu_2\psi^C_{\mu_1}\phi_{\mu_2}^D\sum_{r=1}^n\sum_{r\neq
s\in\ol{1,n}} \frac{z_r+z_r^{-1}}{z_r+z_r^{-1}-z_s-z_s^{-1}}\\
&=&\frac{n(n-1)}{2}\mu_1\mu_2\psi^C_{\mu_1}
\phi_{\mu_2}^D.\hspace{10cm}(6.28)\end{eqnarray*} Set
$$\phi_{\mu_1,\mu_2}^C=\psi^C_{\mu_1}\phi_{\mu_2}^D=\left[\prod_{r=1}^n(z_r-z_r^{-1})^{\mu_1}
\right]\prod_{1\leq i<j\leq n}(z_i+z_i^{-1}-z_j-z_j^{-1})^{\mu_2}\eqno(6.29)$$
for $\mu_1,\mu_2\in\mbb{C}$.
By Theorem 6.3 and (6.27)-(6.29), we have:
\psp

{\bf Theorem 6.4}. {\it The function $\phi^C_{\mu_1,\mu_2}$ satisfies:
\begin{eqnarray*}& &\sum_{r=1}^n(z_r\ptl_{z_r})^2(\phi^C_{\mu_1,\mu_2})+[4\mu_1(1-\mu_1)
\sum_{r=1}^n\frac{z_r^2}{(z_r^2-1)^2}\\ &
&+2\mu_2(1-\mu_2)\left(\sum_{1\leq i<j\leq n}
\frac{z_iz_j}{(z_i-z_j)^2}+\sum_{1\leq i<j\leq
n}\frac{z_iz_j}{(z_iz_j-1)^2}\right)] \phi^C_{\mu_1,\mu_2}\\
&=&\left(n\mu^2_1+n(n-1)\mu_1\mu_2+\frac{n(n-1)(2n-1)\mu_2^2}{6}\right)
\phi^C_{\mu_1,\mu_2},\hspace{4.4cm}(6.30)\end{eqnarray*} which is
an Olshanesky-Perelomov equation}. \psp

Observe that ${\cal W}_{C_n}=\phi^C_{1,1}$.
The Weyl function of $so(2n+1)$ is
$${\cal W}_{B_n}=\left[\prod_{r=1}^n(z_r^{1/2}-z_r^{-1/2})\right]\prod_{1\leq i<j\leq n}
(z_i+z_i^{-1}-(z_j+z_j^{-1})),\eqno(6.31)$$
In terms of determinant,
\begin{eqnarray*}& &{\cal W}_{B_n}=\\
& &\left|\begin{array}{ccc}z_1^{1/2}-z_1^{-1/2}&\cdots& z_n^{1/2}-z_n^{-1/2}\\ (z_1^{1/2}
-z_1^{-1/2})(z_1^{1/2}+z_1^{-1/2})^2&\cdots& (z_n^{1/2}-z_n^{-1/2})(z_n^{1/2}+z_n^{-1/2})^2\\
(z_1^{1/2}-z_1^{-1/2})(z_1^{1/2}+z_1^{-1/2})^4&\cdots&
(z_n^{1/2}-z_n^{-1/2})(z_n^{1/2} +z_n^{-1/2})^4\\
\vdots&\vdots&\vdots\\
(z_1^{1/2}-z_1^{-1/2})(z_1^{1/2}+z_1^{-1/2})^{2(n-1)} &\cdots&
(z_n^{1/2}-z_n^{-1/2})(z_n^{1/2}+z_n^{-1/2})^{2(n-1)}\end{array}\right|.
\hspace{0.6cm}(6.32)\end{eqnarray*} Moreover,
\begin{eqnarray*}\qquad & &(z_r\ptl_{z_r})^2[(z_r^{1/2}-z_r^{-1/2})(z_r^{1/2}+z_r^{-1/2})^{2k}]
\\ &=&\left(k^2+k+\frac{1}{4}\right)(z_r^{1/2}-z_r^{-1/2})(z_r^{1/2}+z_r^{-1/2})^{2k}\\
&
&+2k(2k-1)(z_r^{1/2}-z_r^{-1/2})(z_r^{1/2}+z_r^{-1/2})^{2(k-1)}.\hspace{4.9cm}(6.33)
\end{eqnarray*}
 Thus
$$(\sum_{i=1}^n(z_i\ptl_{z_i})^2)({\cal W}_{B_n})=\frac{n(4n^2-1)}{12}{\cal W}_{B_n}\eqno(6.34)$$
by Lemma 6.1 and (6.33).

Set
$$\psi^B_\mu=\prod_{r=1}^n(z_r^{1/2}-z_r^{-1/2})^{\mu}\qquad\for\;\;\mu\in\mbb{C}.\eqno(6.35)$$
Then
\begin{eqnarray*}\sum_{r=1}^n(z_r\ptl{z_r})^2(\psi^B_\mu)&=&\frac{\mu}{2}\psi^B_\mu
\sum_{r=1}^n\left[z_r\ptl{z_r}\left(\frac{z_r^{1/2}+z_r^{-1/2}}{z_r^{1/2}-z_r^{-1/2}}\right)
+\frac{\mu}{2}\left(\frac{z_r^{1/2}+z_r^{-1/2}}{z_r^{1/2}-z_r^{-1/2}}\right)^2\right]\\ &=&
\frac{\mu}{4}\psi^B_\mu\sum_{r=1}^n\left[\frac{(z_r^{1/2}-z_r^{-1/2})^2-(z_r^{1/2}+
z_r^{-1/2})^2}{(z_r^{1/2}-z_r^{-1/2})^2}+\mu\frac{(z_r^{1/2}+z_r^{-1/2})^2}{(z_r^{1/2}
-z_r^{-1/2})^2}\right]\\  &=&\left[\frac{n\mu^2}{4}+\mu(\mu-1)\sum_{r=1}^n\frac{z_r}
{(z_r-1)^2}\right]\psi^B_\mu.\hspace{4.5cm}(6.36)\end{eqnarray*}
Moreover,
\begin{eqnarray*}\qquad& &\sum_{r=1}^n(z_r\ptl{z_r})(\psi^B_{\mu_1})(z_r\ptl{z_r})
(\phi_{\mu_2}^D)\\ &=&\frac{\mu_1\mu_2}{2}\sum_{r=1}^n\frac{z_r^{1/2}+z_r^{-1/2}}
{z_r^{1/2}-z_r^{-1/2}}\psi^B_{\mu_1}\sum_{r\neq s\in\ol{1,n}}\frac{z_r-z_r^{-1}}
{z_r+z_r^{-1}-z_s-z_s^{-1}}\phi_{\mu_2}^D\\ &=&\frac{\mu_1\mu_2}{2}\psi^B_{\mu_1}
\phi_{\mu_2}^D\sum_{r=1}^n\sum_{r\neq s\in\ol{1,n}}\frac{(z_r^{1/2}+z_r^{-1/2})^2}
{z_r+z_r^{-1}-z_s-z_s^{-1}}\\ &=&\frac{n(n-1)}{4}\mu_1\mu_2\psi^B_{\mu_1}\phi_{\mu_2}^D.
\hspace{9.2cm}(6.37)\end{eqnarray*}
Set
$$\phi_{\mu_1,\mu_2}^B=\psi^B_{\mu_1}\phi_{\mu_2}^D=\left[\prod_{r=1}^n(z_r^{1/2}
-z_r^{-1/2})^{\mu_1}\right]\prod_{1\leq i<j\leq n}(z_i+z_i^{-1}-z_j-z_j^{-1})^{\mu_2}
\eqno(6.38)$$
for $\mu_1,\mu_2\in\mbb{C}$.
By Theorem 6.3 and (6.36)-(6.38), we obtain:
\psp

{\bf Theorem 6.5}. {\it The function $\phi^B_{\mu_1,\mu_2}$ satisfies:
\begin{eqnarray*}& &\sum_{r=1}^n(z_r\ptl_{z_r})^2(\phi^B_{\mu_1,\mu_2})+
[\mu_1(1-\mu_1)\sum_{r=1}^n\frac{z_r}{(z_r-1)^2}\\ & &+2\mu_2(1-\mu_2)
\left(\sum_{1\leq i<j\leq n}\frac{z_iz_j}{(z_i-z_j)^2}+\sum_{1\leq i<j\leq n}\frac{z_iz_j}
{(z_iz_j-1)^2}\right)]\phi^C_{\mu_1,\mu_2}\\ &=&\left(\frac{n\mu^2_1}{4}+\frac{n(n-1)}{2}
\mu_1\mu_2+\frac{n(n-1)(2n-1)\mu_2^2}{6}\right)\phi^B_{\mu_1,\mu_2},\hspace{4.2cm}(6.39)
\end{eqnarray*}
which is an Olshanesky-Perelomov equation}.

\section{Olshanesky-Perelomov Equation of Type C}

In this section, we want to show that the trace functions of
certain intertwining operators among $sp(2n)$-modules give rise to
solutions of the equation (1.13) with $K_3=0$, following Etingof's
idea in [Ep].

Recall the sympletic Lie algebra $sp(2n)$ in (4.1) and the notions
$\{C_{r,s}\}$ in (4.8)-(4.10). The Casimier element
\begin{eqnarray*}\qquad \Omega&=&\sum_{i,j=1}^nC_{i,j}C_{j,i}+\sum_{1\leq s<r\leq n}
(C_{n+r,s}C_{s,n+r}+C_{s,n+r}C_{n+r,s})\\ &
&+2\sum_{p=1}^n(C_{n+p,p}C_{p,n+p}+C_{p,n+p}
C_{n+p,p}),\hspace{6.2cm}(7.1)\end{eqnarray*} which is a quadratic
central element in $U(sp(2n))$. Denote by  $T^\ast$  the set of
all
 linear functions (weights) on the toral Cartan subalgebra $T$ in (4.3).  Let
$$M=\bigoplus_{\lmd\in T^\ast}M^\lmd\eqno(7.2)$$
be any weight $sp(2n)$-module with the representation $\pi_{_M}$
such that
$$\pi_{_M}(\Omega)=\mu \mbox{Id}_M\qquad\mbox{for some}\;\;\mu\in \mbb{C},\eqno(7.3)$$
where
$$M^\lmd=\{w\in M\mid h(w)=\lmd (h)w\;\for\;h\in T\}.\eqno(7.4)$$
In particular, we can take $M$ to be any highest weight module; say $M_\lmd$ Section 4. The
 module $M$ is not necessarily irreducible. We use the notations in (4.29) and (4.30). Denote
 the space of formal Laurent series in $\{x_1,x_2,...,x_n\}$ with the coefficients in $M$ by
$$\td{M}=\{\sum_{\vec i\in\mbb{Z}^n}w_{\vec i}x^{\vec i}x^\ast\mid w_n\in M\}.\eqno(7.5)$$
We extend the $\pi_{_M}$ to $\td{M}$ by
$$\pi_{_M}(u)(\sum_{\vec i\in\mbb{Z}^n}w_{\vec i}x^{\vec i}x^\ast)=\sum_{\vec i\in\mbb{Z}^n}
\pi_{_M}(u)(w_{\vec i})x^{\vec i}x^\ast\qquad\for\;\;u\in sp(2n).\eqno(7.6)$$
Moreover, we define another representation $\pi'$ of $sp(2n)$ on $\td{M}$ by
$$\pi'(C_{p,q})(\sum_{\vec i\in\mbb{Z}^n}w_{\vec i}x^{\vec i}x^\ast)=\sum_{\vec i\in\mbb{Z}^n}
\left(C_{p,q}(w_{\vec i}x_p\ptl_{x_q}(x^{\vec i}x^\ast)+\frac{\dlt_{p,q}}{2}w_{\vec i}x^{\vec i}
x^\ast\right),\eqno(7.7)$$
$$\pi'(C_{p,n+q})(\sum_{\vec i\in\mbb{Z}^n}w_{\vec i}x^{\vec i}x^\ast)=-\frac{1}{1+\dlt_{p,q}}
\sum_{\vec i\in\mbb{Z}^n}w_{\vec i}x_px_qx^{\vec i}x^\ast,\eqno(7.8)$$
$$\pi'(C_{n+p,q})(\sum_{\vec i\in\mbb{Z}^n}w_{\vec i}x^{\vec i}x^\ast)=
\frac{1}{1+\dlt_{p,q}}\sum_{\vec i\in\mbb{Z}^n}w_{\vec i}\ptl_{x_p}\ptl_{x_q}
(x^{\vec i}x^\ast)\eqno(7.9)$$
for $p,q\in\ol{1 n}$. Set
$$\pi=\pi_{_M}+\pi'.\eqno(7.10)$$

Suppose that $\Phi:M\rta \td{M}$ is a linear map such that
$$\Phi(\pi_{_M}(u)(w))=\pi(u)\Phi(w)\eqno(7.11)$$
for $u\in sp(2n)$ and $w\in M$. View $\Phi$ as a function in $\{x_1,x_2,...,x_n\}$ taking value
 in the spaces of  linear transformations on $M_\lmd$. Define the  trace function:
$${\cal E}=(x^\ast)^{-1}\mbox{tr}_M\Phi z_1^{\pi_{_M}(C_{1,1})} z_2^{\pi_{_M}(C_{2,2})}
\cdots  z_n^{\pi_{_M}(C_{n,n})}.\eqno(7.12)$$ We want to prove
that the function
$$\Psi={\cal W}_{C_n}{\cal E}\eqno(7.13)$$
satisfies the Olshanesky-Perelomov equation (1.13) with $K_3=0$. Note that (6.22) implies
\begin{eqnarray*}z_i\ptl_{z_i}({\cal W}_{C_n})&=&\sum_{i=1}^n\left[\frac{z_i+z_i^{-1}}{z_i
-z_i^{-1}}+\sum_{i\neq r\in\ol{1,n}}\frac{z_i-z_i^{-1}}{z_i+z_i^{-1}-z_r-z_r^{-1}}\right]
{\cal W}_{C_n}\\ &=&\sum_{i=1}^n\left[\frac{z_i^2+1}{z_i^2-1}+\sum_{i\neq r\in\ol{1,n}}
\left(\frac{z_i}{z_i-z_r}+\frac{1}{z_iz_r-1}\right)\right]{\cal W}_{C_n}.\hspace{3cm}
(7.14)\end{eqnarray*}

Set
$${\cal E}^{i_1,i_2}=(x^\ast)^{-1}\mbox{tr}_M\Phi \pi_{_M}( C_{i_1,i_2}C_{i_2,i_1})
z_1^{\pi_{_M}(C_{1,1})}z_2^{\pi_{_M}(C_{2,2})}\cdots
z_n^{\pi_{_M}(C_{n,n})}\eqno(7.15)$$ for $i_1,i_2\in\ol{1,n}$.
Then
$${\cal E}^{i,i}=(z_i\ptl_{z_i})^2{\cal E}\eqno(7.16)$$
and for $i_1\neq i_2$,
\begin{eqnarray*} x^{\ast}{\cal E}^{i_1,i_2}&=&\mbox{tr}_M\pi(C_{i_1,i_2})\Phi \pi_{_M}
(C_{i_2,i_1})z_1^{\pi_{_M}(C_{1,1})}z_2^{\pi_{_M}(C_{2,2})}\cdots
z_n^{\pi_{_M}(C_{n,n})}
\\ &=& \mbox{tr}_M\pi'(C_{i_1,i_2})\Phi \pi_{_M}(C_{i_2,i_1})z_1^{\pi_{_M}(C_{1,1})}
z_2^{\pi_{_M}(C_{2,2})}\cdots  z_n^{\pi_{_M}(C_{n,n})}\\ &
&+\mbox{tr}_M\pi_{_M}(C_{i_1,i_2}) \Phi
\pi_{_M}(C_{i_2,i_1})z_1^{\pi_{_M}(C_{1,1})}z_2^{\pi_{_M}(C_{2,2})}\cdots
z_n^{\pi_{_M} (C_{n,n})}\\ &=&\mbox{tr}_M\pi'(C_{i_1,i_2})\Phi
\pi_{_M}(C_{i_2,i_1})z_1^{\pi_{_M}(C_{1,1})}
z_2^{\pi_{_M}(C_{2,2})}\cdots  z_n^{\pi_{_M}(C_{n,n})}\\ &
&+\mbox{tr}_M \Phi
\pi_{_M}(C_{i_2,i_1})z_1^{\pi_{_M}(C_{1,1})}z_2^{\pi_{_M}(C_{2,2})}\cdots
z_n^{\pi_{_M}(C_{n,n})}\pi_{_M}(C_{i_1,i_2})\\
&=&\mbox{tr}_M\pi'(C_{i_1,i_2}) \Phi
\pi_{_M}(C_{i_2,i_1})z_1^{\pi_{_M}(C_{1,1})}z_2^{\pi_{_M}(C_{2,2})}\cdots
z_n^{\pi_{_M}(C_{n,n})}\\ &
&+\frac{z_{i_1}}{z_{i_2}}\mbox{tr}_M\Phi \pi_{_M}(C_{i_2,i_1}
C_{i_1,i_2})z_1^{\pi_{_M}(C_{1,1})}z_2^{\pi_{_M}(C_{2,2})}\cdots  z_n^{\pi_{_M}(C_{n,n})}\\
 &=& \mbox{tr}_M\pi'(C_{i_1,i_2})\Phi \pi_{_M}(C_{i_2,i_1})z_1^{\pi_{_M}(C_{1,1})}
 z_2^{\pi_{_M}(C_{2,2})}\cdots  z_n^{\pi_{_M}(C_{n,n})}\\ & &+\frac{z_{i_1}}{z_{i_2}}
 \mbox{tr}_M\Phi \pi_{_M}(C_{i_1,i_2}C_{i_2,i_1}+[C_{i_2,i_1},C_{i_1,i_2}])z_1^{\pi_{_M}
 (C_{1,1})}z_2^{\pi_{_M}(C_{2,2})}\cdots
 z_n^{\pi_{_M}(C_{n,n})}\hspace{3cm}\end{eqnarray*}
\begin{eqnarray*}&=& \mbox{tr}_M\pi'
 (C_{i_1,i_2})\Phi \pi_{_M}(C_{i_2,i_1})z_1^{\pi_{_M}(C_{1,1})}z_2^{\pi_{_M}(C_{2,2})}\cdots
  z_n^{\pi_{_M}(C_{n,n})}\\ & &+\frac{z_{i_1}}{z_{i_2}}\mbox{tr}_M\Phi \pi_{_M}(C_{i_1,i_2}
  C_{i_2,i_1}+C_{i_2,i_2}-C_{i_1,i_1})z_1^{\pi_{_M}(C_{1,1})}z_2^{\pi_{_M}(C_{2,2})}\cdots
   z_n^{\pi_{_M}(C_{n,n})}\\ &=&
\mbox{tr}_M\pi'(C_{i_1,i_2})\Phi \pi_{_M}(C_{i_2,i_1})z_1^{\pi_{_M}(C_{1,1})}z_2^{\pi_{_M}
(C_{2,2})}\cdots  z_n^{\pi_{_M}(C_{n,n})}+\frac{z_{i_1}}{z_{i_2}}x^\ast {\cal E}^{i_1,i_2}\\
 & &+\frac{z_{i_1}}{z_{i_2}}(z_{i_2}\ptl_{z_{i_2}}-z_{i_1}\ptl_{z_{i_1}})x^\ast {\cal E}.
 \hspace{9.2cm}(7.17)\end{eqnarray*}
Thus
\begin{eqnarray*}x^\ast{\cal E}^{i_1,i_2}&=&\frac{z_{i_2}}{z_{i_2}-z_{i_1}}\mbox{tr}_M
\pi'(C_{i_1,i_2})\Phi
\pi_{_M}(C_{i_2,i_1})z_1^{\pi_{_M}(C_{1,1})}z_2^{\pi_{_M}(C_{2,2})}
\cdots  z_n^{\pi_{_M}(C_{n,n})}\\ &
&+x^\ast\frac{z_{i_1}}{z_{i_2}-z_{i_1}}(z_{i_2}
\ptl_{z_{i_2}}-z_{i_1}\ptl_{z_{i_1}}) {\cal
E}.\hspace{7cm}(7.18)\end{eqnarray*} Moreover,
\begin{eqnarray*}& &\mbox{tr}_M\pi'(C_{i_1,i_2})\Phi \pi_{_M}(C_{i_2,i_1})
z_1^{\pi_{_M}(C_{1,1})}z_2^{\pi_{_M}(C_{2,2})}\cdots
z_n^{\pi_{_M}(C_{n,n})}\\ &=&
\mbox{tr}_M\pi'(C_{i_1,i_2})\pi(C_{i_2,i_1})\Phi
z_1^{\pi_{_M}(C_{1,1})}z_2^{\pi_{_M} (C_{2,2})}\cdots
z_n^{\pi_{_M}(C_{n,n})}\\ &=&
 \mbox{tr}_M\pi'(C_{i_1,i_2})\pi'(C_{i_2,i_1})\Phi z_1^{\pi_{_M}(C_{1,1})}z_2^{\pi_{_M}
 (C_{2,2})}\cdots  z_n^{\pi_{_M}(C_{n,n})}\\ & &+ \mbox{tr}_M\pi'(C_{i_1,i_2})\pi_{_M}
 (C_{i_2,i_1})\Phi z_1^{\pi_{_M}(C_{1,1})}z_2^{\pi_{_M}(C_{2,2})}\cdots  z_n^{\pi_{_M}
 (C_{n,n})}\\ &=&x_{i_1}\ptl_{x_{i_2}}x_{i_2}\ptl_{x_{i_1}} (x^\ast {\cal E})+
 \mbox{tr}_M\pi'(C_{i_1,i_2})\Phi z_1^{\pi_{_M}(C_{1,1})}z_2^{\pi_{_M}(C_{2,2})}\cdots
 z_n^{\pi_{_M}(C_{n,n})}\pi_{_M}(C_{i_2,i_1})\\ &=&-\frac{1}{4}x^\ast {\cal E}+ \frac{z_{i_2}}
 {z_{i_1}}\mbox{tr}_M\pi'(C_{i_1,i_2})\Phi \pi_{_M}(C_{i_2,i_1}) z_1^{\pi_{_M}(C_{1,1})}
 z_2^{\pi_{_M}(C_{2,2})}\cdots  z_n^{\pi_{_M}(C_{n,n})}\hspace{2.1cm}(7.19)\end{eqnarray*}
by (4.30). Hence
$$\mbox{tr}_M\pi'(C_{i_1,i_2})\Phi \pi_{_M}(C_{i_2,i_1})
z_1^{\pi_{_M}(C_{1,1})}z_2^{\pi_{_M}(C_{2,2})}\cdots
z_n^{\pi_{_M}(C_{n,n})}=-
\frac{z_{i_1}}{4(z_{i_1}-z_{i_2})}x^\ast{\cal E}.\eqno(7.20)$$
Substituting it to (7.18), we obtain
$${\cal E}^{i_1,i_2}=\frac{z_{i_1}}{z_{i_2}-z_{i_1}}(z_{i_2}\ptl_{z_{i_2}}-z_{i_1}
\ptl_{z_{i_1}}) {\cal
E}+\frac{z_{i_1}z_{i_2}}{4(z_{i_1}-z_{i_2})^2}{\cal E}\qquad
\for\;\;i_1,i_2\in\ol{1,n},\;i_1\neq i_2.\eqno(7.21)$$

Let
$${\cal E}^{s,n+r}=(x^\ast)^{-1}\mbox{tr}_M\Phi \pi_{_M}( C_{s,n+r}C_{n+r,s})z_1^{\pi_{_M}
(C_{1,1})}z_2^{\pi_{_M}(C_{2,2})}\cdots
z_n^{\pi_{_M}(C_{n,n})}\eqno(7.22)$$ and
$${\cal E}^{n+r,s}=(x^\ast)^{-1}\mbox{tr}_M\Phi \pi_{_M}(C_{n+r,s}C_{s,n+r})z_1^{\pi_{_M}
(C_{1,1})}z_2^{\pi_{_M}(C_{2,2})}\cdots
z_n^{\pi_{_M}(C_{n,n})}\eqno(7.23)$$ for $1\leq s\leq r\leq n$. By
similar calculations as (7.17)-(7.21), we have
$${\cal E}^{s,n+r}=\frac{z_rz_s}{z_rz_s-1}(z_r\ptl_{z_r}+z_s\ptl_{z_s}){\cal E}+\frac{z_rz_s}
{4(z_rz_s-1)^2}{\cal E}.\eqno(7.24)$$
$${\cal E}^{n+r,s}=\frac{1}{z_rz_s-1}(z_r\ptl_{z_r}+z_s\ptl_{z_s}){\cal E}+\frac{z_rz_s}
{4(z_rz_s-1)^2}{\cal E}\eqno(7.25)$$
for $1\leq s<r\leq n$ and
$${\cal E}^{p,n+p}=\frac{z_p^2}{z_p^2-1}z_p\ptl_{z_p}{\cal E}+\frac{3z_p^2}{16(z_p^2-1)^2}
{\cal E}.\eqno(7.26)$$
$${\cal E}^{n+p,p}=\frac{1}{z_p^2-1}z_p\ptl_{z_p}{\cal E}+\frac{3z_p^2}{16(z_p^2-1)^2}
{\cal E}\eqno(7.27)$$
for $p\in\ol{1,n}$.

Now
\begin{eqnarray*}\mu \Psi&=&{\cal W}_{C_n}(x^\ast)^{-1}\mbox{tr}_M\Phi \Omega z_1^{\pi_{_M}
(C_{1,1})} z_2^{\pi_{_M}(C_{2,2})}\cdots  z_n^{\pi_{_M}{C_{n,n}}}
\\ &=&{\cal W}_{C_n}[\sum_{i_1,i_2=1}^n{\cal E}^{i_1,i_2}+\sum_{1\leq s\leq r\leq n}
(\dlt_{r,s}+1)({\cal E}^{n+r,s}+{\cal E}^{s,n+r})]\\ &=& {\cal
W}_{C_n}[\sum_{i=1}^n(z_i\ptl_{z_i})^2({\cal
E})+\sum_{i_1,i_2\in\ol{1,n}\:i_1\neq i_2}
\left(\frac{z_{i_1}}{z_{i_2}-z_{i_1}}(z_{i_2}\ptl_{z_{i_2}}-z_{i_1}\ptl_{z_{i_1}})
{\cal E} +\frac{z_{i_1}z_{i_2}}{4(z_{i_1}-z_{i_2})^2}{\cal
E}\right)\\ & &+\sum_{1\leq s<r\leq n}
\left(\frac{z_rz_s+1}{z_rz_s-1}(z_r\ptl_{z_r}+z_s\ptl_{z_s}){\cal
E}+\frac{z_rz_s}{2(z_rz_s-1)^2}{\cal E}\right)\\&
&+2\sum_{p=1}^n\left(\frac{z_p^2+1}{z_p^2-1}z_p\ptl_{z_p}{\cal
E}+\frac{3z_p^2}{8(z_p^2-1)^2}{\cal E}\right)]\\
&=&\left[\sum_{1\leq i_1<i_2\leq
n}\frac{1}{2}\left(\frac{z_{i_1}z_{i_2}}{(z_{i_1}-z_{i_2})^2}+
\frac{z_{i_1}z_{i_2}}{(z_{i_1}z_{i_2}-1)^2}\right)+\frac{3}{4}\sum_{p=1}^n
\frac{z_p^2}{(z_p^2-1)^2}\right]\Psi\\
& &+ 2\sum_{r,s\in\ol{1,n};\;r\neq
s}\left(\frac{z_r}{z_r-z_s}+\frac{1}{z_rz_s-1}\right)z_r\ptl_{z_r}{\cal
E}+2\sum_{p=1}^n\frac{z_p^2+1}{z_p^2-1}z_p\ptl_{z_p}{\cal E}\\ &
&+{\cal W}_{C_n}\sum_{i=1}^n(z_i\ptl_{z_i})^2({\cal E})\\&=&
\left[\sum_{1<i_1<i_2\leq
n}\frac{1}{2}\left(\frac{z_{i_1}z_{i_2}}{(z_{i_1}-z_{i_2})^2}+\frac{z_{i_1}z_{i_2}}
{(z_{i_1}z_{i_2}-1)^2}\right)+\frac{3}{4}\sum_{p=1}^n\frac{z_p^2}{(z_p^2-1)^2}\right]\Psi\\
& &+\sum_{i=1}^n[{\cal W}_{C_n}(z_i\ptl_{z_i})^2({\cal
E})+2z_i\ptl_{z_i}({\cal W}_{C_n})(z_i\ptl_{z_i})({\cal E})]\\ &=&
\left[\sum_{1<i_1<i_2\leq
n}\frac{1}{2}\left(\frac{z_{i_1}z_{i_2}}{(z_{i_1}-z_{i_2})^2}+\frac{z_{i_1}z_{i_2}}
{(z_{i_1}z_{i_2}-1)^2}\right)+\frac{3}{4}\sum_{p=1}^n\frac{z_p^2}{(z_p^2-1)^2}\right]\Psi\\
&
&+\sum_{i=1}^n(z_i\ptl_{z_i})^2(\Psi)-\sum_{i=1}^n(z_i\ptl_{z_i})^2({\cal
W}_{C_n}){\cal E}\\ &=&\left[\sum_{1<i_1<i_2\leq
n}\frac{1}{2}\left(\frac{z_{i_1}z_{i_2}}{(z_{i_1}-z_{i_2})^2}+\frac{z_{i_1}z_{i_2}}
{(z_{i_1}z_{i_2}-1)^2}\right)+\frac{3}{4}\sum_{p=1}^n\frac{z_p^2}{(z_p^2-1)^2}\right]\Psi\\
&
&+\sum_{i=1}^n(z_i\ptl_{z_i})^2(\Psi)-\frac{n(n+1)(2n+1)}{6}\psi\hspace{6.5cm}(7.28)
\end{eqnarray*}
by (6.25), (7.12)-(7.16), (7.21) and (7.24)-(7.27). Therefore, we
have: \psp

{\bf Theorem 7.1}. {\it The function $\Psi$ in (7.13) satisfies the following
Olshanesky-Perelomov equation of type C}:
\begin{eqnarray*}& &\sum_{i=1}^n(z_i\ptl_{z_i})^2(\Psi)+\left[\sum_{1<i_1<i_2\leq n}
\frac{1}{2}\left(\frac{z_{i_1}z_{i_2}}{(z_{i_1}-z_{i_2})^2}+\frac{z_{i_1}z_{i_2}}
{(z_{i_1}z_{i_2}-1)^2}\right)+\frac{3}{4}\sum_{p=1}^n\frac{z_p^2}{(z_p^2-1)^2}\right]\Psi\\
 &=&\left(\mu+\frac{n(n+1)(2n+1)}{6}\right)\Psi.\hspace{8.7cm}(7.29)\end{eqnarray*}
\pse

In particular, when $M=M_\lmd$ in (4.23),
$$\Psi={\cal W}_{C_n}E_C=
z_n^{\lmd_n+1}\left[\prod_{i=1}^{n-1}z_i^{n-i+1/2}\right]{\cal
X}_C(1/2,...,1/2;-\lmd_n)\{\xi^A_{r_2,r_1},\xi^C_{n+s_2,s_1}\}\eqno(7.30)$$
(cf. (4.115)) satisfies the above equation with
$$\mu= \sum_{i=0}^{n-1}\left(\lmd_n-\frac{i}{2}\right)\left(\lmd_n+\frac{3i}{2}+2\right).\eqno(7.31)$$
\vspace{1cm}

\noindent{\Large \bf References}

\hspace{0.5cm}

\begin{description}

\item[{[AAR]}] G. Andrews, R. Askey and R. Roy, {\it Special
Functions}, Cambridge University Press, 1999.

\item[{[BO]}] R. Beerends and E. Opdam, Certain hypergeometric series related to the root
 system $BC$, {\it Trans. Amer. Math. Soc. } {\bf 119} (1993), 581-609.

\item[{[BH]}] F. Beukers and G. Heckman, Monodromy for a hypergeometric function $_nF_{n-1}$,
{\it Invent. Math.} {\bf 95} (1989), no. 2, 325-354.

\item[{[BL]}]  D. Briten and F. Lemire,  A classification of simple Lie modules
having a 1-dimensional weight modules, {\it Trans. Amer. Math. Soc.} {\bf 299} (1987), no.2,
683-697.

\item[{[C]}] F. Calogero, Solution of the one-dimensional $n$-body problem with quadratic and
 /or inversely quadratic pair potentials, {\it J. Math. Phys.} {\bf 12} (1971), 419-432.

\item[{[Ep]}] P. Etingof, Quantum integrable systems and representations of
Lie algebras, {\it J. Math. Phys.} {\bf 36} (1995), no. 6,2637-2651.

\item[{[Eh]}] H. Exton, Multiple hypergeometric functions and applications, Halsted Pree
(John Wiley \& Sons, Inc,), New York, 1976.

\item[{[F]}] S. L. Fernando, Simple weight modules of complex reductive Lie algebras,
{\it Ph. D. Thesis, Univ. of Wisconsin at Madison}, 1983.

\item[{[GG]}] I. M. Gel'fand and M. I. Graev, GG-functions and their relation to general
hypergeometric functions, {\it Russian Math. Surveys} {\bf 52}(1997), no. 4, 639-684.

\item[{[Hg1]}] G. Heckman, Root systems and hypergeometric functions II, {\it Compositio. Math.}
 {\bf 64} (1987), 353-373.

\item[{[Hg2]}] G. Heckman, Heck algebras and hypergeometric functions, {\it Invent. Math.}
{\bf 100} (1990), 403-417.

\item[{[Hg3]}] G. Heckman, An elementary approach to the hypergeometric shift operators of
 Opdam, {\it Invent. Math.} {\bf 103} (1990), 341-350.

\item[{[HO]}] G. Heckman and E. Opdam, Root systems and hypergeometric functions I,
{\it Compositio. Math.} {\bf 64} (1987), 329-352.

\item[{[Hj]}] J. E. Humphreys, {\it Introduction to Lie Algebras and Representation Theory},
 Springer-Verlag New York Inc., 1972.

\item[{[J]}] J. C. Jantzen, Kontravariante formen auf induzierten Darstellungen habeinfacher
Lie-algebren, {\it Math. Ann.} {\bf 226} (1977), no. 1, 53-65.

\item[{[K]}] D. Kazhdan, B. Konstant and S. Steinberg, Hamiltonian group actions and dynamical
system of Calogero type, {\it Commun. Pure Appl. Math.} {\bf 31} (1978), 481-507.

\item[{[M]}] M. Marvan, Reducibility of zero curvature
representations with application to recursion operators, {\it
arXiv:nlin.SI/0306006}.

\item[{[NFP]}] J. F. N\'{u}nez, W. G. Fuertes and A. M. Perelonov,
Some results on the eigenfunctions of the quantum trigonometric Calogero-Sutherland model
related to the Lie algebra $D_4$,  {\it arXiv:math-ph/0305012}, 7 May 2003.

\item[{[O1]}]  E. Opdam, Root systems and hypergeometric functions III, {\it Compositio. Math.}
 {\bf 67} (1988), 21-49.

\item[{[O2]}]  E. Opdam, Root systems and hypergeometric functions IVI, {\it Compositio. Math.}
 {\bf 67} (1988), 191-209.

\item[{[O3]}] E. Opdam, Some applications of hypergeometric shift operators,
{\it Invent.  Math.} {\bf 98} (1989), 1-18.

\item[{[O4]}]  E. Opdam, An analogue of the Gauss summation formula for hypergeometric
 functions related to root systems, {\it Math. Z.} {\bf 212} (1993), 313-336.

\item[{[O5]}]  E. Opdam, Harmonic analysis for certain representations of graded Hecke algebras,
{\it Acta Math.} {\bf 175} (1995), 75-121.

\item[{[O6]}]  E. Opdam, Cuspital hypergeometric functions, {\it Methods Appl. Anal.} {\bf 6}
(1999), 67-80.

\item[{[OP]}] M. Olslanetsky and A. Perelomov, Completely integrable Hamiltonian systems
connected with semisimple Lie algebras, {\it Invent. Math.} {\bf 37} (1976), 93-108.

\item[{[S]}] B. Sutherland, Exact results for a quantum many-body problem in
one-dimension, {\it Phys. Rev. A} {\bf 5} (1972), 1372-1376.

\item[{[V1]}] D,-N. Verma, Structure of certain induced representations of complex semisimple
 Lie algebras, {\it thesis, Yale University,} 1966.

\item[{[V2]}] D,-N. Verma, Structure of certain induced representations of complex semisimple
 Lie algebras, {\it Bull. Amer. Math. Soc.} {\bf 74}(1968), 160-166.

\item[{[X]}] X. Xu, Differential equations for singular vectors of $sl(n)$, {\it Preprint}.

\end{description}

\end{document}